\documentclass[11pt,a4paper]{article}
\usepackage[utf8]{inputenc}
\usepackage[english]{babel}
\usepackage[T1]{fontenc}
\usepackage{amsmath}
\usepackage{amsthm}
\usepackage{amsfonts}
\usepackage{amssymb}
\usepackage{multicol}
\usepackage{comment}
\usepackage{graphicx}

\usepackage[all,cmtip]{xy} 

\usepackage[left=2cm,right=2cm,top=3cm,bottom=3cm]{geometry}
\usepackage{pgf,tikz}

\usetikzlibrary{arrows}
\usepackage{scalerel}[2016/12/29]
\newcommand\Hahn{\scaleobj{1.5}{\normalfont H}}
\DeclareMathOperator{\hprod}{\raisebox{-2.5pt}{\Hahn}}

\def\restriction#1#2{\mathchoice
              {\setbox1\hbox{${\displaystyle #1}_{\scriptstyle #2}$}
              \restrictionaux{#1}{#2}}
              {\setbox1\hbox{${\textstyle #1}_{\scriptstyle #2}$}
              \restrictionaux{#1}{#2}}
              {\setbox1\hbox{${\scriptstyle #1}_{\scriptscriptstyle #2}$}
              \restrictionaux{#1}{#2}}
              {\setbox1\hbox{${\scriptscriptstyle #1}_{\scriptscriptstyle #2}$}
              \restrictionaux{#1}{#2}}}
\def\restrictionaux#1#2{{#1\,\smash{\vrule height .8\ht1 depth .85\dp1}}_{\,#2}}

\newcommand*{\upwedge}{\mathbin{\raisebox{0ex}{\textasciicircum}}}

\DeclareMathOperator{\bdn}{bdn}

\DeclareMathOperator{\tp}{tp}
\DeclareMathOperator{\val}{val}
\DeclareMathOperator{\scal}{scal}

\DeclareMathOperator{\WD}{WD}

\DeclareMathOperator{\RV}{RV}

\DeclareMathOperator{\Th}{Th}

\DeclareMathOperator{\Card}{Card}

\DeclareMathOperator{\supp}{supp}
\DeclareMathOperator{\IP}{IP}
\DeclareMathOperator{\NIP}{NIP}

\author{Pierre Touchard\footnote{The author was partially supported through the “Oberwolfach Leibniz Fellows” program – Project-ID2043r, by Mathematisches Forschungsinstitut Oberwolfach in 2020 and supported by a grant of the University of Campania ‘Luigi Vanvitelli’ in the framework of V:ALERE 2019 (GoAL project)}}
\title{On Model Theory of Valued Vector Spaces\footnotetext{2020\textit{Mathematics Subject Classification}: Primary: 03C10
 03C45; Secondary: 03C60, 15A03}}

\theoremstyle{plain}
\newtheorem{theorem}{Theorem}[section]

\newtheorem{thmx}{Theorem}

\newtheorem{lemma}[theorem]{Lemma}
\newtheorem{proposition}[theorem]{Proposition}
\newtheorem{claim}{Claim}
\newtheorem{remark}[theorem]{Remark}

\newtheorem{fact}[theorem]{Fact}

\theoremstyle{definition}
\newtheorem{definition}[theorem]{Definition}
\newtheorem{corollary}[theorem]{Corollary}

\theoremstyle{remark}

\definecolor{black}{rgb}{0,0,0}
\newtheorem*{example}{Example}
\newtheorem*{examples}{Examples}

\newtheorem*{observation}{Observation}
\begin{document}

\maketitle
\begin{abstract}
    In analogy to valued fields, we study model-theoretic properties of valued vector spaces with variable base field by proving transfer principles down to the skeleton and down to the value set and base field. For instance, we give a formula which computes its burden in terms of the burden of its base field and its value set. To do this, we study these transfer principles in the context of lexicographic products of structures.
\end{abstract}
\section*{Introduction}
Relative quantifier elimination is a basic but powerful tool in model theory. We say that a theory $T$ eliminates quantifiers relative to a sort $\Sigma$ if it admits quantifier elimination in the language enriched with predicates for all definable sets in $\Sigma$.
Often, this allows us to prove relatives statement, which are typically of the form: `a model $M$ has  Property $P$ if and only if a certain algebraic condition is satisfied and $\Sigma(M)$ has Property $P$'.
Relative quantifier elimination and more generally relative statements are more flexible notions compared to absolute quantifier elimination and absolute statements for at least two distinct reasons. The first is that some natural classes of structures such as linear orders or Boolean algebras do not have quantifier elimination in a ‘reasonable’ language; if a (partial) theory interprets a sort ranging over such a class of structures then one can only hope to eliminates quantifiers relative to that sort. For instance, there is no reasonable language where all totally ordered abelian groups eliminate quantifiers because already for linear orders, such a language does not exist. However, Gurevich and Schmitt -- and later Clucker and Halupczok -- show  quantifier elimination relative to the spine -- that is to say, relative to the linear order of definable convex subgroups. Another reason is that in order to study a given structure, one can use relative quantifier elimination to proceed by steps: that is, one can first eliminate quantifiers relative to $\Sigma$ and then, after computing the induced structures, eliminates quantifiers in $\Sigma$. This approach can be very fruitful: this simplifies the search of a `reasonable language' and may also allow to eliminate quantifier in an effective way. One of the most important examples comes from valued field theory and the celebrated Theorem of Pas: it states that Henselian valued fields of equicaracteristic $0$ equipped with an angular component eliminates quantifiers relative to their value group and their residue field. However, in general no angular map is definable in the traditional language of valued fields; and naming such a map changes the structure \footnote{for instance, the ultraproduct $\prod_\mathcal{U}\mathbb{Q}_p$ has burden 1 in the language of valued fields, and $2$ in the language of Pas. See \cite{CS19} and \cite{Tou20a}).}. Only recently \cite[Theorem 5.15]{ACGZ20} a language was found where the theory of Henselian valued fields of equicharacteristic $0$ (with no named angular component) admits quantifier elimination relative to the residue field and value group. The authors of \cite{ACGZ20} use first a result of Flenner \cite{Fle11} -- which
eliminates quantifiers relative to the $\RV$-sort\footnote{Recall that the $\RV$-sort of a valued field $\mathcal{K}$ is the quotient of the multiplicative group $K^\star$ with the subgroup $1+\mathfrak{m}$, where $\mathfrak{m}$ denotes the maximal ideal of the valuation ring.}. Then with a clever use of $p.p.$-elimination of quantifiers in abelian structures, they describe a natural language where this $\RV$-sort eliminates quantifiers relative to the value group and the residue field.  In \cite{CS19} (resp. \cite{ACGZ20}), the authors reduce the study of dp-minimality, resp. distality, in Henselian valued fields of equicharacteristic $0$ to the $\RV$-sort, and then to the value group and residue field. For this second step, they observe that it is enough to see the $\RV$-sort as an enriched pure exact sequence of abelian groups. This two step reduction method has inspired many transfer principles in the context of (enriched) henselian valued fields (see for instance \cite{HM21},\cite{HK21},\cite{Tou20a},\cite{Tou20b}). In this paper, we study simpler objects: vector spaces and valued vector spaces (with value preserving scalar multiplication). Of course, model theory of two-sorted vector spaces is well studied, and the reader can refer to Pierce \cite{Pie09} and Kuzichev \cite{Kuz92} for original proofs of relative quantifier elimination. The theory of two-sorted vector spaces is one of the easiest example of relative quantifier elimination. We study transfer principles in such structures as an introduction and a warm up. Some of the model theory of vector spaces equipped with a valuation is already known, mainly due to \cite{KK97}, for that we cite as a reference on the subject. However, it seems that there is no result of relative quantifier elimination in valued vector spaces with variable base fields. We take the opportunity to address this and to study other transfer principles in valued vector spaces following a  method similar to \cite{CS19} and \cite{ACGZ20}.

In this study of valued vector space, we follow a similar strategy. This has allowed us to state transfer principles for quantifier elimination, model completness, NIPness and burden.  Although relative model completness and relative NIPness of a structure can often be deduced directly using relative quantifier elimination, the notion of burden however (and that of distality as in \cite{ACGZ20}) seems to justify the reduction method by step. By proving transfer principles for burden in valued vector spaces, we hope to show that the two-step reduction method can be adapted to contexts other than valued fields in order to study some subtle model theoretic properties.

We briefly describe the main results. Let $\mathrm{L}$ be the language of valued vector spaces over a non-specified field:
  \[\mathrm{L}:= \{(V,+,0, \scal : K \times V\rightarrow V),(K,+,\cdot, 0, 1),(\Gamma,<,\infty), \val:V \rightarrow \Gamma\}.\]
 A valued vector space $\mathcal{V}$ admits a natural expansion in the language
    \[\mathrm{L}_{W,\lambda}:= \mathrm{L} \cup \lbrace W_n \subset V^n\}_{n\in \mathbb{N}}\cup \{\lambda_{n,i}:V^{n+1} \rightarrow K\rbrace_{1\leq i \leq n \in \mathbb{N}},\]
where we interpret $W_n$ as the set of $K$-linearly value-independent $n$-tuples of vectors and $\lambda_{n,i}$ as the corresponding coordinate projections (see Section 3). We also enrich the value set with predicate symbols for the sets $D_n:=\{\gamma \in \Gamma \ \vert \ \dim_K(B_\gamma)=n\}$ where for $\gamma \in \Gamma$, $B_\gamma$ denotes the $K$-vector spaces $\{ v\in V \ \vert \ \val_{v}\geq \gamma \} / \{ v\in V \ \vert \ \val_{v}> \gamma \}$, called \textit{the rib of $\mathcal{V}$ at $\gamma$}  (see Section 3 for more details). Then we have:

\setcounter{thmx}{2}
 \begin{thmx}
        The theory of valued vector spaces eliminates quantifiers relative to the value set $(\Gamma,D_n,<)$, and to the base field $(K,+,\cdot,0,1)$ in the language $\mathrm{L}_{W,\lambda}\cup\{D_n \subset \Gamma\}$.
    \end{thmx}

 We deduce from Theorem \ref{TheoremQuantifierEliminationRelativeVariableBaseField2} some transfer principles. We emphasise here the one for burden:

{
    \renewcommand{\thetheorem}{\ref{PropositionValuedVectorSpaceBurdenTranfer}}
    \begin{proposition}
        Let $\mathcal{V}$ be a valued vector space with variable base field $K$ and value set $\Gamma$. 
        Then  \[\bdn(\mathcal{V}) \leq \max \{\kappa^{d}_{inp}(K), \bdn(\Gamma)\}\]
        if the $K$-dimension of $B_\gamma$ is bounded by $d$ for all $\gamma \in \Gamma$ and  
        \[\bdn(\mathcal{V})= \sup_n(\kappa^{n}_{inp}(K),\bdn(\Gamma))\]
        if the $K$-dimension of the $B_\gamma$'s are not bounded by any integer.

    \end{proposition}

    \addtocounter{theorem}{-1}
}

We prove these results in two steps. By analogy with valued fields, one may define  in a valued vector space $\mathcal{V}$ what is called \emph{the skeleton}   and reduce some model theoretic questions to it, and then reduce them to the base field and and value set. The skeleton $S(V)$ of $V$ can be defined as the quotient of $ D(V):=\{ (\gamma,v) \in \Gamma \times V \ \vert \ \val(v) \geq \gamma \}$ by the definable equivalence relation:
    \[(\gamma,v) \simeq (\gamma',v') \text{ if and only if } \gamma=\gamma' \text{ and } \val(v-v')>\gamma. \]
A similar definition can be found in \cite{Kuh}.

  We define the language $\mathrm{L}_S$ with a sort for the skeleton:
    \[\mathrm{L}_S:= \{(V,+,0),(\Gamma,<),(K,+,\cdot,0,1), (S(V),+_\bullet), \scal_\bullet, \val_S, \pi\},\]
       
    {\setcounter{thmx}{0}
  \begin{thmx}
        The theory of valued $K$-vector spaces eliminates quantifiers relative to the skeleton in the language $\mathrm{L}_S$.  
    \end{thmx}

    \setcounter{thmx}{0}
}

As we will see, this skeleton can be seen - at least to some extent - as an enriched (generalised) lexicographic product of the value set and vector spaces in the sense of Meir \cite{Mei19}. As noted by Hils, the lexicographic product can be seen as a non-algebraic equivalent of a short pure exact sequence of abelian groups. The study of valued vector spaces and their skeletons is thus an opportunity to push the analogy forwards, which was originally the motivation for the writing of these notes. There are two major differences between valued fields and valued vector spaces:
while a main challenge in proving relative quantifier elimination of valued fields is to linearise terms occuring in formulas (see Flenner's cell decomposition in $\cite{Fle11}$), in vector spaces however, all terms are of course linear. Elimination of quantifiers relative to the skeleton admits, for this reason, a very simple proof. However, the lack of multiplicative structure brings a new kind of difficulty: their is no notion of residue field.  Then, in order to state transfer principles, one has to consider instead $\vert \Gamma \vert$-many structures, namely the ribs of $\mathcal{V}$. One needs then some caution to prove a transfer down to the value set ( also called `\textit{spine}') and the ribs . Addressing these difficulties may not be very challenging, but this approach using the skeleton and lexicographic products should bring new perspectives to prove other transfer principles in valued vector spaces and also in other contexts such as valued groups and valued modules. With this end in mind, we discuss in detail transfer principles in lexicographic products. We notably prove the following:

{
    \renewcommand{\thetheorem}{\ref{TheoremRelativeQuantifierEliminationEnrichedLexProd}}
\begin{theorem}
 Consider the lexicographic sum $\mathcal{S}:=\mathcal{M}[\mathfrak{N}]$ of a class of $\mathrm{L}_\mathfrak{N}$-structures $\mathfrak{N}:= \{\mathcal{N}_a\}_{a\in \mathcal{M}}$ with respect to an $\mathrm{L}_\mathcal{M}$-structure $\mathcal{M}$ in the language
    \[\mathrm{L}_{\mathcal{M}[\mathfrak{N}]}:= (S,\mathrm{L}_{\bullet,\mathfrak{N}}) \cup (\mathcal{M},\mathrm{L}_\mathcal{M})\cup \{v: S\rightarrow \mathcal{M}   \}.\]
    
    Assume that for all sentences $\phi\in \mathrm{L}_\mathfrak{N}$, the set $\{a\in \mathcal{M} \ \vert \ \mathcal{N}_a \models \phi \}$ is $\emptyset$-definable in $\mathcal{M}$.
    Assume also that $\mathfrak{N}$ interprets uniformly a structure $\mathcal{K}$ (in a language $\mathrm{L}_K$) with a non-empty  set of \emph{surjective} maps $\{\lambda_i: \mathcal{N}_a^{k_i} \rightarrow K\}_{a\in \mathcal{M},i\in I}$.
    We can enriched $\mathcal{S}$ with a sort for $\mathcal{K}$ and surjective maps $\lambda_{\bullet,i}: {S}^k \rightarrow K$ define as follows:
    \[ \forall (a_1,b_1), \dots, (a_k,b_k) \in S, \ \lambda_{\bullet,i}((a_1,b_1), \dots, (a_k,b_k)): = \begin{cases} \lambda_i^{\mathcal{N}_a}(b_1,\dots,b_k) & \text{ if } a_k=\cdots=a_k=a \\
    0 & \text{ otherwise.}\end{cases}\]
    where $0$ is a constant in $K$. 
    Finally, assume that:
        \begin{itemize}
            \item [$(RQE)_\mathfrak{N}$] For all $\rho\in S_1^{\Th(\mathcal{M})}$, $\mathcal{N}_\rho$ admits quantifier elimination relative to $K$ in $\mathrm{L}_{\mathfrak{N}} \cup \mathrm{L}_K \cup \{\lambda_i \}_{i\in I}$.

        \end{itemize}
        Then $(\mathcal{M}[\mathfrak{N}],\mathcal{K},\{\lambda_{\bullet,i}\}_{i\in I})$ eliminates quantifiers relative to $\mathcal{M}$ and $K$ in $\mathrm{L}_{\mathcal{M}[\mathfrak{N}]}\cup \mathrm{L}_K \cup \{\lambda_{\bullet,i}\}_{i\in I}$.   
\end{theorem}

    \addtocounter{theorem}{-1}
}


As an introduction and a warm up, we will give in Section 1 a new proof of relative quantifier elimination in two-sorted vector spaces and some of its consequences: vector spaces are model complete (resp. NIP) relative to the base field. We also prove a transfer principle for burden. In Section 2, we study lexicographic products $\mathcal{M}[\mathfrak{N}]$: we recall the definition and transfer principles obtained by Meir for quantifier elimination, model completness and NIPness \cite{Mei16}. We extend these results to a slightly more general setting, notably by proving Theorem \ref{TheoremRelativeQuantifierEliminationEnrichedLexProd}. Then we compute the burden in terms of the burden of the structure $\mathcal{M}$ and that of the structures $\mathcal{N}_a$, $a\in \mathcal{M}$. 
In Section 3, we study the valued vector spaces with value-preserving scalar multiplication. We prove first Theorem \ref{TheoremQuantifierEliminationRelativeVariableBaseField2}. Then we prove other transfer principles in valued vector spaces down to the base field and value group model completness, NIPness and burden. If the reader is only interested in the first three, they can skip the formalism of Section 2.


\setcounter{thmx}{0}

 As we said, the reader aware of the work in \cite{CS19} and \cite{ACGZ20} will be familiar with the method used in this text. To them, we give here a table to hopefully clarify all the analogies.\\
 
 \vspace{1.0 cm}
 \hspace{-2.3 cm}
 {        \small
    \begin{tabular}{ | c | c |}
    \hline
        Valued $K$-Vector Space $V$ of Value set $\Gamma$ & Henselian Valued Fields $K$ of value group $\Gamma$ and residue $k$.\\
        \hline \hline
        \begin{tabular}{c}
             Value set with induced structure\\
             $(\Gamma, (D_n)_{n\in \mathbb{N}}, < )$ where   \\
              $D_n:=\{\gamma\in \Gamma \ \vert \ \dim_K(B_\gamma)=n \}$
              \vspace{0.2cm}
        \end{tabular}  & 
        \begin{tabular}{c}
        Value group with induced structure \\
        $(\Gamma, + , 0 , < )$     
        \end{tabular} \\
         \hline
          \begin{tabular}{c}

          The ribs with  induced structure (K-vector space) \\ 
       $(B_\gamma, 0, +, \scal: K \times B_\gamma \rightarrow B_\gamma  ), \qquad \gamma \in \Gamma$ \\
       where $B_\gamma:= \{x\in V; \val(x)\geq \gamma\} / \{x\in V; \val(x)> \gamma\}$
       \vspace{0.2cm}
       \end{tabular}
       &
        \begin{tabular}{c}
        The residue field with induced structure (field)\\
       $(k,+,\cdot,1,0)$ \\
       where  $k \simeq B_{\geq \gamma}(0)/ B_{>\gamma}(0)$ for any $\gamma \in \Gamma$.
       \end{tabular}
       \\
               \hline
               \begin{tabular}{c}
               Skeleton \\
               $S(V) := \{ (\gamma,v) \in \Gamma \times V \ \vert \ \val(v) \geq \gamma \} / \simeq$ \\
               where $(\gamma,v) \simeq (\gamma',v') \text{ if and only if } \gamma=\gamma' \text{ and } \val(v-v')>\gamma.$
               \vspace{0.2cm}
               \end{tabular}
               &
               \begin{tabular}{c}
               $\RV$-sort\\
               $\RV(K):= K^\star/1+\mathfrak{m}$
               \end{tabular}
               \\
       \hline
       \begin{tabular}{c}
            Skeleton with  induced structure\\
            $S(V) \simeq \{\Gamma[B_\gamma]_{\gamma \in \Gamma},(K,+,\cdot,0,1), (\Gamma,<),W_{\bullet,n},\lambda_{\bullet,n,i}\}$
            \vspace{0.2cm} 
        \end{tabular}
        
            & 
        \begin{tabular}{c}
            $\RV$-sort with induced structure\\
            $1 \rightarrow k^\star \rightarrow \RV^\star \rightarrow \Gamma \rightarrow 0$  
        \end{tabular}\\
       \hline
            \begin{tabular}{c}
                Relative quantifier elimination (Theorems \ref{TheoremQuantifierEliminationRelativetoSkeleton}  \ref{TheoremQuantifierEliminationRelativeVariableBaseField} and \ref{TheoremQuantifierEliminationRelativeVariableBaseField2})\\
                
                \begin{tikzpicture}
                    \node {$(V,K, \Gamma)$ }
                        child{ 
                        node {$S(V)=\Gamma[B_\gamma]$}
                            child { 
                            node { $\Gamma$ } 
                            }
                            child { node { $B_{\gamma} $ }
                                child [missing] 
                                child [missing] 
                                child { node {} }
                            }
                            child { node { $B_{\gamma'}$ } 
                                child { node {$K$} }
                            }
                            child { node { $B_{\gamma''}$ } 
                                child { node {
                                } }
                                    child [missing]
                                    child [missing] 
                            }
                            child { node {$\ldots$}
                                    child { node {
                                } }
                                    child [missing]
                                    child [missing]
                                    child [missing]
                                    child [missing]
                            }
                        };
                \end{tikzpicture}
            \end{tabular}
  &
            \begin{tabular}{c}
                Relative quantifier elimination\\
                \begin{tikzpicture}
                    \node {$(K,\Gamma,k)$ }
                        child { node { $\RV(K)$ } 
                            child { node {$\Gamma$} }
                            child { node {$k$} 
                                child [missing]
                                }
                        };
                \end{tikzpicture}
            \end{tabular}
            \\
            \hline
       \end{tabular}
}
\newpage

\section*{Notation and prerequisites}

We use standard notation from model theory. One may refer to \cite{TZ12} for the most basic notions, such as \textit{quantifier elimination},\textit{ model completeness}, \textit{bi-interpretability}, \textit{stable embeddedness} etc. We refer the reader to \cite{Sim15} for NIP theories and related notions and to \cite{Che14} for \textit{burden}, \textit{inp-patterns} and related notions. Occasionaly, we use Adler's convention  \footnote{One can consider `new cardinals' of the form $\lambda_-$ where $\lambda$ is a limit cardinal. Some notions of dimension, such as the burden, extend naturally to this extended set of cardinals, denoted by $\Card^\star$.}\cite{Adl07}. When we do, this will be indicated to the reader with a superscript $\star$: $\Card^\star$, $\bdn^\star$, $\sup^\star$ etc. 
 For more details, notably about the arithmetic of $\Card^\star$, see \cite[Section 1.1.3]{Tou20a}). We will also use Rideau-Kikushi's terminology in \cite[Annex A]{Rid17}. The reader can refer to this annex for the definition of \textit{relative quantifier elimination}, \textit{closed sorts}, \textit{enrichment} and \textit{resplendence}.

A structure is typically denoted by a curvy letter or a capital Greek letter: $\mathcal{V},\mathcal{K}$,$\Gamma ...$ The corresponding base set is then denoted by $V, K ,\Gamma ...$ If a tuple of elements in a set $S$ is denoted by a single letter $a$, them $a_1,\dots,a_{\vert a\vert}$ denote its coordinates. If $S$ is a sort in a language $\mathrm{L}$, we mean by $S$-variable a variable in the sort $S$ and by  \textit{$S$-sorted-quantifier} a quantifier ranging over the sort $S$; that is $\forall x\in S$ or $\exists x\in S$.
By \textit{$S$-sorted-quantifier-free $\mathrm{L}$-formula}, or $S$-quantifier-free $\mathrm{L}$-formula, so we mean an $\mathrm{L}$-formula with no $S$-sorted-quantifiers. 

\section{Vector Spaces}

We discuss here some model theory of vector spaces. Most of the results in this section are folklore, but we haven't found reference for all of them. We take the occasion to write them down with a modern presentation (relative statement). These results will be used in Section 3 when we will discuss valued vector spaces. The reader can refer to the article of Pierce \cite{Pie09} or that of Kuzichev \cite{Kuz92}. One will also find some development in \cite{Gra99}. 

We cite here a theorem of Pierce, which describe a minimal language for the structure of vector spaces with variable base field:
\begin{fact}[{\cite[Theorem 1.1]{Pie09}}]\label{FactPierce}
    A $K$-vector space $\mathcal{V}$ interprets the field $K$ in the language $\{V,0,+,W_2\}$, where $W_2$\footnote{the $\neg W_2$ is called parallelism relation and denoted by $\vert\vert$ in \cite{Pie09}.} is interpreted as the relation of $2$-independence.
\end{fact}

We focus on the opposite thematic: we will describe an interpretable language rich enough to admit quantifier elimination relative to the base field. We then prove other transfer principles in vector spaces down to the base fields.   

\subsection{Relative Quantifier Elimination in Vector Spaces with variable base fields. }
We describe a well-known result of relative quantifier elimination of vector spaces in the two-sorted language. We want to emphasise here a simple observation: there is a language where vector spaces eliminate quantifiers relative to the base field, and such that the field sort is closed. Such language has then good syntactical properties, with which it is natural to work in order to prove transfer principles.


We consider the theory of vector spaces $\mathcal{V}$ over a non-fixed field $K$ in the 2-sorted language $\mathrm{L}=\lbrace V,+,0 \} \cup \{ K,+,\cdot,0,1\} \cup \{\scal: K\times V  \rightarrow V\}$, where $\scal(s,v)$ is interpreted by the scalar multiplication $s \cdot v$.  To obtain relative quantifier elimination, we have to name: 
\begin{itemize}
    \item for all $n \in \mathbb{N}^\star$, the $n$-ary predicate $W_n \subset V^n$ for $K$-linear independent $n$-tuples of vectors.
    \item for $1\leq i\leq n  \in \mathbb{N}^\star$, the $i^{\text{th}}$ \textit{lambda function of order $n$}, denoted by $\lambda_{n,i}: V^{n+1}\rightarrow K$:
\[ \lambda_{n,i}(v,v_1, \cdots, v_n):= \begin{cases}
\lambda_i & \text{ if } v=\sum_{1\leq
i\leq n}\lambda_i v_i \text{ and } W_n(v_1,\cdots,v_n),\\
0 &\text{ if } \neg W_n(v_1,\cdots,v_n)\vee W_{n+1}(v,v_1,\cdots,v_n).
\end{cases}\]

\end{itemize}

Notice that the graph of the lambda-functions and the predicates $W_n$ are defined with quantifiers in $K$ only.  
We denote by $\mathrm{L}_{W,\lambda}$ the language
\[\mathrm{L} \cup \lbrace W_n \subset V^n\}_{n\in \mathbb{N}}\cup \{\lambda_{n,i}:V^{n+1} \rightarrow K\rbrace_{1\leq i \leq n \in \mathbb{N}}.\]

We remark that in the language $\mathrm{L}_{W,\lambda}$, the symbols for scalar multiplication $\scal$ and for the addition $+$ in $\mathcal{V}$ are no more required to express definable sets:
\begin{remark}\label{RemarkBiinterpretabilityLWLambda}
    The structures
    $\mathcal{V}=\lbrace V,+,0 \} \cup \{ K,+,\cdot,0,1\} \cup \{\scal: K\times V  \rightarrow V\}$
    
    and
    
    \[\{ V\} \cup \{ K,+,\cdot,0,1\} \cup \lbrace W_n \subset V^n\}_{n\in \mathbb{N}^\star}\cup \{\lambda_{n,i}:V^{n+1} \rightarrow k\}_{1\leq i \leq n\in \mathbb{N}^\star}\]
    are bi-interpretable.
\end{remark}

This follows from basic linear algebra. If $\bar{x}$ and $\bar{y}$ are respectively tuple of variables in $V$ and $K$, then one has:
\begin{align}
\sum_i y_i \cdot x_i =0 \ \Leftrightarrow \ \bigvee_{k\leq n } \bigvee_{\sigma \in  \mathfrak{S}_n}  & W_k(x_{\sigma(1)},\dots, x_{\sigma(k)}) \wedge \bigwedge_{j>k}\neg W_k(x_{\sigma(j)},x_{\sigma(1)},\dots, x_{\sigma(k)}) \label{Equationbiinterpretabilityproof}\\
\tag*{} & \wedge \bigwedge_{i\leq k} y_{\sigma(i)}+  \sum_{n\geq j>k} \lambda_{k,i}(x_{\sigma(j)},x_{\sigma(1)},\dots,x_{\sigma(k)})y_{\sigma(j)} =0. 
\end{align}.

In general, if $x$ and $\bar{x}'$ are $V$-variables and $\phi(x,\bar{x}')$ an $\mathrm{L}_{W,\lambda}$-formula. Then:

\begin{align}
\phi(\sum_i y_i \cdot x_i, \bar{x}')  \ \Leftrightarrow \ & \bigvee_{k\leq n } \bigvee_{\sigma \in  \mathfrak{S}_n}   W_k(x_{\sigma(1)},\dots, x_{\sigma(k)}) \wedge \bigwedge_{j>k}\neg W_k(x_{\sigma(j)},x_{\sigma(1)},\dots, x_{\sigma(k)}) \label{EquationWMALinIndependent} \\
\tag*{} & \wedge \phi\left(\sum_i \left(y_{\sigma(i)}+  \sum_{n\geq j>k} \lambda_{k,i}(x_{\sigma(j)},x_{\sigma(1)},\dots,x_{\sigma(k)})y_{\sigma(j)}\right) x_{\sigma(j)} , \bar{x}'\right).
\end{align}.

In particular, we have :

\begin{align}
W_{\vert x'\vert+1}(\sum_i y_i \cdot x_i, \bar{x}')  \ \Leftrightarrow \ & \bigvee_{k\leq n } \bigvee_{\sigma \in  \mathfrak{S}_n}   W_k(x_{\sigma(1)},\dots, x_{\sigma(k)}) \wedge \bigwedge_{j>k}\neg W_k(x_{\sigma(j)},x_{\sigma(1)},\dots, x_{\sigma(k)}) \label{EquationW} \\
\tag*{} & \wedge \bigvee_{i\leq k} W_{\vert x'\vert+1}(x_{\sigma(j)} , \bar{x}') \wedge y_{\sigma(i)}+  \sum_{n\geq j>k} \lambda_{k,i}(x_{\sigma(j)},x_{\sigma(1)},\dots,x_{\sigma(k)})y_{\sigma(j)} \neq 0.
\end{align}
and similarly we can express lambda-terms $\lambda(\sum_i y_i \cdot x_i, \bar{x}')$ in terms of $\lambda( x_i, \bar{x}')$ using (\ref{EquationWMALinIndependent}) and the fact that if the $x_i$'s are linearly independent element of the sub-vector space generated by $\bar{x}'$, then 

\begin{align}
\lambda(\sum_i y_i \cdot x_i, \bar{x}')= \sum_i  y_i \cdot \lambda( x_i, \bar{x}') \label{EquationLambda}.
\end{align}

Let us denote by $\mathrm{L}_{W,\lambda}'$ the language $\mathrm{L}_{W,\lambda}$ without the addition and the scalar multiplication: \[\mathrm{L}_{W,\lambda}'= \{ V \} \cup \{ K,+,\cdot,0,1\} \cup \lbrace W_n\}_{n}\cup \{\lambda_{n,i}  \rbrace_{1\leq i \leq n}.\]

\begin{theorem}\label{TheoremRelativeQE}
Let $\mathcal{V}$ be a vector space.
\begin{itemize}
    \item Then $\mathcal{V}$ eliminates quantifiers relative to $K$  in the languages $\mathrm{L}_{W,\lambda}'$ and $\mathrm{L}_{W,\lambda}$, 
    \item $\mathcal{V}$ eliminates $V$-sorted quantifiers in the language $\mathrm{L}$, and
    \item the theory of $\mathcal{V}$ is determinate by the theory of $\mathcal{K}$ and its dimension.
\end{itemize}
\end{theorem}

By the remark above, one deduces easily the second statement from the first: every formula $\phi(\bar{x}, \bar{y})$ in the language $\mathrm{L}$ where $\bar{x}$ is a tuple of $V$-variables and $\bar{y}$ is a tuple of $K$-variables is equivalent to a Boolean combination  of formulas of the form:

\[\phi_K(\lambda(\bar{x}),\bar{y})\]
and 

\[W_{n}(x_{\sigma(1)},\cdots, x_{\sigma(n)}),\]

where $\phi_k(\bar{z})$ is a formula in the language of fields and where $\lambda(\bar{x})$ stands for finitely many terms of the form 
\[\lambda_{i,n}(x_{\sigma(1)}, \cdots,x_{\sigma(n+1)})\]
for some integers $i\leq n$ and some functions $\sigma:\mathbb{N} \rightarrow \mathbb{N}$. Then, one has to observe that any formula of the form
$\phi_K(\lambda_{i,n}(x,x_1,\cdots,x_n),\bar{y})$ is equivalent to:
\begin{align*}
\left ( \exists s_1,\dots, s_n\in K \ \phi_K(s_i, \bar{y}) \wedge\ \sum_{1\leq i\leq n} s_i\cdot x_i =x \wedge W(x_1,\cdots,x_n) \right ) \\
\vee  \left (  W_{n+1}(x,x_1,\dots,x_n) \vee  \neg W_n(x_1,\cdots,x_n) \right ) \wedge
\phi_K(0, \bar{y}). 
\end{align*}
and of course that $W_n(x_1,\cdots,x_n)$ is equivalent to
\[\ \forall s_1, \cdots, s_{n} \in K, \ \sum_{1\leq i\leq n} s_i\cdot x_i =0 \ \Leftrightarrow s_1=\cdots=s_n=0.\]
We obtained an $\mathrm{L}$-formula free from any $V$-sorted-quantifiers. Before proving relative quantifier elimination in the language $\mathrm{L}_{W,\lambda}'$, let us state a direct corollary:

\begin{corollary}\label{Corollarykpure}
The field $K$ is stably embedded and the induced structure on $K$ is the pure structure of fields. Any finite-dimensional vector subspaces $V'$ is stably embedded, and its induced structure is the pure structure of vector space.
\end{corollary}

\begin{proof}
This is due to the fact that in the language $\mathrm{L}_{W,\lambda}'$, the sort $K$ is closed: by relative quantifier elimination, a definable set in $K^{\vert x_K \vert}$ with parameters $a\in K$ and $v \in V$ is given by a formula $\phi(x_K, v,a)$ of the form $\phi_K(x_K, t(v), a)$ where $\phi_K(x_K,y_K,z_K)$ is a formula in the language of rings and $t(v)$ is a tuple of terms in $\mathcal{V}$. Then we simply have
\[\phi_K(K^{\vert x_K \vert}, t(v), a) = \phi(K^{\vert x_K \vert}, v,a).\]

We deduce the statement on finite-dimensional vector subspaces using the identification $K^{\dim_K(V')} \simeq V'$ with the choice of a basis. 
\end{proof}

\begin{proof}[Proof of the theorem]
We prove quantifier elimination by back and forth, in the language $\mathrm{L}_{W,\lambda}$. Let $\mathcal{M}=(V,K)$ and $\mathcal{N}=(V',K')$ two models.
Let $f=(f_{V_A},f_{K_A}):(V_A,K_A)\rightarrow (V_B,K_B)$ be an isomorphism between a finitely generated substructure $A=(V_A,K_A)$ of $\mathcal{M}$ and a finitely generated substructure $B=(V_B,K_B)$ of $\mathcal{N}$, and assume that $f_{K_A}$ is elementary. Notice that $A$ and $B$ are closed under the $\lambda$-functions. We want to extend $f$ to a global embedding $\tilde{f}:\mathcal{M} \rightarrow \mathcal{N}$.\\
\textbf{Step 1}: we extend $f_K:K_A \rightarrow K_B$ to $K$ as a field embedding.  As $f_{K_A}$ is elementary and $K'$ is $\vert K\vert^+$-saturated, we know that such function can be extended to a field-embedding $f_K:K\rightarrow K'$.
Consider $W$ the $K$-vector space generated by $V_A$. Naturally, we extend $f$ to $W$ by setting for $s_1,\dots,s_n \in K$ and $v_1,\dots,v_n \in V_A$ linearly independant:
\[f(\sum_{i\leq n} s_i \cdot v_i)= \sum_{i\leq n} f_K(s_i) \cdot f(v_i).\]
The fact that it is indeed a partial isomorphism is guarantied by  
(\ref{Equationbiinterpretabilityproof}), (\ref{EquationW}) and (\ref{EquationLambda}).\\

\textbf{Step 2}: By Step 1, we may assume that $A=(V_A,K)$ is a $K$-subvector space. If $V_A=V$, there is nothing left to do. Assume there is $w\in V \setminus V_A$. Then $w$ is $K$-linearly independent from $V_A$. We find  $w'\in V' \setminus V_B$ (using  the fact that $\dim_{K}(V)>\dim_{K}(V_A)=\dim_{K'}(V_B)<\dim_{K'}(V')$ ). We can extend $f$ to the $K$-vector space $W$ generated by $V_A$ as follows:
\[\begin{array}{ccccc}
f_W & : & V_A \oplus K\cdot w & \to & V_B \oplus K'\cdot w' \\
 & & a+s\cdot w & \mapsto & f_V(a)+f_K(s)\cdot w'. 
\end{array}\]

Using again (\ref{Equationbiinterpretabilityproof}), (\ref{EquationW}) and (\ref{EquationLambda}), one sees that it defines a partial isomorphism. 
We can of course repeat step $2$ until we get an embedding of $V$ in $V'$. This concludes our proof of relative quantifier elimination in the language $\mathrm{L}_{W,\lambda}$.
We deduce relative quantifier elimination in the language $\mathrm{L}_{W,\lambda}'$ using (\ref{Equationbiinterpretabilityproof}), (\ref{EquationW}) and (\ref{EquationLambda}) again.
\end{proof}

\subsection{Transfer principles in vector spaces relative to the base field}
We start with some easy corollaries:
\begin{corollary}
    Let $\mathcal{V}=(V,K)$ be a vector space. Then $\mathcal{V}$ is model complete relative to $K$ in the language $\mathrm{L}_{W,\lambda}'$. More precisely: assume that $\mathcal{W}=(W,L)$ is an extension of $\mathcal{V}$ in $\mathrm{L}_{W,\lambda}'$, then $\mathcal{V} \preceq \mathcal{W}$ if and only if $K \preceq L$ in the language of rings.
\end{corollary}
This is also completely folklore. We include a proof for completeness
\begin{proof}
    Let $\mathcal{W}=(W,L)$ be an extension of $(V,K)$ and assume that $L \succeq K$ as fields. We apply the Tarsky-Vaught criterion: consider any $\mathrm{L}_{W,\lambda}'(\mathcal{V})$-formula $\phi(x,\bar{v},\bar{a})$ with parameters $\bar{v} \in V$ and $\bar{a}\in K$  and assume that $\mathcal{W} \models  \phi(w, \bar{v},\bar{a})$ for some $w\in W$. We have to show that $\mathcal{V} \models  \ \phi(v,\bar{v},\bar{a})$ for some $v\in V$. By quantifier elimination,  the formula $\exists x \ \phi(x,\bar{x_V},\bar{x_K})$ is equivalent to a disjunction of formulas of the form:
    \[ \phi_K( t_0(\bar{x_V}), \dots , t_{k-1}(\bar{x_V}), \bar{x_K} ) \wedge \phi_W(\bar{x_V})\]
where $\phi_K(x_{K,0},\dots,x_{K,k-1} , \bar{y_K})$ is a field formula, $t_i(x,\bar{x_V})$ are $\lambda$-terms and $\phi_W(\bar{x_V})$ is a conjunction of $W_n(x,\bar{x_V}')$ for subtuples $\bar{x_V}'$ of $\bar{x_V}$ and negation of such formulas. As $K \preceq L$, $L \models \phi_K( t_0(\bar{v}), \dots , t_{k-1}(\bar{v}), \bar{a} )$ if and only if $K \models \phi_K( t_0(\bar{v}), \dots , t_{k-1}(\bar{v}), \bar{a} )$. As $\phi_W(\bar{v})$ is quantifier-free, we get that 
\[\mathcal{W} \models \phi_K( t_0(\bar{v}), \dots , t_{k-1}(\bar{v}), \bar{a} ) \wedge \phi_W(\bar{v}) \leftrightarrow \mathcal{V} \models \phi_K( t_0(\bar{v}), \dots , t_{k-1}(\bar{v}), \bar{a} ) \wedge \phi_W(\bar{v}).\]
This shows that $\mathcal{V} \models \exists x \ \phi(x,\bar{v},\bar{a}) $ and concludes the proof.
\end{proof}

    This statement and its proof can be easily generalised: 
    \begin{remark}\label{RemarkModelCompletnessRQECLoseSort}
        Let $\mathrm{L}$ be a language, $\Sigma \cup  \Pi$ be a partition of sorts such that $\Sigma$ is closed.  Assume $\mathcal{M}$ is an $\mathrm{L}$-structure which admits quantifier elimination relative to $\Sigma$. Then $\mathcal{M}$ is model complete relative to $\Sigma(M)$:
        if $\mathcal{M} \subseteq \mathcal{N}$ is an extension in the language $\mathrm{L}$, then $\mathcal{M} \preceq \mathcal{N}$ if and only if $\Sigma(\mathcal{M}) \preceq \Sigma(\mathcal{N})$ in $\restriction{\mathrm{L}}{\Sigma}$ the language restricted to $\Sigma$.
    \end{remark}

We state now a transfer principle for the independent property (see \cite[Definition 2.10]{Sim11}). 

\begin{corollary}\label{CorollaryNIPtransferVectorSpace}
    Let $\mathcal{V}=(V,K)$ be a vector space in the two sorted language $\mathrm{L}$ (or $\mathrm{L}_{W,\lambda}$ or any other interpretable language). Then $\mathcal{V}$ is NIP if and only $K$ is NIP.
\end{corollary}

\begin{proof}
    If $\mathcal{V}$ is NIP, then so is any structure interpretable in $\mathcal{V}$. In particular $K$ is NIP.
    Assume reciprocally that $K$ is NIP. By quantifier elimination, it is enough to check that every atomic formula $\phi(x_V;y_V,y_K)$ where $x_V,y_V$ are variables from $V$ and $y_K$ are variables from $K$ with $\vert x_V \vert =1$  are NIP. In the language $\mathrm{L}_{W,\lambda}$, these are:
    \begin{itemize}
        \item $\phi_K( t_0(x_V,y_V), \dots , t_{k-1}(x_V,y_V),y_K)$ where $\phi_K(x_{K,0},\dots,x_{K,k-1} , y_K)$ is an atomic formula in the language of fields. 
        \item $W_n(x_V,y_V')$ for a subtuple $y_V'$ of $y_V$.
    \end{itemize}
    The first are NIP since the field $K$ is NIP. Clearly, the second is NIP of VC-dimension exactly $n$. 
\end{proof}

Again, we can make a more general statement:
    \begin{remark}\label{RemarkNIPRQECLoseSort}
        Let $\mathrm{L}$ be a language, $\Sigma \cup  \Pi$ be a partition of sorts such that $\Sigma$ is closed.  Assume $\mathcal{M}$ is an $\mathrm{L}$-structure which admits quantifier elimination relative to $\Sigma$. Assume that atomic formulas in $\restriction{\mathrm{L}}{\Pi}$, the language restricted to $\Pi$, are $\NIP$. Then $\mathcal{M}$ is $\NIP$ if and only if $\Sigma(M)$ is NIP.
    \end{remark}

We state now another result of classification regarding the burden (see \cite[Definitions 2.1 \& 2.8]{Che14}).  

\begin{theorem}\label{TheoremBurdenVectorSpace}
    Let $V$ be a $K$-vector space, where $K$ is an infinite field. 
    \begin{itemize}
        \item Assume $V$ is of dimension $d$. Then $\bdn(V,K)= \kappa_{inp}^d(K):= \bdn(K^d, \pi_j:K^d \rightarrow K, j\leq d)$.
        \item Assume $V$ is of infinite dimension. Then $\bdn(V,K)= \sup_{d\in \mathbb{N}}\kappa_{inp}^d(K)$.
    \end{itemize}
    
\end{theorem}
We may more precisely use Adler's convention: we have then 
\[\bdn^\star(V,K)= {\sup_{d\in \mathbb{N}}}^\star{\kappa_{inp}^d}^\star(K).\]
    In particular, if $K$ is strong (that is, of burden less or equal to $\aleph_{0,-}$), then so is $(V,K)$.

Also notice that if $K$ is finite, then one sees that $\bdn(V)=1$ if $V$ is of infinite dimension, and $\bdn(V)=0$ otherwise.

\begin{proof}
    If $V$ is of finite dimension $d$, then the structures $(V,K)$ and $\{K^d, K, \pi_i:K^d \rightarrow K, i \leq d\}$ are bi-interpretable on unary sets (modulo a
    choice of a $K$-basis for $V$). Then, one has $\bdn(V,K)= \bdn(K^d)=\kappa_{inp}^d(K)$.
    Assume now that the dimension of $V$ is infinite. For the same reason, $\bdn(V) \geq \sup_{d\in \mathbb{N}}\kappa_{inp}^d(K)$. 
    Consider an inp-pattern $P(x_K)$ in $(V,K)$, where $x_K$ is a $K$-variable. By stable embeddedness (Corollary \ref{Corollarykpure}), the array of parameters in $P(x_K)$ can be chosen in $K$; we have then in fact a pattern in $K$, of depth bounded by $\bdn(K)$.
    More interestingly, consider an inp-pattern $P(x)=\{\phi_i(x,y_i), (a_{\alpha,j})_{j<\omega},k_i\}_{\alpha<\mu}$ with a variable $x$ in $V$ and of depth $\mu$. We may assume the array of parameters to be mutually indiscernible. 
    By Proposition \ref{TheoremRelativeQE} and 
    elimination of the disjunction in inp-patterns (see for instance \cite{Tou20a}), we may assume that every formula $\phi(x,\bar{b},\bar{c})$ with variable $x\in V$ and parameters $\bar{b}\in V$, $\bar{c}\in K$ is a conjunction of:
    \begin{itemize}
        \item $\phi_K(\lambda(x,\bar{b}),\bar{c})$ called the $K$-part, where $\phi_K(\bar{z})$ is a formula in the language of fields and $\lambda(x,\bar{b})$ is a tuples of terms of the form $\lambda_{i,n}(x,b_{\sigma(1)}, \cdots,b_{\sigma(n)})$
for some integers $i\leq n$ and functions $\sigma:\{1,\dots,n\} \rightarrow \vert \bar{b}\vert$,
        \item $W(x,\bar{b})$ called the $W$-part, a conjunction of $W_{n+1}(x, b_{\sigma(1)},\dots,b_{\sigma(n)})$ or negation of such
for some integer $n$ and functions $\sigma: \{1,\dots,n\} \rightarrow \vert \bar{b}\vert$.
    \end{itemize}

Let $d$ be a realisation of the first column of the pattern $P(x)$:
\[d \models \{\phi_\alpha(x,a_{\alpha,0})\}_{\alpha<\mu}. \]

\begin{claim}\label{ClaimLambdaTerms}
    We may assume that if a given term $\lambda_{n,i}(x,y_{\alpha}^1,\dots,y_{\alpha}^n)$ occurs in (the $K$-part of)
    
    a line $\phi_\alpha(x,\bar{y}_\alpha)$, then 
    \[\phi_\alpha(x,\bar{y}_\alpha) \vdash W_n(y_{\alpha}^1,\dots,y_{\alpha}^n) \wedge \neg W_{n+1}(x,y_{\alpha}^1,\dots,y_{\alpha}^n).\]
    
\end{claim}
 \begin{proof}
 
 Assume that a given term $\lambda_{n,i}(x,y^1,\dots,y^n)$ occurs in the $K$-part  
\[\phi_K(\lambda_{n,i}(x,y^1,\dots,y^n),\lambda(x,\bar{y}_V),\bar{y}_K)\]
of the formula of the first line (we drop the index $\alpha=0$ for clarity).

We simply notice that it is equivalent to:
    \begin{align*}
& W_n(y^1,\dots,y^n) \wedge \neg W_{n+1}(x,y^1,\dots,y^n) \wedge \phi_K(\lambda_{n,i}(x,\bar{y}_V),\lambda(x,\bar{y}_V),\bar{y}_K) \\
     \bigvee & \neg \left(W_n(y^1,\dots,y^n)\wedge \neg W_{n+1}(x,y^1,\dots,y^n) \right)  \wedge \phi_K(0,\lambda(x,\bar{y}_V),\bar{y}_K).
    \end{align*}
    
    Then $(d,a_{0,0})$ satisfies either one or the other disjunct. We may remove the other term, as paths will still be consistent by mutual indiscernibility.
 \end{proof}

    For each line $\{\phi_{\alpha}(x,a_{\alpha,j})\}_{j< \omega}$, denote by $V_{\alpha,j}$ the intersection of the vector spaces generated by $b_{\alpha,j}^1, \dots, b_{\alpha,j}^n$ for each occurring terms $\lambda(x,b_{\alpha,j}^1, \dots, b_{\alpha,j}^n)$
and 
\[V_{\alpha}:= \bigcap_{j}  V_{\alpha,j}.\]
By indiscernibility, $V_{\alpha,j}$ has dimension $d_{\alpha} \in \mathbb{N}$ independent of $j$; and there is a $k_\alpha\in \mathbb{N}$ such that 
$\bigcap_{i<k_\alpha}V_{\alpha,j_{i}}= V_{\alpha}$ for all $j_0< \cdots < j_{k_\alpha} \in \mathbb{N}$.

We distinguish two cases:\\
\textbf{Case 1}: There is $\alpha<\mu$ such that $d \in V_{\alpha}$. Then we have in fact an inp-pattern in $V_{\alpha}$ (as by Corollary \ref{Corollarykpure}, $V_{\alpha}$ is stably embedded). We may conclude that $\mu\leq \kappa_{inp}^{\dim_K(V_\alpha)}$.\\
\textbf{Case 2}: There is no such $\alpha$. By mutual indiscernibility, for every path $f:\mu \rightarrow \omega$, we may find a solution
\[b_f \models \{\phi_{\alpha}(x,a_{\alpha,f(\alpha)})\}_{\alpha<\mu}\]
such that for all $\alpha<\mu$, $b_f \notin V_\alpha$.
In other words, the $W$-part form an inp-pattern 
\[P_W(x):=\{x\in V_{\alpha,j} \setminus V_{\alpha}\}_{\alpha<\mu,j<\omega}\]
where each line $\alpha$ is $k_\alpha$-inconsistent.
We show that such a pattern must be in fact of finite depth, bounded by  $\min_{\alpha < \mu}d_{\alpha}$.
\begin{claim}
    In the pattern $P_W(x)$, the dimension $d_{n}$ of $V_{n,j}$ is at least $n$.
\end{claim}
\begin{proof}
We show by induction on $n$ that $\bigcap_{m\geq n} V_{m,0}\setminus V_m$ has dimension at least $n$. The case $n=1$ is given by consistency of paths.
Assume $V_{n,0}\setminus V_n \cap \bigcap_{m>n} V_{m,0}\setminus V_m$ has at least $n$ linearly independent vectors. By mutual indicernibility, so is $V_{n,j}\setminus V_n \cap \bigcap_{m>n} V_{m,0}\setminus V_m$ for every $j< \omega$. As $\{x\in V_{n,j} \setminus V_{n}\}_{j<\omega}$ is $k_{n}$-inconsistent, one can find at least $n+1$ linearly independent vectors in  $\bigcap_{m>n} V_{m,0}\setminus V_m$. This concludes the induction.
\end{proof}

Of course, any two lines are interchangeable, which lead to conclude that $\mu\leq d_{\alpha}$ for all $\alpha<\mu$, and $\mu$ is in particular finite. This concludes our proof.

 

\end{proof}

\begin{example}
Consider for $n\in \mathbb{N}$, the two-sorted structure of $\mathbb{C}$-vector space $V_n$ of dimension $n$. Then $\bdn(V_n)=n$ for all $n\in \mathbb{N}$. However, $V:=\prod_\mathcal{U} v_n$ is simply a $K$-vector space of infinite dimension, and has burden $\aleph_{0,-}$. The reason is that the patterns witnessing $\bdn(V)>n$ are non uniform. Indeed, these are (after simplification) of the form $\phi(x,\lambda_i) \equiv \exists_{k\neq i,k<n} \lambda_k \in K \ x= \sum_{k<n} \lambda_k v_k $ where $v_0,\dots,V_{n-1}$ are $n$ linearly independent  vectors. In particular, if we want to increase the number of lines of the pattern, one must change all formulas of the pattern.
\end{example}

\section{Lexicographic products}
The lexicographic product of relational structures was studied by Meir in \cite{Mei16}. It is method to construct an $\mathrm{L}$-structure $\mathcal{M}[\mathcal{N}]$ from two $\mathrm{L}$-structures $\mathcal{M}$ and $\mathcal{N}$, where $\mathrm{L}$ is a relational language. To some aspects, the study of lexicographic products is similar to that of pure short exact sequences of abelian groups 
$0 \rightarrow A \rightarrow B \rightarrow C \rightarrow 0$. The reader can for instance compare Theorem \ref{ThmBdnLexicoProd} with \cite[Theorem 2.2]{Tou20a}. We will push the analogy further in the next paragraph when we will prove a reduction principle for valued vector spaces down to an (enriched) lexicographic products.

\subsection{Definition}
We recall here the definitions related to lexicographic product given by Meir in \cite{Mei19}. We will then take some conventions and change the language to fit better the context of this paper. Let $\mathrm{L}$ be a relational language.
\begin{definition} [{Meir, \cite[definition 0.0.1]{Mei19}}]
    Let $\mathcal{M}$ be an $\mathrm{L}$-structure and $\mathfrak{N}=\{\mathcal{N}_a\}_{a\in \mathcal{M}}$ be a collection of $\mathrm{L}$-structures indexed by $\mathcal{M}$. We consider the language $\mathrm{L}_{\mathbb{U},s} =\mathrm{L} \cup \{R^{\mathbb{U}} \}_{R \in \mathrm{L}}  \cup  \{s\}$ where $R^{\mathbb{U}}$ are new unary predicates and $s$ is a binary predicate. 
    The \textit{generalised lexicographic product} $\mathcal{M}[\mathfrak{N}^s]^{\mathbb{U}}$ of $\mathcal{M}$ and $\mathfrak{N}$ is the $\mathrm{L}_{\mathbb{U},s}$-structure of base set $\cup \mathfrak{N} := \cup_{a\in \mathcal{M}} \{a\}\times \mathcal{N}_a$ where the relations are interpreted as follows:
     \begin{itemize}
         \item $s^{\mathcal{M}[\mathfrak{N}^s]^{\mathbb{U}}} := \{ ((a,b),(a,b')) \ \vert \ a\in M, b,b' \in \mathcal{N}_a\},$ 
         \item if $R\in \mathrm{L}$ is an $n$-ary predicate, 
         \begin{align*}
            R^{\mathcal{M}[\mathfrak{N}^s]^{\mathbb{U}}} :=& \\
            & \left\{\left((a,b_1),\ldots, (a,b_n)\right) \ \vert \ a\in M \ \text{and } \mathcal{N}_a \models R(b_1,\ldots,b_n)\right\} \ \cup \\
            & \left \{\left((a_1,b_1),\ldots, (a_n,b_n)\right) \ \vert \ {\bigvee}_{0\leq i < j \leq n}a_i \neq a_j \ \text{and } \mathcal{M} \models R(a_1,\ldots,a_n) \right\},
         \end{align*}
         \item if $R\in \mathrm{L}$ is an $n$-ary predicate,
         \begin{align*}
             {R^{\mathbb{U}}}^{\mathcal{M}[\mathfrak{N}^s]^{\mathbb{U}}} :=& \bigg\{(a,b) \ \vert \ \mathcal{M}\models R(\underbrace{a, \ldots, a}_{n \text{ times}})\bigg \}.
         \end{align*} 
     \end{itemize}
     
     We called $\mathcal{M}[\mathfrak{N}^s]^{\mathbb{U}}$ the \textit{lexicographic product} (omitting `generalised' ) of $\mathcal{M}$ and $\mathfrak{N}$, or \textit{the lexicographic sum of $\mathfrak{N}$ with respect to $\mathcal{M}$}.
\end{definition}

     We will later abbreviate the notation and write $\mathcal{M}[\mathfrak{N}]$ instead of $\mathcal{M}[\mathfrak{N}^s]^{\mathbb{U}}$ to denote the structure (in any bi-interpretable language) of the lexicographic sum of $\mathfrak{N}$ with respect to $\mathcal{M}$. 
     If $\mathcal{N}_{a}= \mathcal{N}$ is constant for all $a\in \mathbb{M}$, we write $\mathcal{M}[\mathcal{N}]$ instead of $\mathcal{M}[\mathfrak{N}]$. Following a well-known terminology for ordered abelian groups, we occasionally call the $\mathcal{N}_a$'s the \emph{ribs} of $\mathcal{M}[\mathfrak{N}]$ and $\mathcal{M}$ the \emph{spine} of $\mathcal{M}[\mathfrak{N}]$.

    \begin{example}
        Consider the language of graphs $\mathrm{L}= \{R\}$, and let $\mathcal{M}$ and $\mathfrak{N}$ be the graphs:
        
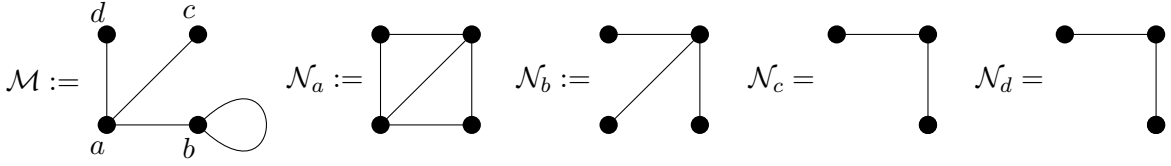
\begin{center} 
\begin{tikzpicture}[scale=1.2]\clip(-1.2,-0.4) rectangle (12,1.5);

\draw (-0.7,0.5) node {$\mathcal{M}:=$};
\fill (0,0) circle (0.1);
\draw (-0.1,-0.25) node {$a$};
\fill (1,0) circle (0.1);
\draw (0.9,-0.25) node {$b$};
\fill (1,1) circle (0.1);
\draw (0.9,1.25) node {$c$};
\fill (0,1) circle (0.1);
\draw (-0.1,1.25) node {$d$};
\draw (0,0) -- (1,0);
\draw (0,0) -- (0,1);
\draw (0,0) -- (1,1);
\draw (1,0) .. controls (2,1) and (2,-1) .. (1,0);

\draw (2.4,0.5) node {$\mathcal{N}_a:=$};
\fill (4,0) circle (0.1);
\fill (3,0) circle (0.1);
\fill (3,1) circle (0.1);
\fill (4,1) circle (0.1);
\draw (3,0) -- (3,1) -- (4,1) -- (4,0) -- (3,0) -- (4,1);
\draw (4,1) ;
\begin{scope}[shift={(2.5,0)}]
\draw (2.4,0.5) node {$\mathcal{N}_b:=$};
\fill (4,0) circle (0.1);
\fill (3,0) circle (0.1);
\fill (3,1) circle (0.1);
\fill (4,1) circle (0.1);
\draw (3,1) -- (4,1);
\draw (3,0) -- (4,1) -- (4,0);
\end{scope}
\begin{scope}[shift={(5,0)}]
\draw (2.4,0.5) node {$\mathcal{N}_c=$};
\fill (4,0) circle (0.1);
\fill (4,0) circle (0.1);
\fill (3,1) circle (0.1);
\fill (4,1) circle (0.1);
\draw (3,1) -- (4,1);
\draw (4,1) -- (4,0);
\end{scope}
\begin{scope}[shift={(7.5,0)}]
\draw (2.4,0.5) node {$\mathcal{N}_d=$};
\fill (4,0) circle (0.1);
\fill (4,0) circle (0.1);
\fill (3,1) circle (0.1);
\fill (4,1) circle (0.1);
\draw (3,1) -- (4,1);
\draw (4,1) -- (4,0);
\end{scope}
\end{tikzpicture}
\end{center}

We obtain the following graph:
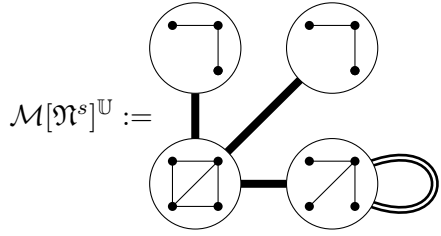
\begin{center} 
\begin{tikzpicture}[scale=0.6]\clip(-4,-1) rectangle (10,5);
\draw[line width=3pt] (0.5,1)--(0.5,2.5);
\draw[line width=3pt] (1.5,0.5)--(2.5,0.5);
\draw[line width=3pt] (0.5,0.5)--(3.5,3.5);
\draw[line width=1pt,double distance = 1pt] (3.5,0.5) .. controls (6.5,2.5) and (6.5,-1.5) .. (3.5,0.5);
\draw (-2,2) node {$\mathcal{M}[\mathfrak{N}^s]^{\mathbb{U}}:=$};
\begin{scope}[shift={(0,0)}] 
\filldraw[fill=white] (0.5,0.5) circle (1);
\fill (1,0) circle (0.1);
\fill (0,0) circle (0.1);
\fill (0,1) circle (0.1);
\fill (1,1) circle (0.1);
\draw (0,0) -- (0,1) -- (1,1) -- (1,0) -- (0,0) -- (1,1);

\end{scope}
\begin{scope}[shift={(3,0)}]
\filldraw[fill=white] (0.5,0.5) circle (1);
\fill (0,0) circle (0.1);
\fill (1,0) circle (0.1);
\fill (0,1) circle (0.1);
\fill (1,1) circle (0.1);
\draw (0,1) -- (1,1);
\draw (0,0) -- (1,1) -- (1,0);

\end{scope}
\begin{scope}[shift={(0,3)}]
\filldraw[fill=white] (0.5,0.5) circle (1);
\fill (1,0) circle (0.1);
\fill (1,0) circle (0.1);
\fill (0,1) circle (0.1);
\fill (1,1) circle (0.1);
\draw (0,1) -- (1,1);
\draw (1,1) -- (1,0);

\end{scope}
\begin{scope}[shift={(3,3)}]
\filldraw[fill=white] (0.5,0.5) circle (1);
\fill (1,0) circle (0.1);
\fill (1,0) circle (0.1);
\fill (0,1) circle (0.1);
\fill (1,1) circle (0.1);
\draw (0,1) -- (1,1);
\draw (1,1) -- (1,0);

\end{scope}
\end{tikzpicture}
\end{center}
where any point in a circle is related to any point from any linked circle. Notice that the loop is given by the predicate $R^{\mathbb{U}}$ (and not by $R$). 
\end{example}

    Within $\mathcal{M}[\mathfrak{N}^s]^{\mathbb{U}}$, the structure $\mathcal{M}$ can be seen as an imaginary. It is indeed the quotient of $\mathcal{M}[\mathfrak{N}^s]^\mathbb{U}$ modulo the equivalence relation $s(\mathbf{x},\mathbf{y})$. We denote by $v: \mathcal{M}[\mathfrak{N}^s]^{\mathbb{U}} \rightarrow \mathcal{M}$ the projection. It is not an $\mathrm{L}$-homomorphism.  Indeed, if $R$ is an $n$-ary predicate of $\mathrm{L}$, one has:
    \[v(R^{\mathcal{M}[\mathfrak{N}^s]^\mathbb{U}})= R^{\mathcal{M}}\cup \left\{ (a,\dots,a) \ \vert \ \text{ $R^{\mathcal{N}_a}$ is not empty}\right\}.
    \]
    
  Nonetheless, we recover the $\mathrm{L}$-structure when we consider the additional symbol $R^{\mathbb{U}}$:
        \begin{align}\label{EquationInducedStructureSpine}
       R^{\mathcal{M}}=v(R^{\mathcal{M}[\mathfrak{N}^s]^\mathbb{U}}) \setminus \{(a,\ldots,a) \ \vert \ a\in \mathcal{M}\}  \cup \{(a,\ldots,a) \ \vert \ a\in v({R^{\mathbb{U}}}^{\mathcal{M}[\mathfrak{N}^s]^\mathbb{U}}) \}. \end{align}

    For $a\in M$, the equivalence class of $a$ can be identified to $\mathcal{N}_a$. Here, by definition 
    \[\begin{array}{cll}
         f_a: & \mathcal{N}_a \rightarrow \mathcal{M}[\mathfrak{N}^s]^{\mathbb{U}} \\
            & b \mapsto  (a,b)
    \end{array}\]
    is an $\mathrm{L}$-isomorphism onto $\{a\}\times N_a$. We will see that there is no additional structure on these sets (Corollary \ref{CoroStEm}).

    As noticed by Meir, we may consider the lexicographic sum of a class of structures $\mathfrak{N}$ in a language $\mathrm{L}_\mathfrak{N}$ with respect to a structure $\mathcal{M}$ given in another language $\mathrm{L}_\mathcal{M}$:
    one simply needs to consider the union $\mathrm{L}:= \mathrm{L}_\mathcal{M} \cup \mathrm{L}_\mathfrak{N}$ and interprets in $\mathcal{M}$ (resp. in $\mathcal{N}_a$) all predicates $P \in \mathrm{L} \setminus \mathrm{L}_\mathcal{M}$ (resp. $P \in \mathrm{L} \setminus \mathrm{L}_\mathfrak{N}$) by the trivial set. 
    One can also consider non-relational languages: one can replace all function symbols by predicates interpreted by the graph of the functions. Without loss of generality, we may choose the opposite convention and assume that $\mathrm{L}_\mathcal{M}$ and $\mathrm{L}_\mathfrak{N}$ are disjoint, non necessary relational. This is the convention that we will often take. We rewrite the definition for sake of clarification.

\begin{definition} \label{DefinitionLexicographicProduct2}
    Let $\mathcal{M}$ be an $\mathrm{L}_\mathcal{M}$-structure and $\mathfrak{N}=\{\mathcal{N}_a\}_{a\in \mathcal{M}}$ be a collection $\mathfrak{N}$ of $\mathrm{L}_\mathfrak{N}$-structures indexed by $\mathcal{M}$. The lexicographic sum --denoted by $\mathcal{M}[\mathfrak{N}]$-- of $\mathfrak{N}$ with respect to $\mathcal{M}$ is the multisorted structure with:
    
    \begin{itemize}
        \item a main sort with base set  $S:=\cup_{a\in M} \{a\}\times \mathcal{N}_a$, 
        \item a sort for the structure $\mathcal{M}$,
        \item the natural projection map $v: S\rightarrow \mathcal{M}$.
    \end{itemize}
     To distinguish symbols 
    in the main sort $S$ to symbols in the ribs $\mathcal{N}_a$, we denote by $\mathrm{L}_{\bullet, \mathfrak{N}} := \{P_\bullet\}_{P\in \mathrm{L}_\mathfrak{N}}$ a copy of $\mathrm{L}_\mathfrak{N}$\footnote{This notation is inspired by the notation used in \cite{CH11} for the language $\mathcal{L}_{syn}$.}.
    The set $S$ is equipped with an $\mathrm{L}_{\bullet,\mathfrak{N}}$-structure: 
    \begin{itemize}
        \item if $P\in \mathrm{L}_\mathfrak{N}$ is an $n$-ary predicate, then
         \begin{align*}
            P_\bullet^{\mathcal{M}[\mathfrak{N}]} :=  \left\{\left((a,b_1),\ldots, (a,b_n)\right) \ \vert \ a\in M \ \text{and } \mathcal{N}_a \models R(b_1,\ldots,b_n)\right\}.
         \end{align*}
         \item if $f\in \mathrm{L}_\mathfrak{N}$ is an $n$-ary function symbol, then
         \begin{align*}
            f_\bullet^{\mathcal{M}[\mathfrak{N}]}((a,b_1),\ldots, (a,b_n)) := \begin{cases} (a,f^{\mathcal{N}_a}(b_1,\cdots,b_n) ) & \text{ if } a=a_1=\cdots=a_n,\\
            u & \text{ otherwise}.
            \end{cases}.
         \end{align*}
        \end{itemize}
    where $u$ is a constant which stand for `undetermined'. We can also pick a constant in the language.          
  We denote by $\mathrm{L}_{\mathcal{M}[\mathfrak{N}]}$ the multisorted language 
  \[(S,\mathrm{L}_{\bullet,\mathfrak{N}}) \cup (M,\mathrm{L}_\mathcal{M}) \cup \{v:S \rightarrow M\}.\]

\end{definition}

    This operation on first order structures is not adequate for structures sharing a common algebraic structure: the lexicographic sum of groups with respect to a group is (almost) never a group. However, we want to consider lexicographic product of algebraic structures: in this paper, we will consider the lexicographic sum of vector spaces with respect to a linear order. 

    These two definitions are equivalent: the languages are different but the structures on the lexicographic products is the same. For instance, consider the lexicographic product $\mathcal{M}[\mathcal{N}]$ of two linear orders $\mathcal{M}=(M,<)$ and $\mathcal{N}=(N,\prec)$ in two different languages, we obtain a structure bi-interpretable with the standard lexicographic product $(M\times N, <_{lex})$:  $(a,b)<_{lex}(a',b')$ if and only if $(a,b)<(a',b')$ or $s((a,b),(a',b')) \wedge (a,b) \prec (a',b')$. More generally, we have:
        \begin{remark}
    Let $\mathcal{M}$ be an $\mathrm{L}_\mathcal{M}$-structure and $\mathfrak{N}=\{\mathcal{N}_a\}_{a\in \mathcal{M}}$ be a collection $\mathfrak{N}$ of $\mathrm{L}_\mathfrak{N}$-structures. Consider the lexicographic product  $\mathcal{M}[\mathfrak{N}]$ in the multisorted language $\mathrm{L}_{\mathcal{M}[\mathfrak{N}]}$, and the lexicographic product  $\mathcal{M}[\mathfrak{N}^s]^{\mathbb{U}}$ in the one sorted language $\mathrm{L}_{\mathbb{U},s}$ (where we replaced each function symbol by a predicate for the graph and set $\mathrm{L}:=\mathrm{L}_\mathcal{M} \cup \mathrm{L}_\mathfrak{N}$). These two structures are bi-interpretable.
    
\end{remark}

    \subsection{Relative quantifier elimination in lexicographic products}
    Let $\mathrm{L}$ be a relational language.
    Here is the result of quantifier elimination obtained by Meir:

    \begin{theorem}\label{TheoremMeir1}[Meir]
    Let $\mathcal{M}$ be an $\mathrm{L}$-structure and $\mathfrak{N}$ be a collection of $\mathrm{L}$-structures.
    Assume that $Th(\mathcal{M})$ and $Th(\mathfrak{N})$ -- the common theory of all $\mathcal{N}_a$ -- eliminates quantifiers in $\mathrm{L}$. Then $\mathcal{M}[\mathfrak{N}^s]^{\mathbb{U}}$ eliminates quantifiers in $\mathrm{L}_{\mathbb{U},s}$.

        

        
        
    \end{theorem}
        
        
    
    On can refer to {\cite[Theorem 2.6]{Mei16}} in the case where $\mathcal{M}$ is transitive (for the action of the automorphism group). The above version of the theorem can be found in Meir's thesis \cite[Theorem 1.1.4]{Mei19}. 
    
    Notice that since $\mathrm{L}$ is a relational language (in particular, with no constant), a theory eliminates quantifiers only if it is complete. We state now a slightly more general theorem which does not necessary require that the $\mathcal{N}_a$ are elementary equivalent: 
    
    \begin{theorem}\label{TheoremMeir2}
        Assume that for every sentence $\phi \in \mathrm{L}$, the set $\{a\in \mathcal{M} \ \vert \ \mathcal{N}_a \models \phi \}$ is $\emptyset$-definable in $\mathcal{M}$.
        Assume also:
        \begin{itemize}
            \item[$(QE)_\mathfrak{N}$] For each $\mathrm{L}$-formula  $\phi(\bar{x})$ there are $k \in \mathbb{N}$ and quantifier-free $\mathrm{L}$-formulas $\phi_i(\bar{x})$, $i<k$ , such that for all models $\mathcal{N}$ of $Th(\mathfrak{N})$ there is $i<k$ such that:
        $\mathcal{N} \models \forall \bar{x} \ \phi(\bar{x}) \leftrightarrow \phi_i(\bar{x})$. \label{ConditionWeakEliminationQuantifier}
            
        \item[$(QE)_{\mathcal{M}}$] $\mathcal{M}$ eliminates quantifiers in $\mathrm{L}$.  
        \end{itemize}
        Then $\mathcal{M}[\mathfrak{N}^s]^{\mathbb{U}}$ eliminates quantifiers $\mathrm{L}_{\mathbb{U},s}$.   
\end{theorem}

    Before proving this more general statement, let us consider an elementary extension of $\mathcal{M}[\mathfrak{N}^s]^{\mathbb{U}}$. It is of the form $\tilde{\mathcal{M}}[\tilde{\mathfrak{N}}]$ where $\tilde{\mathfrak{N}}$ is a collection of $L$-structures $\tilde{\mathcal{N}}_{\tilde{a}}$ such that:
    \begin{itemize}
        \item $\tilde{\mathcal{M}}$ is an elementary extension of $ \mathcal{M}$
        \item if $a\in M$, then $\tilde{\mathcal{N}}_{a}$ is an elementary extension of $\mathcal{N}_a$.
    \end{itemize}
        We prove now that for an $\tilde{a} \in \tilde{\mathcal{M}}$, $\tilde{\mathcal{N}}_{\tilde{a}}$ is elementary equivalent to an ultraproduct $\prod_{\mathcal{U}}\mathcal{N}_a$. The type $\tp(\tilde{a})=:\rho\in \mathcal{S}_1^{Th(\mathcal{M})}$ can be seen as a filter $\{\phi(\mathcal{M}) \ \vert \ \phi(x)\in \rho \}$ on the set of definable (unary) subsets of $\mathcal{M}$. Let $\mathcal{U}$ be any ultrafilter on $\mathcal{M}$ containing $\rho$. It makes sense to consider the ultraproduct  $\prod_{\mathcal{U}}\mathcal{N}_a$. By assumption, the set $\{a\in \mathcal{M} \ \vert \ \mathcal{N}_a \models \phi \}$ is $\emptyset$-definable in $\mathcal{M}$ for all sentences $\phi$. By Łoś' theorem, this ultraproduct is independent of the choice of $\mathcal{U}$ extending $\rho$.
        
        \begin{fact}
            Let $\rho\in \mathcal{S}_1^{Th(\mathcal{M})}$.
            We denote by $\mathcal{N}_\rho$, $\prod_{\rho}\mathcal{N}_a$ or $\prod_{\rho}\mathfrak{N}$ the ultraproduct 
            $\prod_{\mathcal{U}}\mathcal{N}_a$
            for any ultrafilter as above.
              If $\tilde{\mathcal{M}}[\tilde{\mathfrak{N}}]$ is an elementary extension and $\tilde{a}$ is an element of $\tilde{\mathcal{M}}$ of type $\rho$, then we have that 
            \[\tilde{\mathfrak{N}}_{\tilde{a}} \equiv \mathcal{N}_\rho.\]
        \end{fact}

 Now, using $(QE)_{\mathcal{M}}$ and compactness, one can see that the condition $(QE)_\mathfrak{N}$ can be replaced by the following: 
      \begin{itemize}
            \item[$(QE')_\mathfrak{N}$] For all $\rho \in S_1^{Th(\mathcal{M})}$, $\mathcal{N}_\rho$ eliminates quantifiers in $\mathrm{L}$.
        \end{itemize}

    We can represent this theorem with the following `reduction diagram':
     \begin{center}
                \begin{tikzpicture}
                    \node{$\mathcal{M}[\mathfrak{N}^s]^{\mathbb{U}}$ }
                        child {  node {$\mathcal{M}$}}
                        child { node {$\mathcal{N_\rho} $ } edge from parent node[right,draw=none] {$\rho\in S_1^{\Th(\mathcal{M})}$}
                        };
                \end{tikzpicture}
            \end{center}
    We have a reduction diagram in the sense that if all structures in the leaves admit quantifier elimination, then so does the structure at the root. However, for a given $\rho \in S_1^{\Th(\mathcal{M})}$, there can be many ribs elementary equivalent to $\mathcal{N_\rho}$, and they do not need to be isomorphic. In particular, the ribs don't have necessarily the same cardinality.

    \begin{proof}
    It can be deduced from the proof of Theorem \ref{TheoremMeir1} with the following modification: consider a formula $\phi(\bar{v})$ of the form $\exists w \ \theta_i(\bar{v},w)$ where $\theta_i(\bar{v},w)$ are atomic and negated atomic formulas. 
    As in the proof of \cite[Theorem 1.1.4]{Mei19}, we may write equation (1.4) and assume that there are two formulas $\phi_M(\bar{v})$ and $\phi_N(\bar{v})$ such that for all tuples $(\bar{a},\bar{b})$:
    \begin{align}\label{equationMN}
    \mathcal{M}[\mathfrak{N}^s]^{\mathbb{U}} \models \phi((\bar{a},\bar{b})) \text{ if and only if } \mathcal{M}\models  \phi_M(\bar{a}) \text{ and } \mathcal{N}_{a_0}\models \bar \phi_N(\bar{b})  \text{ where $a_0$ is the first element of $\bar{a}$. }\end{align} 
    
    Then, by our quantifier elimination assumption, there are an $\mathrm{L}$-formula $\phi_M'(\bar{v})$ and two finite sets of quantifier-free $\mathrm{L}$-formulas $\{\phi_{N,0}(\bar{v}), \dots ,\phi_{N,k-1}(\bar{v})\}$ and 
    $\{\phi_{M,0}(z), \dots ,\phi_{M,k-1}(z)\}$ such that:
    \begin{itemize}
        \item for all $a\in \mathcal{M}$, there is $i<k$ such that $\mathcal{N}_a \models \forall \bar{v} \ \phi_{N,i}(\bar{v}) \leftrightarrow \phi_N(\bar{v})$,
        \item $\mathcal{M} \models \forall \bar{v} \ \phi_M(\bar{v}) \leftrightarrow \phi_M'(\bar{v}),$
        \item for all $i<k$, $\{ a\in \mathcal{M} \ \vert \ \mathcal{N}_a\models \forall \bar{v} \ \phi_{N,i}(\bar{v}) \leftrightarrow \phi_N(\bar{v})\} = \phi_{M,i}(\mathcal{M})$.
    \end{itemize}
    
    Then for all tuples $(\bar{a},\bar{b})$,
        \begin{align}
    \mathcal{M}[\mathfrak{N}^s]^{\mathbb{U}} \models \phi((\bar{a},\bar{b})) \text{ if and only if} & \text{  there is $i<k$ such that } \mathcal{M}\models  \phi_{M}'(\bar{a}) \wedge \phi_{M,i}(a_0)  \\ \notag
    & \text{and } \mathcal{N}_{a_0}\models \phi_{N,i}(\bar{b}).\end{align} 
    
    Then, we may conclude as in the proof of \cite[Theorem 1.1.4]{Mei19} using \cite[Lemma 1.1.7]{Mei19}  and \cite[Lemma 1.1.12]{Mei19}.
    
    \end{proof}

    We discuss here some immediate corollaries. Let us denote by $\restriction{\mathcal{M}}{\mathbb{U}}$ the set $v(\mathcal{M}[\mathfrak{N}^s]^{\mathbb{U}})$ -- the projection of $\mathcal{M}[\mathfrak{N}^s]^{\mathbb{U}}$ modulo the equivalence relation $s$-- with the full induced structure, and for $a\in M$, let $\{a\}\times \mathcal{N}_a \subset \mathcal{M}[\mathfrak{N}]$ be the set of elements with first coordinate equal to $a$, equipped with its full induced structure. 
    \begin{corollary}\label{CoroStEm}
        Assume that for every sentence $\phi \in \mathrm{L}$, the set $\{a\in \mathcal{M} \ \vert \ \mathcal{N}_a \models \phi \}$ is $\emptyset$-definable in $\mathcal{M}$.
        The structures $\restriction{\mathcal{M}}{\mathbb{U}}$ and $\{a\}\times \mathcal{N}_a$ for $a\in M$ are stably embedded and setwise orthogonal. 
        The structures $\restriction{\mathcal{M}}{\mathbb{U}}$ and $\mathcal{M}$ on the one hand, and $\{a\}\times \mathcal{N}_a$ and $\mathcal{N}_a$ on the other hand, have the same definable sets.
    \end{corollary}
    \begin{proof}
        We enriched the language $\mathrm{L}$ with predicates $P_{\phi,\mathcal{M}}$ and $P_{\phi(x),\mathfrak{N}}$ for each parameter-free formula $\phi(x)$, where 
        
        \begin{itemize}
            \item $P_{\phi(x),\mathcal{M}}$ is interpreted by $\phi(\mathcal{M}^{\vert x \vert})$ in $\mathcal{M}$ and by the empty set in $\mathcal{N}_{a}$ for all $a\in \mathcal{M}$.
            \item $P_{\phi(x),\mathfrak{N}}$ is interpreted by $\phi(\mathcal{N}_a^{\vert x \vert})$ in $\mathcal{N}_a$ for all $a\in M$ and by the empty set in $\mathcal{M}$.
        \end{itemize}
        Then $\mathcal{M}[\mathfrak{N}^s]^{\mathbb{U}}$ eliminates quantifiers by Theorem \ref{TheoremMeir2} in the corresponding enrichment of $\mathrm{L}_{\mathbb{U},s}$. Without loss, we may assume that  $\mathcal{M}[\mathfrak{N}^s]^{\mathbb{U}}$ eliminates quantifiers in $\mathrm{L}_{\mathbb{U},s}$. It is then enough to look at the structure on $\restriction{\mathcal{M}}{\mathbb{U}}$ and $\{a\}\times \mathcal{N}_a$ induced by atomic formulas. We see that $\restriction{\mathcal{M}}{\mathbb{U}}$ and $\mathcal{M}$ has the same definable sets using (\ref{EquationInducedStructureSpine}). It is immediate that $\{a\}\times \mathcal{N}_a$ and $\mathcal{N}_a$ are the same structure. Orthogonality is also immediate. 
    \end{proof}

    \begin{examples}
    One can easily see that in all the examples below, $Th(\mathfrak{N})$ satisfies condition $(QE)_{\mathfrak{N}}$.
        \begin{enumerate} 
            \item Consider for each integer $n$ a model $\mathcal{N}_n$ of the theory of $n$ cross-cutting equivalence relations, in the a common language $\mathrm{L}=\{E_i\}_{i\in \mathbb{N}}$ (if $i\geq n$, then $E_i^{\mathcal{N}_n}= \emptyset$). Consider $\mathbb{N}$ as a countable model of $T_\infty$ and $\rho_{\infty}$ the unique non realised type in $\mathbb{N}$.  Then $\prod_{\rho_{\infty}}\mathcal{N}_n$ has $\omega$ cross-cutting equivalence relations and $\mathbb{N}[\mathcal{N}_n]$ admits quantifier elimination in $\mathrm{L}$.
            \item Again, consider for each integer $n$ a model $\mathcal{N}_n$ theory of $n$ cross-cutting equivalence relations, but a disjoint language $\mathrm{L}_n$ and set $\mathrm{L}=\bigcup_n \mathrm{L}_n$. Then $\prod_{\rho_{\infty}}\mathcal{N}_n$ has no equivalence relation.

    \end{enumerate}
    We take the occasion to construct some linear orders which admit quantifier elimination in some natural - but non-trivial - languages.
    \begin{enumerate}
            \setcounter{enumi}{2}
            
            \item Set \[\mathcal{N}_{n}= \begin{cases}(\mathbb{Z},<) & \text{ if $n$ is even}\\ (\mathbb{R},<) & \text{ otherwise.} \end{cases}\] Then the lexicographic product $S:= \mathbb{N}[\mathcal{N}_n]$ eliminates quantifiers in the language $\{<,d_n^0,P^0\}$ where $D^0$ is a unary predicate for $2\mathbb{N}[\mathbb{Z}] \subset S $ and $d_n^0$ is a binary predicates for pairs of elements $\{(2k,a), (2k,a+n)\}_{k,a\in \mathbb{N}}$.  
            \item Consider a discrete linear order $(D_0,<)$ and a dense linear order $(D_0',<)$. We define by induction linear orders $D_{k+1}$ and $D_{k+1}'$ as respectively the lexicographic products $\mathbb{Z}[\mathcal{N}_n]$ and $\mathbb{R}[\mathcal{N}_r']$ where 
             \[\mathcal{N}_{n}= \begin{cases}D_k & \text{ if $n$ is even}\\ D_k' & \text{ otherwise.}\end{cases}\]
             and 
             \[\mathcal{N}_{r}'= \begin{cases}D_k & \text{ if $r$ is rational}\\ D_k' & \text{ otherwise.}\end{cases}\]
             By induction, $D_{k}$ and $D_{k}'$ admit quantifier elimination in the language $\{<, D_l , d_{l,n}\}_{l<k, n<\omega}$ where $D_l$ is a unary predicate for the union of intervals which are models of $D_l$ and $d_{l,n}$ is a binary predicate for pairs of elements in $D_l$ separated by $n$ or less models of $D_{l-1}$ (with the convention that models of $D_{-1}$ are singletons).
            \item With the notation above, let us consider the lexicographic product $\mathbb{N}[D_n]$. It admits quantifier elimination in the language $\{<, D_l , d_{l,n}\}_{l<\omega, n<\omega}$. 
        \end{enumerate}
    \end{examples}
    
    
We present now one more generalisation, stated in the convention of Definition \ref{DefinitionLexicographicProduct2}. It concerns lexicographic product $\mathcal{M}[\mathfrak{N}]$ enriched with a uniformly interpretable structure $\mathcal{K}$ on the component $\mathcal{N}_a$, $a\in \mathcal{M}$.

\begin{definition}[Uniform interpretation]
Let $\mathrm{L}_\mathcal{M}$, $\mathrm{L}_\mathfrak{N}$ and $\mathrm{L}_\mathcal{K}$ be disjoint  (non necessary relational) languages.
    Let $\mathfrak{N}$ be a class of $\mathrm{L}_\mathfrak{N}$-structures $\mathcal{N}_a$, and $\mathcal{K}$ an $\mathrm{L}_K$-structure. We say that $\mathfrak{N}$ interprets $\mathcal{K}$ uniformly with a set of maps $\{\lambda_i: \mathcal{N}^{k_i} \rightarrow K\}_{\mathcal{N} \in \mathfrak{N},i\in I}$ if 
     \begin{itemize}
                \item there are $\mathrm{L}_\mathfrak{N}$-formulas $D(\bar{x})$ and $E(\bar{x})$ such that for all $\mathcal{N} \in \mathfrak{N}$, $E(\mathcal{N})$ is an equivalence relation on $D(\mathcal{N}^{\vert x \vert})$ and  $D(\mathcal{N}^{\vert x \vert}) / E(\mathcal{N}) \simeq K$ .
                \item for all function symbol $f$ of arity $k\in \mathbb{N}$ in $\mathrm{L}_K$, there is an $\mathrm{L}_\mathfrak{N}$-formula $\phi_f(\bar{y})$ such that for all $\mathcal{N} \in \mathfrak{N}$, $\phi_f(\mathcal{N}^{\vert y \vert}) \subseteq D(\mathcal{N}^{\vert x \vert})^{k+1} $ interprets the graph of $f$.
                \item for all relation symbol $R$ of arity $k\in \mathbb{N}$ in $\mathrm{L}_R$, there is an $\mathrm{L}_\mathfrak{N}$-formula $\phi_R(\bar{y})$ such that for all $\mathcal{N} \in \mathfrak{N}$, $\phi_R(\mathcal{N}^{\vert \bar{y} \vert}) \subseteq D(\mathcal{N}^{\vert \bar{x} \vert})^k$ interprets $R$.
                \item for all $i\in I$,  there is an $\mathrm{L}$-formula $\phi_{\lambda_i}(\bar{y})$ such that for all $\mathcal{N} \in \mathfrak{N}$, $\phi_{\lambda_i}(\mathcal{N}_a^{\vert \bar{x} \vert}) \subseteq \mathcal{N}^{k_i} \times D(\mathcal{N}^{\vert \bar{x} \vert})$ interprets the graph of $\lambda_i: \mathcal{N}^{k_i} \rightarrow K$.
            \end{itemize}
\end{definition}

If a structure $\mathcal{K}$ is uniformly interpreted in $\mathfrak{N}$ with a set of maps $\{\lambda_i: \mathcal{N}^{k_i} \rightarrow K\}_{\mathcal{N} \in \mathfrak{N},i\in I}$, we can enriched $\mathcal{M}[\mathfrak{N}]$ with a sort for $\mathcal{K}$ and function symbols for maps $\lambda_{\bullet,i}: {\mathcal{M}[\mathfrak{N}]}^k \rightarrow K$ define as follows:
    \[ \forall (a_1,b_1), \dots, (a_k,b_k) \in \mathcal{M}[\mathfrak{N}], \ \lambda_{\bullet,i}((a_1,b_1), \dots, (a_k,b_k)): = \begin{cases} \lambda_i^{\mathcal{N}_a}(b_1,\dots,b_k) & \text{ if } a_1=\cdots=a_k=a \\
    0 & \text{ otherwise.}\end{cases}\]
    where $0$ is a constant\footnote{if there is no such constant, one can formally add one point in $K$, and name it by a symbol $0$.} in $K$. 
Such enrichment will be called $\mathcal{K}$-enrichment or simply enrichment of 

The following is our main example and will discuss in detail in Section 3:
\begin{example}\label{ExampleEnrichedLexicographicProductVectorSpaces}
     Let $K$ be a field. Assume that $(\Gamma, <)$ is an order set and $\mathfrak{N}:= \{B_\gamma\}_{\gamma \in \Gamma}$ is a family of $K$-vector spaces, seen as structures in the one-sorted language $\mathcal{L}_W:= \{(V,0,+), (W_n)_n \}$ where $W_n$ denotes a predicates for $K$-independent $n$-tuples. By Fact \ref{FactPierce}, the field $K$ together with the lambda-functions $(\lambda_{n,i})_{i\leq n}$ are interpretable in  $\mathrm{L}_W$ . 
     Then $\Gamma[B_\gamma]_\gamma$ can be enriched with a sort for the field $\mathcal{K}$ and with maps 
     \[\lambda_{\bullet n,i}:=  ( (a,b),(a_1,b_1), \dots, (a_k,b_k) )\mapsto \begin{cases} \lambda_i^{\mathcal{N}_a}(b,b_1,\dots,b_k) & \text{ if } a_1=\cdots=a_k=a \\
    0 & \text{ otherwise.}\end{cases}\]
     for $i\leq n \in\mathbb{N}$.
\end{example}


\begin{theorem}\label{TheoremRelativeQuantifierEliminationEnrichedLexProd}
 Consider the lexicographic sum $\mathcal{S}:=\mathcal{M}[\mathfrak{N}]$ of a class of $\mathrm{L}_\mathfrak{N}$-structures $\mathfrak{N}:= \{\mathcal{N}_a\}_{a\in \mathcal{M}}$ with respect to an $\mathrm{L}_\mathcal{M}$-structure $\mathcal{M}$ in the language
    \[\mathrm{L}_{\mathcal{M}[\mathfrak{N}]}:= (S,\mathrm{L}_{\bullet,\mathfrak{N}}) \cup (\mathcal{M},\mathrm{L}_\mathcal{M})\cup \{v: S\rightarrow \mathcal{M}   \}.\]
    
    Assume that for all sentences $\phi\in \mathrm{L}_\mathfrak{N}$, the set $\{a\in \mathcal{M} \ \vert \ \mathcal{N}_a \models \phi \}$ is $\emptyset$-definable in $\mathcal{M}$.
    Assume also that $\mathfrak{N}$ interprets uniformly a structure $\mathcal{K}$ (in a language $\mathrm{L}_K$) with a non-empty  set of \emph{surjective} maps $\{\lambda_i: \mathcal{N}_a^{k_i} \rightarrow K\}_{a\in \mathcal{M},i\in I}$.
        Finally, assume that:
        \begin{itemize}
            \item [$(RQE)_\mathfrak{N}$] For all $\rho\in S_1^{\Th(\mathcal{M})}$, $\mathcal{N}_\rho$ admits quantifier elimination relative to $K$ in $\mathrm{L}_{\mathfrak{N}} \cup \mathrm{L}_K \cup \{\lambda_i \}_{i\in I}$.

        \end{itemize}
        Then $(\mathcal{M}[\mathfrak{N}],\mathcal{K},\{\lambda_{\bullet,i}\}_{i\in I})$ eliminates quantifiers relative to $\mathcal{M}$ and $K$ in $\mathrm{L}_{\mathcal{M}[\mathfrak{N}]}\cup \mathrm{L}_K \cup \{\lambda_{\bullet,i}\}_{i\in I}$.   
\end{theorem}
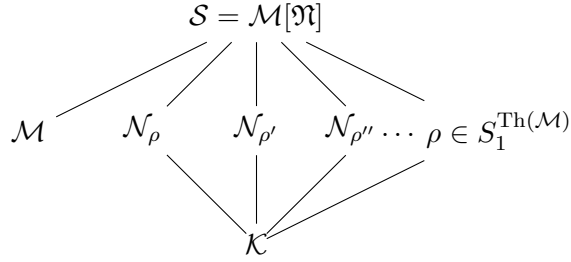
\begin{center}
\begin{tikzpicture}
                    \node {$\mathcal{S}=\mathcal{M}[\mathfrak{N}]$}
                            child { 
                            node { $\mathcal{M}$ } 
                            }
                            child { node { $\mathcal{N}_\rho $ }
                                child [missing] 
                                child [missing] 
                                child { node {} }
                            }
                            child { node { $\mathcal{N}_{\rho'}$ } 
                                child { node {$\mathcal{K}$} }
                            }
                            child { node { $\mathcal{N}_{\rho''}\ldots$ } 
                                child { node {
                                } }
                                    child [missing]
                                    child [missing] 
                            }
                            child { node {$ \quad \rho \in S_1^{\Th(\mathcal{M})} $}
                                    child { node {
                                } }
                                    child [missing]
                                    child [missing]
                                    child [missing]
                                    child [missing]
                        };
                \end{tikzpicture}
 \end{center}

Notice that this is a generalisation of the previous theorems, as we may always take $\mathcal{K}$  empty: the condition of relative quantifier elimination for the $\mathcal{N}_a$'s says that the $\mathcal{N}_a$ (absolutly) eliminates quantifier. 
It is important to notice that if $\mathcal{K}$ is not empty, then the structure $(\mathcal{M}[\mathfrak{N}],\mathcal{K},\{\lambda_{\bullet,i}\}_{i\in I})$ is an enrichment of $\mathcal{M}[\mathfrak{N}]$. In some sense, we have identified the interpretable copies of $\mathcal{K}$. Since the notation considerably change, we redo the proof completely, with less details.

 As Meir in \cite{Mei19}, we will also use the notion of $s$-diagrams (adapted to our notation): 
    \begin{definition}
        Let ${x}, {y}^0, \ldots  ,{y}^{n-1}$ be single variables.  A formula $\delta({x},{y}^0,\ldots,{y}^{n-1})$ is a \textit{ complete $s$-diagram in ${x}$ over ${y}^0, \ldots  ,{y}^{n-1}$} if it is a maximally consistent conjunction of $ v({x})=v({y}^l)$ or $v({x})\neq v({y}^l)$ for $l<n$. 
    \end{definition}

\begin{proof}
    We may consider the relational language obtained by replacing all function symbols in $\mathrm{L}_\mathcal{M}$ by a predicate interpreted by the graph of the function, and similarly with $\mathrm{L}_\mathfrak{N}$. Since $v: \mathcal{S} \rightarrow \mathcal{M}$ and (at least one) $\lambda_i$ are surjective, we only need to eliminate $\mathcal{S}$-sorted quantifier in formula with main-sorted variables, as we may replace a variable $x_\mathcal{M}$ in $\mathcal{M}$ (resp. $x_\mathcal{K}$ in $\mathcal{K}$) by a term $v(x)$ (resp. $\lambda_{\bullet,i}(x)$) where $x$ is a $\mathcal{S}$-sorted variable. Let $\Phi(x,x_1,\dots,x_n)$ be a $\mathcal{S}$-sorted-quantifier-free formula in $\mathrm{L}_{\mathcal{M}[\mathfrak{N}]}$ with main-sorted variables.  Let us show that  $\exists x \ \Phi(x,x_1,\dots,x_n)$
is $\mathcal{S}$-sorted-quantifier-free definable.
Consider the set $\Delta(x,x_1,\dots,x_n)$ of all $s$-diagrams of $x$ over $x_1,\dots,x_n$.
We may write \[\Phi(x,x_1,\dots,x_n) \equiv \bigvee_{\delta(x,x_1,\dots,x_n)\in \Delta} \Phi(x,x_1,\dots,x_n) \wedge \delta(x,x_1,\dots,x_n).\]
Using that for all $i \in I$ and $P\in \mathrm{L}_{\mathfrak{N}}$, $\mathcal{S}\models \lambda_{i,\bullet}(x,x_{i_1},\dots,x_{i_k}) \neq 0$ and $\mathcal{S}\models P_\bullet(x,x_{i_1},\dots,x_{i_k}) $ only if $v(x)=v(x_{i_1})=\cdots=v(x_{i_k})$, we may consider the terms of the disjunction one by one and remove all variables $x_i$ such that $\delta(x,x_1,\dots,x_n) \vdash v(x)\neq v(x_i)$.

We reset the notation and assume that $\Phi(x,x_1,\dots,x_n) \vdash v(x)=v(x_1) = \cdots = v(x_n)$.
Then there are some $\mathcal{S}$-sorted-quantifier-free $\mathrm{L}_{\bullet,\mathfrak{N}}\cup \mathrm{L}_K\cup \{\lambda_{\bullet,i}\}_{i\in I}$-formulas $\phi_{\bullet,k}(x,x_1,\dots,x_n) $ and $\mathrm{L}_\mathcal{M}$-formulas $\theta_i(x)$ such that
\[ \ \equiv \ \bigvee_k \phi_{\bullet,k}(x,x_1,\dots,x_n)\wedge \theta_i(v(x)) \wedge v(x)=v(x_1)= \cdots = v(x_n).\]

Again by considering each disjunct one at the time and reset the notation, we may simply assume that 
\[\Phi(x,x_1,\dots,x_n) \ \equiv \ \phi_{\bullet}(x,x_1,\dots,x_n)\wedge \theta(v(x)) \wedge v(x)=v(x_1)= \cdots = v(x_n).\]
We denote by $\phi(x_{},x_{1},\dots,x_{n})$
the $\mathrm{L}_\mathfrak{N}\cup \mathrm{L}_K\cup \{\lambda_i\}_{i\in I}$-formula obtained from $\phi_\bullet$ by replacing each occurrences of $P_\bullet$  by $P$ and each occurrences of $\lambda_{\bullet,i}$ by $\lambda_i$.

Fix an element $a\in \mathcal{M}$.
By relative quantifier elimination, there is an $\mathcal{N}$-sorted-quantifier-free $\mathrm{L}_\mathfrak{N}\cup \mathrm{L}_K\cup \{\lambda_i\}_{i\in I}$-formula $\psi_a(x_1,\dots,x_n)$ such that 
\[(\mathcal{N}_a,\mathcal{K},\{\lambda_i\}_{i\in I}) \models  \exists x \  \phi(x,x_1,\dots,x_n) \leftrightarrow \psi_a(x_1,\dots,x_n).\] 
We have for all $c_1=(a,b_1),\dots,c_n=(a,b_n) \in \{a\}\times \mathcal{N}_a$. 
    \begin{align*}
         \mathcal{S} \models \exists x \  \Phi(x,c_1,\dots,c_n)  & \Leftrightarrow  (\mathcal{N}_a, \mathcal{K},\{\lambda_i\}_{i\in I}) \models \exists x \ \phi(x,b_{1},\dots,b_{n}) \text{ and } \mathcal{M}\models \theta(a),\\
         & \Leftrightarrow (\mathcal{N}_a, \mathcal{K},\{\lambda_i\}_{i\in I}) \models \psi_a(b_{1},\dots,b_{n}) \text{ and } \mathcal{M}\models \theta(a) .
    \end{align*}
   
    Denote by $\psi_{\bullet,a}(x_{1},\dots,x_{n})$ the $\mathrm{L}_{\mathcal{M}[\mathfrak{N}]}$-formula obtained by replacing in $\psi_a(x_{1},\dots,x_{n})$ all occurrences of $P$ (resp. $\lambda_i$) by $P_\bullet$ (resp. $\lambda_{\bullet,i})$.
    We have:
    
    \begin{align}
         (\mathcal{N}_a, \mathcal{K},\{\lambda_i\}_{i\in I}) \models \psi_a(b_{1},\dots,b_{n})& 
         \Leftrightarrow (S, \mathcal{K}, \{\lambda_{\bullet,i}\}_{i\in I}) \models \psi_{\bullet,a}(c_{1},\dots,c_{n}).
    \end{align}

All together, this give us: 
\begin{align}\label{EquationEliminationEQuantifierEnrichedLexPro}
    \forall c_1,\dots,c_n \in \{a\}\times \mathcal{N}_a \quad 
    \mathcal{S} \models \exists x \  \Phi(x,c_1,\dots,c_n)  \Leftrightarrow \mathcal{S} \models \psi_{\bullet,a}(c_{1},\dots,c_{n}) \wedge \theta(v(c_1)). \end{align}

Notice that the above formula $\psi_a$ depends only on $a$. 

\begin{claim}
    For a given $\mathrm{L}_\mathfrak{N}$-formula $\psi(x_1,\cdots,x_n)$, the set of elements $a\in \mathcal{M}$ such that we may take $\psi_a=\psi$ in (\ref{EquationEliminationEQuantifierEnrichedLexPro}) is $\emptyset$-definable in $\mathcal{M}$.
    We denote this formula by $\delta_\psi(x)$.
\end{claim}
\begin{proof}
   Indeed, by uniform interpretation of $\mathcal{K}$ in the $\mathcal{N}_a$'s, there is an $\mathrm{L}_\mathfrak{N}$-formula 
$ \tilde{\psi}(x_1,\cdots,x_n)$ such that for all $a\in \mathcal{M}$ and all $b_1,\dots,b_n \in \mathcal{N}_a$, 
\[(\mathcal{N}_a, \mathcal{K},\{\lambda_i\}_{i\in I}) \models \psi(b_{1},\dots,b_{n}) \quad \Leftrightarrow \quad  \mathcal{N}_a \models \tilde{\psi}(b_{1},\dots,b_{n}).\]

Similarly, there is an $\mathrm{L}_\mathfrak{N}$-formula $\tilde{\phi}(x,x_1,\dots,x_n)$
such that for all $a\in \mathcal{M}$ and  $b,b_1,\dots,b_n \in \mathcal{N}_a $
\[(\mathcal{N}_a, \mathcal{K},\{\lambda_i\}_{i\in I}) \models \phi(b,b_{1},\dots,b_{n}) \quad \Leftrightarrow \quad  \mathcal{N}_a \models \tilde{\phi}(b,b_{1},\dots,b_{n}).\]
Then if  $\theta$ is the sentence
\[\forall x_1,\dots,x_n \ \exists x \ \tilde{\phi}(x,x_{1},\dots,x_{n}) \ \leftrightarrow \ \tilde{\psi}(x_{1},\dots,x_{n}),\]
we have by assumption that 
\[\{a \in \mathcal{M} \ \vert \ \mathcal{N}_a \models \theta\}\]
is $\emptyset$-definable.

\end{proof}

Assume that
there is no $\mathcal{N}$-sorted-quantifier-free  $\mathrm{L}_\mathfrak{N}\cup \mathrm{L}_K\cup\{\lambda_i\}_{i\in I}$-formula $\psi$ such that 
 \[\mathcal{S} \models \exists x \  \Phi(x,x_1,\dots,x_n) \leftrightarrow \psi_\bullet(x_1,\dots,x_n) \wedge \theta(v(x_1)) \wedge v(x_1) = \cdots v(x_n).\]
 This means that for any finite set of formula $\Psi$, there is an element $a_\Psi\in \mathcal{M}$ such that $ \mathcal{M} \models \neg \delta_\psi (a_\Psi)$ for all $\psi\in \Psi$.
 It follows that there is a type $\rho \in S_1^{\Th(\mathcal{M})}$  such that for all $\mathcal{N}$-sorted-quantifier-free $\mathrm{L}_\mathfrak{N}\cup \mathrm{L}_K\cup\{\lambda_i\}_{i\in I}$-formula $\psi(x_1,\dots,x_n)$ such that 
 \[(\mathcal{N}_\rho, \mathcal{K},\{\lambda_i\}_{i\in I}) \models \exists x \ \phi(x,x_1,\dots,x_n) \not\leftrightarrow \psi(x_1,\dots,x_n). \]
 This is a contradiction with quantifier elimination in $\mathcal{N}_\rho$ relative to $\mathcal{K}$. Therefore, we have for some $\mathcal{N}$-sorted-quantifier-free  $\mathrm{L}_\mathfrak{N}\cup \mathrm{L}_K\cup\{\lambda_i\}_{i\in I}$-formula $\psi$
  \[\mathcal{S} \models \exists x \  \Phi(x,x_1,\dots,x_n) \leftrightarrow \psi_\bullet(x_1,\dots,x_n) \wedge \theta(v(x_1)) \wedge v(x_1) = \cdots v(x_n).\]
 This concludes our proof.
\end{proof}
It follows from the proof the following:
\begin{remark}
    In the theorem above, we may add the condition that  $\mathcal{M}$ is $\vert \mathrm{L}_\mathfrak{N} \vert^+$-saturated and replace the condition
    \begin{itemize}
        \item[$(RQE)_\mathfrak{N}$] For all $\rho\in S_1^{\Th(\mathcal{M})}$, $\mathcal{N}_\rho$ admits quantifier elimination relative to $K$ in $\mathrm{L}_{\mathfrak{N}} \cup \mathrm{L}_K \cup \{\lambda_i \}_{i\in I}$,
        
            \end{itemize}
            by the condition
            \begin{itemize}
        \item[$(RQE')_\mathfrak{N}$] For all $a\in \mathcal{M}$, $\mathcal{N}_a$ admits quantifier elimination relative to $K$ in $\mathrm{L}_{\mathfrak{N}} \cup \mathrm{L}_K \cup \{\lambda_i \}_{i\in I}$.
            \end{itemize}
\end{remark}

Finaly, let us restate Corollary $\ref{CoroStEm}$:
 \begin{corollary}
        Under the same assumption as in Theorem \ref{TheoremRelativeQuantifierEliminationEnrichedLexProd}, we have:
        \begin{itemize}
            \item the spine $\mathcal{M}$, equipped with its structure in the language $\mathrm{L}_\mathcal{M}$,  is stably embedded and pure.
            \item for $a\in \mathcal{M}$, the rib $\{a\}\times \mathcal{N}_a$ seen as a structure in the language $\mathrm{L}_{\bullet,\mathcal{N}} \cup \mathrm{L}_\mathcal{K} \cup \{\lambda_{\bullet,i}\}_{i\in I}$ is stably embedded and pure. It is isomorphic to $(\mathcal{N}_a,\mathcal{K},\{\lambda_{i}\}_{i\in I})$.
            \item All the above structures are pairwise orthogonal.
            \item the sort $\mathcal{K}$ is stably embedded and pure. It is also orthogonal to $\mathcal{M}$.
        \end{itemize}
    \end{corollary}

\subsection{Transfers principles in lexicographic products}
    We prove now other transfer principles in lexicographic product, namely for model completeness, NIPness and burden.  For more readability, elements and variables from the sort $S$ will be written in bold character: $\mathbf{a}$, $\mathbf{x}$ etc. We will use the standard font for elements and variables from $\mathcal{M}$ or $\mathcal{N}_a$ and we will identify parameters $\mathbf{a}$ with a pair $(a_\mathcal{M}, a_\mathfrak{N})$.

    The following is a direct corollary of Theorem \ref{TheoremMeir2} / Theorem \ref{TheoremRelativeQuantifierEliminationEnrichedLexProd}:
 \begin{corollary}
        Consider the lexicographic sum $\mathcal{S}$ of a class $\mathfrak{N}:= \{\mathcal{N}_a\}_{a\in \mathcal{M}}$ of $\mathrm{L}_\mathfrak{N}$-structures with respect to an $\mathrm{L}_\mathcal{M}$-structure $\mathcal{M}$  in the language $\mathrm{L}_{\mathcal{M}[\mathfrak{N}]}$:
        \[\mathcal{S}:= (S,\mathrm{L}_{\bullet,\mathfrak{N}}) \cup (M,\mathrm{L}_\mathcal{M})\cup \{v:S\rightarrow M\}.\]
        Consider an extension $\mathcal{S}'$ of $\mathcal{S}$. Then $\mathcal{S}'$ is also a lexicographic sum of a class of $\mathrm{L}_\mathfrak{N}$-structures $\mathfrak{N}':= \{\mathcal{N}_a'\}_{a\in \mathcal{M}}$ with respect to an $\mathrm{L}_\mathcal{M}$-structure $\mathcal{M}'$:
        \[\mathcal{S}':= (S',\mathrm{L}_{\bullet,\mathfrak{N}}) \cup (M',\mathrm{L}_\mathcal{M})\cup \{v:S'\rightarrow M'\}.\]
        
        Assume:
        \begin{itemize}
            \item for every sentence $\phi  \in \mathrm{L}_\mathfrak{N}$, the set $\{a\in \mathcal{M} \ \vert \ \mathcal{N}_a \models \phi \}$ is $\emptyset$-definable in $\mathcal{M}$,
            \item $\mathcal{M}\preceq \mathcal{M}'$,
            \item for all $a \in \mathcal{M}$, $\mathcal{N}_a \preceq \mathcal{N'}_{a}$.
        \end{itemize}
  
        Then we have $\mathcal{S} \preceq \mathcal{S}'$.
    \end{corollary}
\begin{proof}
   We consider the morleyrisation $\mathrm{L}_\mathfrak{N}^{Mor}$ and $\mathrm{L}_\mathcal{M}^{Mor}$ (obtained from $\mathrm{L}_\mathfrak{N}$ resp. $\mathrm{L}_\mathcal{M}$ by adding predicates for each  $\mathrm{L}_\mathfrak{N}$-formula resp. $\mathrm{L}_\mathcal{M}$-formula). Then $\mathcal{S}$ and $\mathcal{S}'$ admits a natural expansion to the language \[(S,\mathrm{L}_{\bullet\mathfrak{N}}^{Mor}) \cup (\mathcal{M},\mathrm{L}_\mathcal{M}^{Mor}) \cup \{v\},\]
   where they eliminates quantifiers by Theorem \ref{TheoremRelativeQuantifierEliminationEnrichedLexProd}.
   Using that $\mathcal{M} \preceq \mathcal{M}'$ and that for all $a\in \mathcal{M}$, $\mathcal{N}_a \preceq \mathcal{N}_a'$, we see that $\mathcal{S}$ is still a substructure of $\mathcal{S}'$ in this language. It follows that $\mathcal{S}\preceq \mathcal{S}'$, which concludes our proof.
   Remark that it automatically follows that for all $\rho \in S_1^{Th(\mathcal{M})}$, $ \mathcal{N}_\rho \preceq \mathcal{N'}_{\rho}$.
\end{proof}
    \addtocounter{theorem}{-1}
     A direct corollary of Theorem \ref{TheoremRelativeQuantifierEliminationEnrichedLexProd} and Remark \ref{RemarkModelCompletnessRQECLoseSort} gives  a similar result in the  enriched setting:
    \begin{corollary}[Enriched version]
     Under the same condition of Theorem \ref{TheoremRelativeQuantifierEliminationEnrichedLexProd},
    an extension of enriched lexicographic product $(\mathcal{M}[\mathfrak{N}],\mathcal{K},\{\lambda_{\bullet,i}\}_{i\in I}) \subseteq (\mathcal{M}'[\mathfrak{N}'],\mathcal{K}',\{\lambda_{\bullet,i}'\}_{i\in I})$ is elementary if and only if $\mathcal{K} \subseteq \mathcal{K}'$ and $\mathcal{M}\subseteq \mathcal{M}'$ are elementary.
    
    \end{corollary}

    The following corollary is due to Meir in the case that the common theory of the ribs $\mathcal{N}_a$ is $\NIP$ \cite[Corollary 1.6.13]{Mei19}. We give now a slightly more general version and give a proof for completeness.

    \begin{corollary}
        Let $\mathcal{S}$ be the lexicographic sum of a class of $\mathrm{L}_\mathfrak{N}$-structures $\mathfrak{N}:= \{\mathcal{N}_a\}_{a\in \mathcal{M}}$ with respect to an $\mathrm{L}_\mathcal{M}$-structure $\mathcal{M}$.
      Assume:
        \begin{itemize}
            \item for every sentence $\phi \in \mathrm{L}_\mathfrak{N}$, the set $\{a\in \mathcal{M} \ \vert \ \mathcal{N}_a \models \phi \}$ is $\emptyset$-definable in $\mathcal{M}$,
            \item $\mathcal{M}$ is NIP,
            \item for all $\rho \in S_1^{Th(\mathcal{M})}$, $\mathcal{N}_\rho$ is NIP.
        \end{itemize}
   Then $\mathcal{S}$ is NIP.

    \end{corollary}

    \begin{remark}
        The assumption that for all types $\rho \in S_1^{Th(\mathcal{M})}$ the structure $\mathcal{N}_\rho$ is NIP  is necessary and can't be replaced by
        the assumption that for all $a \in \mathcal{M}$ the structure $\mathcal{N}_a$ is NIP. This is due to the fact that a ultraproduct of $\NIP$ structures can have $\IP$. For instance, if $G_n$ is the finite bipartite graph on $n \times 2^n$ such that $R(i,J)$ if and only if $i\in J$, then an elementary extension of $\mathbb{N}[G_n]$ contains a model of $\prod_{p_{\infty}}G_n$, which is the random graph, and so has $\IP$.
    \end{remark}

     \begin{proof}
        As in the previous proof, we morleyrise the languages $\mathrm{L}_\mathcal{M}$ and $\mathrm{L}_\mathfrak{N}$ so that we may assume that  $(S,\mathrm{L}_{\bullet,\mathfrak{N}}) \cup (\mathcal{M},\mathrm{L}_\mathcal{M})\cup\{v\}$ eliminates quantifiers in the language $\mathrm{L}_{\mathcal{M}[\mathfrak{N}]}$.\\
        A formula $\phi(\mathbf{x},\mathbf{y})$ with $\vert \mathbf{x} \vert = 1$ is equivalent to a Boolean combination of formulas of the form:
            \begin{itemize}
                \item $\delta(\mathbf{x},\mathbf{y})$, a complete $s$-diagram in $\mathbf{x}$ over $\mathbf{y}^1,\ldots, \mathbf{y}^{\vert\mathbf{y} \vert}$,
                \item $ P_\bullet(\mathbf{x},\mathbf{y})$ for $P\in \mathrm{L}_\mathfrak{N}$,
                \item $ P(v(\mathbf{x}),v(\mathbf{y}))$ for $P\in \mathrm{L}_\mathcal{M}$.
            \end{itemize}
            
            Since $\mathcal{M}$ is NIP, the formula $P(v(\mathbf{x}),v(\mathbf{y}))$ is automatically NIP. It is enough to show that the conjunction of a complete diagram $\delta(\mathbf{x},\mathbf{y})$ and an atomic formula $P_\bullet(\mathbf{x},\mathbf{y})$ with $\vert x\vert=1$ and $P\in \mathrm{L}_\mathfrak{N}$ is $\NIP$. Since $P_\bullet(\mathbf{a},\mathbf{b}) \neq \emptyset$ if and only if $v(\mathbf{a})=v(\mathbf{b}_1)=\cdots=v(\mathbf{b}_{\vert \mathbf{b} \vert})$, we may assume that $\delta(\mathbf{x},\mathbf{y}) \vdash v(\mathbf{x})=v(\mathbf{y}_1)=\cdots=v(\mathbf{y}_{\vert \mathbf{y} \vert})$.
            We prove that $P_\bullet(\mathbf{x},\mathbf{y}) \wedge \delta(\mathbf{x},\mathbf{y})$ is NIP using a characterisation with indiscernible sequences.
            Consider an indiscernible sequence $(\mathbf{a}_i)_{i\in I}$ of singletons and a parameter $\mathbf{b} \in \mathbb{S}^{\vert \mathbf{y}\vert}$ in a elementary extension of $\mathcal{S}$.
        We write $\mathbf{a}_i=(a_{\mathcal{M},i},a_{\mathfrak{N},i})$ and similarly,  $\mathbf{b}^l=(b_\mathcal{M}^l,b_{\mathfrak{N}}^l)$ for $l<\vert \mathbf{y} \vert$
    If $\neg (v(\mathbf{b}_1)=\cdots=v(\mathbf{b}_{\vert \mathbf{y} \vert})$, then 
    \[S\models \neg P_\bullet(\mathbf{a}_i,\mathbf{b}) \wedge \delta(\mathbf{a}_i,\mathbf{b})\]
    for all $i\in I$. 
    Assume $v(\mathbf{b}_1)=\cdots=v(\mathbf{b}_{\vert \mathbf{y} \vert})$ and let $\rho$ be the type of $v(\mathbf{b}_1)$ in $S^{\Th(\mathcal{M})}_1$. Then
    
        \begin{align*}
\mathcal{S} \models  P_\bullet(\mathbf{a}_i,\mathbf{b}) \wedge \delta(\mathbf{a}_i,\mathbf{b})  & \Leftrightarrow \mathcal{N}_\rho \models P(a_{\mathfrak{N},i},b_\mathfrak{N}). 
    \end{align*}
    As $ \mathcal{N}_\rho$ is $\NIP$, the truth value of $P(a_{\mathfrak{N},i},b_\mathfrak{N})$ is eventually constant.
    
    We showed that the truth value of 
    $(P_\bullet(\mathbf{a}_i,\mathbf{b}) \wedge \delta(\mathbf{a}_i,\mathbf{b}))_i$ is eventually constant for all indiscernible sequences $(\mathbf{a}_i)_i$ and for all elements $\mathbf{b}$ in an elementary extension. Thus $P_\bullet(\mathbf{x},\mathbf{y}) \wedge \delta(\mathbf{x},\mathbf{y})$ is $\NIP$.
    \end{proof}
    
    \addtocounter{theorem}{-2}
     We also can deduce with the same proof the enriched version :
    \begin{corollary}[Enriched version]\label{CorollaryNIPTransferLexicoProduct}
           Let $\mathcal{M}$ be a  $\mathrm{L}_\mathcal{M}$-structure  and $\mathfrak{N}:= \{\mathcal{N}_a\}_{a\in \mathcal{M}}$ be a class of $\mathrm{L}_\mathfrak{N}$-structures which interprets uniformly a structure $\mathcal{K}$ with a non-empty set of surjective maps $\{\lambda_i: \mathcal{N}_a^{k_i} \rightarrow K\}_{a\in \mathcal{M},i\in I}$.
        Assume:
        \begin{itemize}
            \item for every sentence $\phi  \in \mathrm{L}_\mathfrak{N}$, the set $\{a\in \mathcal{M} \ \vert \ \mathcal{N}_a \models \phi \}$ is $\emptyset$-definable in $\mathcal{M}$,
            \item $\mathcal{M}$ is NIP,
            \item for all $\rho \in S_1^{Th(\mathcal{M})}$, $(\mathcal{N}_\rho, \mathcal{K}, \lambda_i)$ is NIP.
        \end{itemize}
   Then $(\mathcal{M}[\mathfrak{N}],\mathcal{K},\{\lambda_{\bullet,i}\}_{i\in I})$ is NIP.
    \end{corollary}
    \begin{proof}
    We only have to notice that a formula $\phi(\mathbf{x},\mathbf{y})$ (where $\mathbf{y}=\mathbf{y}^1,\ldots, \mathbf{y}^{\vert\mathbf{y} \vert} \in \mathcal{S}$) is still a disjunction of conjunction of formulas of the form:
            \begin{itemize}
                \item $\delta(\mathbf{x},\mathbf{y}) = \bigwedge_{l<\vert \mathbf{y}\vert} (\neg) v(\mathbf{x})=v(\mathbf{y}^l)$, a complete $s$-diagram in $\mathbf{x}$ over $\mathbf{y}^1,\ldots, \mathbf{y}^{\vert\mathbf{y} \vert}$,
                \item $\phi_{\bullet}(\mathbf{x},\mathbf{y})=\bigwedge (\neg)P_\bullet(\mathbf{x},\hat{\mathbf{y}})$, for finitely many $P\in \mathrm{L}_\mathfrak{N}$ and subtuples $\hat{\mathbf{y}}$ of $\mathbf{y}$,
                \item $\phi_{\mathcal{M}}(\mathbf{x},\mathbf{y})=\bigwedge(\neg) P(v(\mathbf{x}),v(\hat{\mathbf{y}}))$, for finitely many $P\in \mathrm{L}_\mathcal{M}$ and subtuples $\hat{\mathbf{y}}$ of $\mathbf{y}$.
            \end{itemize}
            
            For that, we used uniform interpretability of $\mathcal{K}$ in $\mathfrak{N}$ to replace any formula of the form $\phi_\mathcal{K}(\lambda(\mathbf{x}),\lambda'(\mathbf{x}))$ -- where $\phi_\mathcal{K}$ is a formula in $\mathrm{L}_\mathcal{K}$ and $\lambda,\lambda'$ are tuples of functions in $\{\lambda_i\}_{i\in I}$ --  by a formula of the form  $\phi_{\bullet}(\mathbf{x},\mathbf{y})$.
    \end{proof}
        \addtocounter{theorem}{1}

    \begin{theorem}\label{ThmBdnLexicoProd} 
    Consider the lexicographic product $\mathcal{S}$ of an $\mathrm{L}_\mathcal{M}$-structure $\mathcal{M}$ and a class of $\mathrm{L}_\mathfrak{N}$-structures $\mathfrak{N}:= \{\mathcal{N}_a\}_{a\in \mathcal{M}}$. Assume that for every sentence $\phi \in \mathrm{L}_\mathfrak{N}$, the set $\{a\in \mathcal{M} \ \vert \ \mathcal{N}_a \models \phi \}$ is $\emptyset$-definable in $\mathcal{M}$.
    Then:
        \[\bdn(\mathcal{S})=\sup\left(\bdn(\mathcal{M}), \bdn(\mathcal{N}_\rho), \rho \in \mathcal{S}_1^{Th(\mathcal{M})}\right).\]
        
    If $\mathcal{M}$ is $\vert \mathrm{L}_\mathfrak{N} \vert^+$-saturated, we may simply write:
        \[\bdn(\mathcal{S})=\sup\left(\bdn(\mathcal{M}), \bdn(\mathcal{N}_a), a \in \mathcal{M}\right).\]
    \end{theorem}
    
    A proof of the case $\mathcal{N}_a=\mathcal{N}$ for all $a\in \mathcal{M}$ can be found in \cite[Theorem A.4]{Tou20a}, but for completeness, we prefer to go through it again and stress the modification required.
    Notice that the formula holds also with Adler's convention: the proof gives $\bdn^\star(\mathcal{S})=\sup^\star\left(\bdn^\star(\mathcal{M}), \bdn^\star(\mathcal{N}_a), a \in \mathcal{M}\right)$.

    \begin{proof}
        
 The quantity $\max(\bdn(\mathcal{M}), \bdn(\mathcal{N}_a))$ is obviously a lower bound, as $\mathcal{M}[\mathfrak{N}]$ interprets on a unary set both $\mathcal{M}$ and $\mathcal{N}_\rho$ for any $\rho \in S_1^{\mathcal{M}}$.

        If $\mathcal{M}$ and all $\mathcal{N}_a$'s are finite, then so is  $\mathcal{S}$ and the theorem is trivially true. Assume that $\mathcal{M}$ or some $\mathcal{N}_a$ is infinite.
        As in the previous proofs, we consider the morleyrisation of $\mathrm{L}_\mathcal{M}$ and $\mathrm{L}_\mathfrak{N}$ and the natural expansion of $\mathcal{S}$, so that we may assume that $\mathcal{S}$ eliminates quantifiers by Theorem \ref{TheoremMeir2} in the language $\mathrm{L}_{\mathcal{M}[\mathfrak{N}]}$:
           \[\mathcal{S}:= (S,\mathrm{L}_{\bullet,\mathfrak{N}}) \cup (M,\mathrm{L}_\mathcal{M})\cup \{v:S\rightarrow M\}.\]
        We may assume that  $\mathcal{S}$ is enough saturated. Assume 
        \[\left\{\phi_i(\mathbf{x},\mathbf{y}_i), \left(\mathbf{c}_{i,j}=(a_{i,j},b_{i,j})\right)_{j<\omega}, k_i \right\}_{i<k}\]
        is an inp-pattern in $\mathcal{S}$ of depth $k> \max(\bdn(\mathcal{M}), \bdn(\mathcal{N}_a), \ a\in \mathcal{M}) \geq 1$ with $\vert \mathbf{x} \vert=1$.
        By Ramsey and compactness, we can assume that the sequence of parameters $(a_{i,j},b_{i,j})_{j<\omega}$ for $i<k$ are mutually indiscernible. By quantifier elimination and by elimination of the disjunctions (see for instance \cite{Tou20a}), we can also assume that for $i<k$ the formula $\phi_i(\mathbf{x},\mathbf{y}_i)$ (where $\mathbf{y}_i=\mathbf{y}_i^1,\ldots, \mathbf{y}_i^{\vert\mathbf{y}_i \vert}$) is a conjunction of formulas of the form:
            \begin{itemize}
                \item $\delta_i(\mathbf{x},\mathbf{y}_i) = \bigwedge_{l<\vert \mathbf{y}_i\vert} (\neg) v(\mathbf{x})=v(\mathbf{y}_i^l)$, a complete $s$-diagram in $\mathbf{x}$ over $\mathbf{y}_i^1,\ldots, \mathbf{y}_i^{\vert\mathbf{y}_i \vert}$,
                \item $\phi_{i,\bullet}(\mathbf{x},\mathbf{y}_i)=\bigwedge (\neg)P_\bullet(\mathbf{x},\hat{\mathbf{y}_i})$, for finitely many $P\in \mathrm{L}_\mathfrak{N}$ and subtuples $\hat{\mathbf{y}_i}$ of $\mathbf{y}_i$,
                \item $\phi_{i,\mathcal{M}}(\mathbf{x},\mathbf{y}_i)=\bigwedge(\neg) P(v(\mathbf{x}),v(\hat{\mathbf{y}_i}))$, for finitely many $P\in \mathrm{L}_\mathcal{M}$ and subtuples $\hat{\mathbf{y}_i}$ of $\mathbf{y}_i$.
            \end{itemize}
         Then, the main point of the proof is to remark that no inp-pattern can have a line "talking about" $\mathcal{M}$ and another one "talking about" an $\mathcal{N}_a$.
            
        \textbf{Case 1:} Assume that the complete $s$-diagram in $\mathbf{x}$ of the first line implies $v(\mathbf{x})=v(\mathbf{y}^0_{0,0})$. \\
        
        \begin{claim}
        By consistency of paths, and inconsistency of the lines, the same holds for every line  $i$: there are some $\mathbf{y}_i^{l_i}$, $l_i <\vert \mathbf{y}_i\vert $, such that $\phi_i(\mathbf{x},\mathbf{y}) \vdash v(\mathbf{x})=v(\mathbf{y}_{i}^{l_i})$.
        \end{claim}
        \begin{proof}
        Assume for example that line 1 implies $ v(\mathbf{x})\neq v(\mathbf{y}_1^l)$ for all $l<\vert \mathbf{y}_1\vert$. Since it also implies $ v(\mathbf{x})=v(\mathbf{y}^0_{0,0})$ , we have that $a_{0,0}^0 \neq a_{1,j}^l $ for all $j<\omega$ and $l<\vert \mathbf{y}_1\vert$. By consistency of paths, there exists $\textbf{d}_{j} \in \mathcal{S}$ such that \[\mathcal{S}\models \phi_0(\mathbf{d}_{j},\mathbf{c}_{0,0}) \wedge \phi_1(\mathbf{d}_{j},\mathbf{c}_{1,j}).\] 
        This in particular implies 
         \[ \mathcal{M} \models \phi_{1,\mathcal{M}}(a_{0,0}^0,a_{1,j}). \]
         By definition of $P_\bullet^{\mathcal{S}}$ and as $v(\textbf{d}_j)=a_{0,0}^0 \neq a_{1,j}^l$, this is equivalent to 
         \[ \mathcal{S} \models  \phi_1(\textbf{c}_{0,0}^0,\mathbf{c}_{1,j}).\]
         Then, the line $1$ would be realised by $\mathbf{c}_{0,0}^0$, contradiction. The same argument holds for any line $i>0$.
         \end{proof}
         Without loss of generality, we may assume that for all $i$, $\phi_i(\mathbf{x},\mathbf{y}) \vdash v(\mathbf{x})=v (\mathbf{y}_{i}^0)$ (we assume that $l_i=0$ for $i$'s). It follows from the claim and by consistency of paths that:
        \begin{claim}
            The parameters $\{a_{i,n}^0\}_{i<k,n<\omega}$ are all equal to some parameter $a\in \mathcal{M}$.
        \end{claim}

        It follows that there is a subconjunction $\phi_i'(\mathbf{x},\mathbf{y_i}')$ of $\phi_i(\mathbf{x},\mathbf{y}_i)$, where $\mathbf{y}_i'$ is a subtuple of $\mathbf{y}_i$ containing $\mathbf{y}_i^0$,  with complete (positive) diagram $\delta_i'(\mathbf{x},\mathbf{y}_i')=\bigwedge_{l<\vert \mathbf{y}_i'\vert} v(\mathbf{x})=v({\mathbf{y}_i'}^l)$ and which already forms an inconsistent line. 
        
        Indeed, if for some $l<\vert \mathbf{y}_i \vert$ $\phi_i(\mathbf{x},\mathbf{y_i}) \vdash  v(\mathbf{x}) \neq v(\mathbf{y}_i^l)$ then for all $P\in \mathrm{L}_\mathfrak{N}$ and subtuple $\hat{\mathbf{y}_i}$ of $\mathbf{y}_i$ containing $\mathbf{y}_i^l$: 
        \[\phi_i(\mathbf{x},\mathbf{y}_i) \vdash \neg P_\bullet(\mathbf{x},\hat{\mathbf{y}_{i}})\]

        and if $\phi_i(\mathbf{x},\mathbf{y_i}) \vdash  P(v(\mathbf{x}),v(\hat{\mathbf{y}_i}))$ (resp. $\phi_i(\mathbf{x},\mathbf{y_i}) \vdash \neg P(v(\mathbf{x}),v(\hat{\mathbf{y}_i}))$) for some subtuple $\hat{\mathbf{y}_i}$ of $\mathbf{y}_i$ containing $\mathbf{y}_i^l$, then by indiscernibility over $a$:
        \[\mathcal{M}\models\{P(a,\hat{a_{i,j}})\}_{j<\omega} \quad (\text{resp. }        \mathcal{M}\models\{\neg P(a,\hat{a_{i,j}})\}_{j<\omega} )\] 
        or equivalently, by definition of $P_\bullet^\mathcal{S}$:
        \[\mathcal{S} \models \{P(v(\mathbf{c}),v(\hat{\mathbf{c}_{i,j}}))\}_{j<\omega} \quad (\text{resp. }  \mathcal{S} \models \{\neg P(v(\mathbf{c}),v(\hat{\mathbf{c}_{i,j}}))\}_{j<\omega}\]
        whenever $v(\mathbf{c})=a$. 
        We deduce from it that we can remove from $\phi_i(\mathbf{x},\mathbf{y}_i)$ all atomic formulas involving $\mathbf{y}_i^l$, and we will still have an inconsistent line. 
        
        So we may assume that the $s$-diagram in $\mathbf{x}$  in any lines is positive, that is to say  \[\delta_i(\mathbf{x}) \equiv v(\mathbf{x})=v(\mathbf{y}_i^{1})=\cdots =v(\mathbf{y}_i^{l}).\]
        Thus, the inp-pattern translates to an inp-pattern of $\mathcal{N}_a$, of depth strictly bigger than $\bdn(\mathcal{N}_a)$, namely:
        \[\{\phi_{i,\mathfrak{N}}(x,y_i), (b_{i,j})_{j<\omega}\}_{i<k}\]
        where $\phi_{i,\mathfrak{N}}$ is the formula $\phi_{i,\bullet}$ see as a formula in $\mathrm{L}_{\mathfrak{N}}$ (by replacing occurrences of $P_\bullet$ by $P$).
        This is a contradiction with $\bdn(\mathcal{N}_a)<k$.\\
        \textbf{Case 2:} The complete $s$-diagram of any line is negative.  Then we define the following pattern in $\mathcal{M}$: 
        
        \[\mathcal{P}(x):= \left\{\bigwedge_{l<\vert \textbf{y}_i\vert} x \neq y_l \wedge \phi_{i,\mathcal{M}}(x,y_i), (a_{i,j})_{j<\omega}\right\}_{i<k}.\]
        We can show -- and this is a contradiction-- that it is an inp-pattern of depth strictly greater than $\bdn(\mathcal{M})$.
        Indeed, lines are inconsistent: if $a\in \mathcal{M}$ realises a line 
        \[\left \{\bigwedge_{l<\vert \textbf{y}_i\vert} x \neq a^l_{i,j} \wedge \phi_{i,\mathcal{M}}(x,a_{i,j}) \right\}_{j<\omega},\]
        then for any $b\in  \mathcal{N}_a$, $\mathbf{d}=(a,b)$ satisfies the corresponding line of the original pattern 
        \[\mathcal{S} \models \{\phi_{i}(\mathbf{d},\mathbf{c}_{i,j})\}_{j<\omega},\]
        this is a contradiction. It remains to show that paths are consistent: take $\mathbf{d}=(a,b)$ a realisation of the first column of our original inp-pattern:
        \[\mathcal{S} \models \{\phi_{i}(\mathbf{d},\mathbf{c}_{i,0})\}_{i<\kappa},\]
        then $a$ satisfies the first column of $\mathcal{P}(x)$:
        \[\left \{\bigwedge_{l<\vert \textbf{y}_i\vert} x \neq a^l_{i,0} \wedge \phi_{i,\mathcal{M}}(x,a_{i,0}) \right \}_{i<\kappa},\]
        This concludes our proof.
        
    \end{proof}
    
    Once again, using the same proof, we can state the enriched version:
    \addtocounter{theorem}{-1}
        \begin{theorem}[Enriched version]
        Let $\mathcal{M}$ be a  $\mathrm{L}_\mathcal{M}$-structure  and $\mathfrak{N}:= \{\mathcal{N}_a\}_{a\in \mathcal{M}}$ be a class of $\mathrm{L}_\mathfrak{N}$-structures which interprets uniformly a structure $\mathcal{K}$ with a non-empty set of surjective maps $\{\lambda_i: \mathcal{N}_a^{k_i} \rightarrow K\}_{a\in \mathcal{M},i\in I}$.
     Assume that for every sentence $\phi \in \mathrm{L}_\mathfrak{N}$, the set $\{a\in \mathcal{M} \ \vert \ \mathcal{N}_a \models \phi \}$ is $\emptyset$-definable in $\mathcal{M}$.
    Then:
        \[\bdn(\mathcal{M}[\mathfrak{N}],\mathcal{K},\{\lambda_{\bullet,i}\}_{i\in I})=\sup\left(\bdn(\mathcal{M}), \bdn(\mathcal{N}_\rho,\mathcal{K},\{\lambda_i\}_{i\in I}), \rho \in \mathcal{S}_1^{Th(\mathcal{M})}\right).\]

    \end{theorem}
    
    \begin{example}

\begin{itemize}
    \item     Consider in the language $\mathrm{L}:= \{\lambda_n\}_{n\in \mathbb{N}}$, and for each $n$ the $\mathrm{L}$-structure $\mathcal{M}_n$ where $\lambda_0,\cdots,\lambda_{n-1}$ interprets $n$ cross-cutting equivalence relation, and $\lambda_k$ for $k\geq n$ are empty predicates. Then we have: $\bdn(\mathbb{N}[\mathcal{M}_n])=\bdn(\prod_\mathcal{U} \mathcal{M}_n)= \aleph_0$ for any non-principal ultrafilter.
    
    \item Consider for all integer $n$ disjoint languages $\mathrm{L}_n$ with $n$ binary predicates and $\mathcal{M}_n$ a $\mathrm{L}_n$-structure with $n$ cross-cutting equivalence relations. Set $\mathrm{L}=\bigcup_n \mathrm{L}_n$. Then $\bdn(\mathbb{N}[\mathcal{M}_n]) = \aleph_{0,-}$.
    
\end{itemize}
    \end{example}
    
    \subsubsection{Example and digression on the burden}
    \begin{example}
Let $\mathrm{L}=\{R,B\}$ be the language with two binary predicates, and let $\mathcal{M}$ be a set with two cross-cutting equivalence relations with infinitely many infinite classes. 
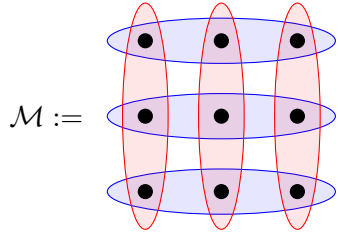
\begin{center} 
\begin{tikzpicture}\clip(-2,-0.7) rectangle (5,2.7);
\draw (-1.3,1) node {$\mathcal{M}:=$};

\draw[blue, fill=blue!10] (1,0) ellipse (1.5 and 0.3);
\draw[blue, fill=blue!10] (1,1) ellipse (1.5 and 0.3);
\draw[blue, fill=blue!10] (1,2) ellipse (1.5 and 0.3);

\fill[red,opacity=0.1] (0,1) ellipse (0.3 and 1.5);
\draw[red] (0,1) ellipse (0.3 and 1.5);
\fill[red,opacity=0.1] (1,1) ellipse (0.3 and 1.5);
\draw[red] (1,1) ellipse (0.3 and 1.5);
\fill[red,opacity=0.1] (2,1) ellipse (0.3 and 1.5);
\draw[red] (2,1) ellipse (0.3 and 1.5);

\fill (0,0) circle (0.1);
\fill (1,0) circle (0.1);
\fill (2,0) circle (0.1);
\fill (0,1) circle (0.1);
\fill (1,1) circle (0.1);
\fill (2,1) circle (0.1);
\fill (0,2) circle (0.1);
\fill (1,2) circle (0.1);
\fill (2,2) circle (0.1);
\end{tikzpicture}
\end{center}
One may show that $\bdn(\mathcal{M})=2$.

We consider the lexicographic product $\mathcal{M}[\mathcal{M}]$ in the one sorted language of Meir $\mathrm{L}_{\mathbb{U},s}$. We leave to the reader to describe $\mathcal{M}[\mathcal{M}]$ (where $R$ and $B$ are no longer equivalence relations). Then, one can verify that it is also of burden $2$. We let the following picture for the intuition:\\
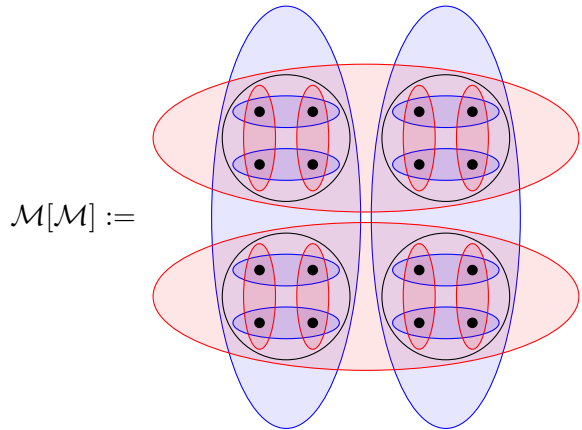
\begin{center}
    
\begin{tikzpicture}[scale=0.7]\clip(-5.5,-2) rectangle (10,6.4);
\draw (-3.5,2) node {$\mathcal{M}[\mathcal{M}]:=$};
\fill[blue,opacity=0.1] (0.5,2) ellipse (1.4 and 4);
\draw[blue] (0.5,2) ellipse (1.4 and 4);
\fill[blue,opacity=0.1] (3.5,2) ellipse (1.4 and 4);
\draw[blue] (3.5,2) ellipse (1.4 and 4);
\fill[red,opacity=0.1] (2,0.5) ellipse (4 and 1.4);
\draw[red] (2,0.5) ellipse (4 and 1.4);
\fill[red,opacity=0.1] (2,3.5) ellipse (4 and 1.4);
\draw[red] (2,3.5) ellipse (4 and 1.4);

\begin{scope}[shift={(0,0)}]
\fill[blue,opacity=0.1] (0.5,0) ellipse (1 and 0.3);
\draw[blue] (0.5,0) ellipse (1 and 0.3);
\fill[blue,opacity=0.1] (0.5,1) ellipse (1 and 0.3);
\draw[blue] (0.5,1) ellipse (1 and 0.3);

\fill[red,opacity=0.1] (0,0.5) ellipse (0.3 and 1);
\draw[red] (0,0.5) ellipse (0.3 and 1);
\fill[red,opacity=0.1] (1,0.5) ellipse (0.3 and 1);
\draw[red] (1,0.5) ellipse (0.3 and 1);
\draw (0.5,0.5) circle (1.2);
\fill (0,0) circle (0.1);
\fill (1,0) circle (0.1);
\fill (0,1) circle (0.1);
\fill (1,1) circle (0.1);
\end{scope}
\begin{scope}[shift={(0,3)}]
\fill[blue,opacity=0.1] (0.5,0) ellipse (1 and 0.3);
\draw[blue] (0.5,0) ellipse (1 and 0.3);
\fill[blue,opacity=0.1] (0.5,1) ellipse (1 and 0.3);
\draw[blue] (0.5,1) ellipse (1 and 0.3);

\fill[red,opacity=0.1] (0,0.5) ellipse (0.3 and 1);
\draw[red] (0,0.5) ellipse (0.3 and 1);
\fill[red,opacity=0.1] (1,0.5) ellipse (0.3 and 1);
\draw[red] (1,0.5) ellipse (0.3 and 1);
\draw (0.5,0.5) circle (1.2);
\fill (0,0) circle (0.1);
\fill (1,0) circle (0.1);
\fill (0,1) circle (0.1);
\fill (1,1) circle (0.1);
\end{scope}
\begin{scope}[shift={(3,0)}]
\fill[blue,opacity=0.1] (0.5,0) ellipse (1 and 0.3);
\draw[blue] (0.5,0) ellipse (1 and 0.3);
\fill[blue,opacity=0.1] (0.5,1) ellipse (1 and 0.3);
\draw[blue] (0.5,1) ellipse (1 and 0.3);

\fill[red,opacity=0.1] (0,0.5) ellipse (0.3 and 1);
\draw[red] (0,0.5) ellipse (0.3 and 1);
\fill[red,opacity=0.1] (1,0.5) ellipse (0.3 and 1);
\draw[red] (1,0.5) ellipse (0.3 and 1);
\draw (0.5,0.5) circle (1.2);
\fill (0,0) circle (0.1);
\fill (1,0) circle (0.1);
\fill (0,1) circle (0.1);
\fill (1,1) circle (0.1);
\end{scope}
\begin{scope}[shift={(3,3)}]
\fill[blue,opacity=0.1] (0.5,0) ellipse (1 and 0.3);
\draw[blue] (0.5,0) ellipse (1 and 0.3);
\fill[blue,opacity=0.1] (0.5,1) ellipse (1 and 0.3);
\draw[blue] (0.5,1) ellipse (1 and 0.3);

\fill[red,opacity=0.1] (0,0.5) ellipse (0.3 and 1);
\draw[red] (0,0.5) ellipse (0.3 and 1);
\fill[red,opacity=0.1] (1,0.5) ellipse (0.3 and 1);
\draw[red] (1,0.5) ellipse (0.3 and 1);
\draw (0.5,0.5) circle (1.2);
\fill (0,0) circle (0.1);
\fill (1,0) circle (0.1);
\fill (0,1) circle (0.1);
\fill (1,1) circle (0.1);
\end{scope}

\end{tikzpicture}

\end{center}

\end{example}

The fact that $\mathcal{M}$ and $\mathcal{M}[\mathcal{M}]$ have the same burden does not means that  the structures cannot be `distinguished' using inp-patterns. To see this, let us recall the following notion:

\begin{definition}
            Let $\mathcal{M}$ and $\mathcal{N}$ be two structures. We say that $\mathcal{N}$ is \emph{interpretable on a unary set} \index{Interpretability on a unary set} in $\mathcal{M}$ if there is a bijection $f:N \rightarrow D/\sim$ where $D$ is a unary definable set in $\mathcal{M}$, $\sim$ is a definable equivalence relation, and the pullback in $\mathcal{M}$ of any graph of function and relation of $\mathcal{N}$ is definable. The structures  $\mathcal{M}$ and $\mathcal{N}$ are said \emph{bi-interpretable on unary sets} if $\mathcal{N}$ is interpretable on a unary set in $\mathcal{M}$ and $\mathcal{M}$ is interpretable on a unary set in $\mathcal{N}$.
\end{definition}

The burden is an invariant of a class of structures bi-interpretable  on unary sets.
We saw with the previous example that the converse does not hold:
even though $\bdn(\mathcal{M}[\mathcal{M}])=\bdn(\mathcal{M})$, $\mathcal{M}[\mathcal{M}]$ is not interpretable in $\mathcal{M}$ on a unary set. This can be seen by analysing formulas and using the Morley rank. One can also use the notion of inp-patterns. Consider the inp-pattern of depth $2$ 
\[\{R(x,a_i),B(x,b_i)\}_{i<\omega}\]
where $\neg s(b_i,a_j)$, $\neg s(b_i,b_j)$ and $\neg s(a_i,a_j)$ for all $i<j\in \omega$. For each choice of a path $(i,j) \in \omega^2$, and realisation $d\models R(x,a_i)\wedge B(x,b_i)$, the type $\tp(d/a_i,b_j)$ is itself of burden $2$. Such pattern can't be found in $\mathcal{M}$. This motivates the following definition\footnote{Many thanks to A. Fornasiero for suggesting a similar definition to us and for the discussion which as inspired this digression.}:

\begin{definition}
    We denote the operation of concatenation of two sequences by $\upwedge$: for two disjoint sets of ordered indices $(I,<)$ and $(J,<)$, then $(a_i)_{(I,<)} \upwedge (a_i)_{(J,<)} = (a_i)_{(I\sqcup J,<)}$ where $i<j$ for all $i\in I$ and $j\in J$.  
    \begin{itemize}
        \item Any type $\rho(x)$ has inp-sequence $\emptyset$.
        \item A type $\rho(x)$ has inp-sequence $\alpha= (n) \upwedge \beta$ where $\beta$ is a sequence of integers (possibly empty), if there is an inp-pattern  $\{\phi_i(x,y_i), (a_{i,j})_{j<\omega}\}_{i<n}$ in $\rho$ of depth $n$ such that for all path, there is some realisation $d \models \{\phi_i(x,a_{if(i)})\}_{i<n}$ such that $\tp(d)$ has inp-sequence $\beta$.
    \end{itemize}
\end{definition}

The type $\{x=x\}$ has inp-sequence $(2,2)$ in $\mathcal{M}[\mathcal{M}]$ but only $(2)$ in  $\mathcal{M}$.

\begin{fact}\label{FactBdnInterpretUnarySet}
        Let $\mathcal{M}$ and $\mathcal{N}$ be two structures, and assume that $\mathcal{N}$ is interpretable on a unary set in $\mathcal{M}$, then $\bdn(\mathcal{N}) \leq \bdn(\mathcal{M})$. In particular, if $\mathcal{M}$ and $\mathcal{N}$ are bi-interpretable on unary sets, then there exists a type of inp-sequence $\alpha$ in $\mathcal{M}$ if and only if  there exists a type of inp-sequence $\alpha$ in $\mathcal{N}$.
\end{fact}

Note that if $\mathcal{M}$ is a model of $\aleph_0$-many cross-cutting equivalence relations in the language $\{E_n\}_{n\in \mathbb{N}}$, then $\mathcal{M}[\mathcal{M}]$ is interpretable in $\mathcal{M}$ on a unary set.
Indeed, first the equivalence relation $s$ is interpreted by $E_0$. 
\[\mathcal{M}[\mathcal{M}]\models s(x,y) \leftrightarrow \mathcal{M}\models E_0(x,y).\]
Then, we use two equivalence relations $E_{2n+1}$ and $E_{2n+2}$ in $\mathcal{M}$ to interpret $E_{n}$:
\[\mathcal{M}[\mathcal{M}]\models E_{n}(x,y) \ \leftrightarrow \mathcal{M}\models (E_{0}(x,y) \wedge E_{2n+1}(x,y) )\bigvee (\neg E_{0}(x,y) \wedge E_{2n+2}(x,y) ) .\]

\section{Valued Vector Spaces with Valuation-Preserving Scalar Multiplication}

In this section, we study scalar preserving valued vector spaces as Kuhlmann and Kuhlmann \cite{KK97}. We borrow from their work some notations and terminologies that we recall. Let $K$ be a field, $V$ a $K$-vector space and $(\Gamma,<,\infty)$ a totally ordered set with last element $\infty$.

A surjective map $\val: V \rightarrow \Gamma$ is called \textit{valuation} if for all $x,y\in V$:
\begin{itemize}
    \item  $\val(x)=\infty$ if and only if $x=0$,
    \item  $\val(x-y) \geq \min (\val(x),\val(y))$. 
\end{itemize}
We will also assume that the scalar multiplication  preserves the valuation: 
\begin{itemize}
    \item (trivial scalar action) $\forall k\in K^\star, \ x\in V \ \val(kx)=\val(x)$.
\end{itemize}
We only consider valued vector spaces with a trivial action of the scalar on the value set. We simply called them valued vector spaces.
For $\gamma \in \Gamma$, we define the $K$-vector spaces:
\begin{itemize}
    \item $V^\gamma:= \{x\in V; \val(x)\geq \gamma\}$,
    \item $V_\gamma:= \{x\in V; \val(x)> \gamma\}$,
    \item and $B(V,\gamma):=V^\gamma / V_\gamma $ with the convention that $B(V,\infty)= \{0\}$.
\end{itemize}
    By convention, we also set $V_\infty := V$. We denote by $\pi_{\delta}: V^\gamma \rightarrow B(V,\gamma)$ the natural projection. We consider the theory of vector spaces $V$ with a fix (resp. non fixed) field of scalars $K$. We will work in the language with fix base field 
    \[\mathrm{L}_K:= \{(V,+,0, \scal_k : V\rightarrow V, k\in K),(\Gamma,<,\infty), \val:V \rightarrow \Gamma\},\]
    where $\scal_k$ is interpreted by the multiplication by the scalar $k\in K$, and in the language of with variable base field:
    \[\mathrm{L}:= \{(V,+,0),(K,+,\cdot, 0, 1),(\Gamma,<,\infty), \scal : K \times V\rightarrow V, \val:V \rightarrow \Gamma\},\]
    where $K$ is a new sort.
    
    \subsection{Quantifier elimination relative to the Skeleton}\label{DefinitionSkeleton}
    Recall that we define the skeleton $S(V)$ of $V$ as the quotient of $ D(V):=\{ (\gamma,v) \in \Gamma \times V \ \vert \ \val(v) \geq \gamma \}$ by the definable equivalence relation:
    \[(\gamma,v) \simeq (\gamma',v') \text{ if and only if } \gamma=\gamma' \text{ and } \val(v-v')>\gamma. \]
    We define the natural projection by $\pi: V \rightarrow S(V), v\mapsto (\val(v),v+V_{\val(v)})$. 
    We denote by $\mathbf{0}=(\infty,0)$ the class of $0$.    We may define the following operation on $S(V):$

\[    \scal_\bullet (k,(\gamma,v+ V_\gamma)) := (\gamma,k\cdot v+ V_\gamma),\]

\[    (\gamma,v+ V_\gamma) +_\bullet (\gamma',v'+ V_{\gamma'}) := \begin{cases}
\mathbf{0} & \text{if } \gamma'\neq \gamma,\\
(\gamma,v + v' + V_\gamma) & \text{otherwise}.
\end{cases}\]
        We may also define the induced valuation map:
        \[\val_S:S(V) \rightarrow \Gamma, (\gamma,v)\mapsto \gamma.\] 
    
    These operation $+_\bullet$ and $\scal_\bullet(k, \cdot)$ on $S(V)$ are well defined i.e. they do not depend of the choice of representative, and are interpretable in the language $\mathrm{L}$ and $\mathrm{L}_{K}$. Then, equipped with these laws and the valuation, the new sort $S(V)$ is a lexicographic sum $\Gamma[B(V,\gamma)]_{\gamma\in\Gamma}$ of the vector spaces $B(V,\gamma)$ with fix base field $K$ (resp. variable base field $K$) with respect to the value set $(\Gamma,<)$. We will see (Theorem \ref{TheoremQuantifierEliminationRelativetoSkeleton} and Corollary \ref{CorollaryInducedStructureSkeleton}) that this is exactly the induced structure on $S(V)$.

    We will define later an another law, denoted by $\oplus$ with more intuitive property and which will also encode the valuation. We use for now this operation $+_\bullet$ as it `fits' better with the lexicographic product point of view.

     The skeleton structure takes a similar role to that of the leading term structure (or $\RV$-sort) in valued fields of equicharacteristic $0$. The $\RV$-sort of a valued field $\mathcal{K}=\{K,\gamma,k\}$ is traditionally defined by 
    \[\RV^\star:= K^\star/(1+\mathfrak{m})\]
    where $\mathfrak{m}$ is the maximal ideal of the valuation ring, but can also be described as follows:
    \[T_1:= \bigcup_{\gamma \in \Gamma} \mathcal{O}_\gamma / \mathfrak{m}_\gamma,\]
    where $\mathcal{O}_\gamma$ the set of elements of value $\geq \gamma$, $\mathfrak{m}_\gamma$ the set of element of value $>\gamma $ (see for example \cite{HHMD06} for a definition of $T_N$). As we said in the introduction, there is no analogous of a `residue field' in valued vector spaces. Instead, one must consider all the vector spaces $B(V,\gamma)$'s. This is of course due to the lack of structure (no internal multiplication in vector spaces).

    We recall informally the following theorem, and leave to the reader to check the precised statement in \cite{KK97}. 
    
    \begin{fact}[F-V. \& S. Kuhlmann] \label{TheoremKuhlmann}
        The theory of valued vector spaces with fixed base field eliminates quantifiers relative to \[\{\Gamma, (B(V,\gamma))_{\gamma \in \Gamma}\}.\]
    \end{fact}
      In \cite{KK97}, the skeleton is defined as this union of sorts
    $\{\Gamma, (B(V,\gamma))_{\gamma \in \Gamma}\}$.
    This is a stronger statement than Theorem \ref{TheoremQuantifierEliminationRelativetoSkeleton} below for fixed base field. However, as observed by F-V. and S. Kuhlmann, this theorem does not hold for variable based fields. As we will see (Theorem \ref{TheoremQuantifierEliminationRelativeVariableBaseField}), the value group $\Gamma$ needs to be indeed enriched with additional predicates.  Notice that in this presentation, one formally has to deal with $\vert \Gamma \vert $-many sorts, which number depends on your model. With our definition of skeleton, this syntactic problem do not arise. 

    We define the language $\mathrm{L}_S$ (resp. $\mathrm{L}_{S,K}$) extending $\mathrm{L}$ (resp. $\mathrm{L}_K$) with a sort for the skeleton:
    \[\mathrm{L}_S:= \mathrm{L} \cup \{(S(V),+_\bullet: S(V)^2 \rightarrow S(V),\mathbf{0}), \scal_{\bullet}: K\times  S(V) \rightarrow S(V), \val_S:S(V)\rightarrow \Gamma, \pi: V \rightarrow S(V)\},\]
        \[\mathrm{L}_{S,K}:= \mathrm{L}_K \cup \{(S(V),+_\bullet: S(V)^2 \rightarrow S(V),\mathbf{0}, \scal_{\bullet,k}:S(V) \rightarrow S(V) ), \val_S:S(V)\rightarrow \Gamma, \pi: V \rightarrow S(V)\}.\]
    
    By interpretation, every valued vector space admits an expansion to this language. Notice also that the function symbol for the valuation map $\val:V \rightarrow \Gamma$ is no more require to express definable set: it is indeed the composition of  $\pi:V \rightarrow S(V)$ and $\val_S: S(V) \rightarrow \Gamma$. 
    In the context of $\mathrm{L}_{S}$-structures (resp. $\mathrm{L}_{S,K}$-structure), $\mathcal{S}$-formula is a formula in the sub-language 
    \[\{(S(V),+_\bullet,\mathbf{0}), (\Gamma,<), (K,+,\cdot,0,1), \scal_{\bullet}, \val_S\}\]
    \[\text{(resp. } \{(S(V),+_\bullet, \scal_{\bullet,k},\mathbf{0}), (\Gamma,<), \val_S\}).\]
    
    Notice that it contains also formulas in the sort $\Gamma$ and $K$ (resp. in $\Gamma$).

  \begin{thmx}\label{TheoremQuantifierEliminationRelativetoSkeleton}
        The theory of valued $K$-vector spaces eliminates quantifiers relative to the skeleton in the language $\mathrm{L}_{S,K}$ and $\mathrm{L}_S$.  
    \end{thmx}
    
    \begin{proof}
     The proof is the same in both languages $\mathrm{L}_{S,K}$ and $\mathrm{L}_S$. For $k\in K$, $\textbf{s}\in S$ and $v\in V$, we simply write $k\cdot v$ instead of $\scal_{k}(v)$ or $\scal(k,v)$ and  $k\cdot \textbf{s}$ instead of $\scal_{\bullet,k}(\textbf{s})$ or $\scal_\bullet(k,\textbf{s})$.
    We proceed as in \cite[Proposition 4.1]{Fle11} by eliminating an existential quantifier `by hand'. First, let $a_i \in V$ be parameters and $\mathbf{z}_i$ be variables in $S(V)$ for all $i<n$. Then we first show that the formula
    \begin{align}
        \exists x\in V\ \left( \bigwedge_{i < n} \mathbf{z}_i=\pi(x-a_i)\right) \label{Formula1QuantifierExists}
    \end{align}
    is equivalent to a formula $\phi(\bar{\mathbf{z}},\bar{a})$ with no quantifiers.
    Let us denote by $B_i$ the ball $B_{>\val(z_i)}(z_i+a_i)$, where $\mathbf{z}_i= \pi(z_i)$. The formula then says that there is an $x$ in $\bigcap_{i<n} B_i $. Since the intersection of $n$ balls is non-empty if and only if the balls have pairwise non-empty intersection, we can assume that $n=2$. We can also assume $\val(z_0) \leq v(z_1)$.
    Then 
    \begin{align*}
        B_0 \cap B_1 \neq \emptyset & \Longleftrightarrow \ B_1 \subset B_0 \\
        &\Longleftrightarrow \ z_1+a_1 \in B_0 \\
        &\Longleftrightarrow \ \val(z_1+a_1-z_0-a_0) > \val(z_0).
    \end{align*}
    If $\val_S(\mathbf{z}_0) > \val(a_1-a_0)$,  then $\val(z_1+a_1-z_0-a_0) = \min(\val(z_0),\val(z_1),\val(a_1-a_0) ) = \val(a_1-a_0)$ and there is no solution.\\ 
    If $\val_S(\mathbf{z}_1)>\val_S(\mathbf{z}_0)$ and $ \val_S(\mathbf{z}_0) \leq  \val(a_1-a_0)$, then have 
    \[\val(z_1+a_1-z_0-a_0) > \val(z_0) \ \Longleftrightarrow \ \mathbf{z}_0+_\bullet \pi(a_1-a_2) = 0\cdot \mathbf{z}_0.\] 
      If $\val_S(\mathbf{z}_1) = \val_S(\mathbf{z}_0) \leq  \val(a_1-a_0)$, then we have 
      \[\val(z_1+a_1-z_0-a_0) > \val(z_0)  \ \Longleftrightarrow \  (\mathbf{z}_0+_\bullet-\mathbf{z}_1)+_\bullet \pi(a_1-a_2) = 0\cdot \mathbf{z}_0.\]

    Thus, the formula (\ref{Formula1QuantifierExists}) is given by a quantifier-free formula.
    Now, consider a formula  $\phi(\bar{u},\bar{\mathbf{w}})$ of the form:
        \begin{align}
            \exists x\in V \ \left( \psi(\pi(t_0(x,\bar{u})),\dots, \pi(t_{k-1}(x,\bar{u})),\bar{\mathbf{w}}) \right),
        \end{align}
    where $\psi(\mathbf{x}_1, \dots, \mathbf{x}_{k-1},\bar{\mathbf{w}})$ is an $S(V)$-formula and where $t_i(x,\bar{u})$ for $i<k$ are $V$-sorted linear terms $\bar{\lambda}_i\cdot \bar{u} + \lambda_i\cdot x $, where $\bar{\lambda}_i,\lambda_i \in K^{\vert \bar{u} \vert+1}$,$\bar{u} \in V^{\vert \bar{u} \vert}$. We can assume that $\lambda_i=1$ as we can write $\pi\left(\bar{\lambda}_i\cdot \bar{u} + \lambda_i\cdot x\right) = \lambda_i\cdot \pi \left( \lambda_i^{-1}\cdot \bar{\lambda}_i\cdot \bar{u} + x \right)$.
    Then $\phi(\bar{u},\bar{\mathbf{w}})$ is equivalent to
    \begin{align}
        \exists \mathbf{z}_0,\dots,\mathbf{z}_{k-1} \in S(V) \ \left( \psi(\mathbf{z}_0,\dots,\mathbf{z}_{k-1},\bar{\mathbf{w}}) \wedge \exists x\  \in V \  \bigwedge_{i\leq n} \mathbf{z_i}= \pi(t_i(x,\bar{u})) \right).
    \end{align}
    
    By the previous case, we may eliminate the existential quantifier  $\exists x \in V$, which give us a formula with only quantifiers in $S(V)$. This concludes the proof of quantifier elimination. 
    
    \end{proof}


    It follows that all formulas are equivalent in $\mathrm{L}_S$ (resp. $\mathrm{L}_{S,K}$) to a formula of the form:
 \[\phi_{S}\left(\pi(t_0(x)),\ldots,\pi(t_{n-1}(x))\right)\]
            where $\phi_S(x_{S,0},\ldots, x_{S,n-1})$ is an $S$-sorted formula and $t_{i}(x)$ for $i<n$ are vector-space-terms in $\mathrm{L}$ (resp. $\mathrm{L}_K$). 
    In particular, notice that a formula of the form 
    $t(x,\bar{v},\bar{k})=0$ (where $t(x,\bar{v},\bar{k})$ is a term) is equivalent to $\val_S(\pi(t(x,\bar{v},\bar{k})))= \infty$.
    
    Notice also that terms in $\mathcal{S}$ can be rather complicated, since the law $+_\bullet$ on $S(V)$ is not associative. There are certain interactions between $\Gamma$ and the $B(V,\gamma)'s$: for every integer $n$, one may indeed define in $\mathrm{L}$ the set of values $\{\gamma \in \Gamma \ \vert \  \dim_K(B(V,\Gamma))  = n\}$. One can also, for instance, uniformly enrich the $B(V,\gamma)$'s with an equivalence relation $E$. This will also enriched the induced structure in $\Gamma$ with predicates such as $\{\gamma \in \Gamma \ \vert \ E\text{ has $n$ classes in } B(V,\Gamma)\}$.
    \begin{corollary}\label{CorollaryInducedStructureSkeleton}
         \begin{itemize}
             \item The structure on the skeleton $S(V)$ induced by the valued vector space in the language $\mathrm{L}_{S,K}$ is given by
             \[\left\{(S(V),+_\bullet,\mathbf{0}),(\Gamma,<), (\scal_{\bullet,k}: S(V) \rightarrow S(V))_{k\in K}, \val_S:S(V) \rightarrow \Gamma, \pi: V \rightarrow S(V) \right\}. \]
         In other terms, the skeleton is a pure lexicographic sum $\Gamma[B(V,\gamma)]$ of one-sorted $K$-vector spaces $( B(V,\gamma),+,(\scal_k)_{k\in K},0)$ with respect to the order set $(\Gamma,<)$.
            \item The structure on the skeleton $S(V)$ induced by the valued vector space in the language $\mathrm{L}_{S}$ is given by
             \[ \left \{ (S(V),+_\bullet,\mathbf{0}),(\Gamma,<) , \scal_{\bullet}:K\times S(V) \rightarrow S(V), \val_S:S(V) \rightarrow \Gamma, \pi: V \rightarrow S(V) \right \}.\]
        In other terms, it is the lexicographic sum $\Gamma[B(V,\gamma)]$ of the groups $( B(V,\gamma),+,0)$ with respect to the order set $(\Gamma,<)$, enriched with the (uniform) action of $K$ on the ribs $\scal_{\bullet}:K\times S(V) \rightarrow S(V)$.
        \end{itemize}
         

    \end{corollary}

    This is a direct corollary of quantifier elimination and closedness of the sort $S(V)$. Remark that in the language $\mathrm{L}_{S,K}$, the ribs $\{B(V,\gamma)\}_{ \gamma \in \Gamma}$, can have different cardinality. In the language $\mathrm{L}_{S}$, they all must be of size at least $\vert K\vert$; and of size exactly $\vert K\vert$ if they have finite dimension.

    One can ask when a given lexicographic sum of vector spaces with respect to a linear order is the skeleton of a valued vector space. As noticed by Kuhlmann and Kuhlmann, this always holds. We recall from \cite{Kuh95} the definition of Hahn product and Hahn sum:

    \begin{definition}
         Let $\Gamma$ be any ordered set and  $(B_\gamma)_{\gamma \in \Gamma}$ be a collection of vector spaces.

        \begin{itemize}
            \item We denote by $\prod_{\gamma \in \Gamma}B_\gamma$ the usual product of vector spaces.
            \item We denote by $\bigoplus_{\gamma \in \Gamma}B_\gamma$ the sub-$K$-vector space of element $s$ of finite support, where the support is given by:
        \[\supp(s):=\{ \gamma \in \Gamma \ \vert \ s(\gamma)\neq 0 \}. \]
        \item We define the canonical valuation
\[        \begin{array}{rrcl}
             \val:& \bigoplus_{\gamma \in \Gamma}B_\gamma & \longrightarrow & \Gamma   \\
             & s & \longmapsto & \min \supp(s),
        \end{array}\]
        with the convention that $\val(0)= \infty$.
        \item The valued vector space $(\bigoplus_{\gamma \in \Gamma}B_\gamma, \Gamma, \val)$ is called the Hahn sum of $(B_\gamma)_{\gamma \in \Gamma}$ and is denoted by $\coprod_{\gamma\in \Gamma} B_\gamma $.

        \item Similarly, one can define the Hahn product, with base set the sub-$K$-vector space of $\prod_{\gamma \in \Gamma}B_\gamma$ of elements of well-ordered support, denoted by $\mathbf{H}_{\gamma\in \Gamma}B_\gamma$. 
        
        \end{itemize}  
    \end{definition}
    
    \begin{fact}
        Let $\Gamma$ be an ordered set and  $(B_\gamma)_{\gamma \in \Gamma}$ be a collection of vector spaces.
    The Hahn product $\hprod_{\gamma \in \Gamma}B_\gamma$ and the Hahn sum $\coprod_{\gamma \in \Gamma}B_\gamma$ admit for skeleton the lexicographic product $\Gamma[B_\gamma]_{\gamma}$.
    \end{fact}

    \subsection{Transfer principles relative to the Skeleton}

    \begin{proposition} \label{PropositionRelativeModelCompletnessSkeleton}
        (Relative Model Completness)
        The theory of valued vector spaces is model complete relative to the skeleton in the language $\mathrm{L}_{S,K}$ and $\mathrm{L}_S$. More precisely, given two models $\mathcal{V} \subseteq \mathcal{W}$, then $\mathcal{V} \preceq \mathcal{W}$  if and only if $S(V) \preceq S(W)$.
    \end{proposition}
    \begin{proof}
         Again, the proof is the same in both languages $\mathrm{L}_{S,K}$ and $\mathrm{L}_S$. It can be deduced directly from Theorem \ref{TheoremQuantifierEliminationRelativetoSkeleton} and Remark \ref{RemarkModelCompletnessRQECLoseSort}.
    \end{proof}
    
    \begin{proposition} \label{PropositionNIPTransferSkeleton}
        (NIP transfer)
        A valued vector space in the language $\mathrm{L}_{S,K}$ (resp. $\mathrm{L}_S$) is NIP if and only if its skeleton is NIP.
    \end{proposition}
    \begin{proof}
    This is deduced from Remark \ref{RemarkNIPRQECLoseSort} and the fact that $S(V)$ is a closed sort.
    \end{proof}

    Following the same proof as in \cite{Tou20a}, we have the transfer principle:
    \begin{proposition} \label{PropositionBdnValuedVectorSpaceRelativeSkeleton}
        Consider a valued vector space $V$ in the language $\mathrm{L}_{S,K}$ (resp. $\mathrm{L}_S$). Then:
         \[\bdn(V)= \bdn(S(V)).\]
    \end{proposition}

    To emphasise the analogy between the $\RV$-structure of valued fields and the skeleton of valued vector spaces, we define a new language for the skeleton of valued vector spaces, with no sort for the value sort. We extend the usual sum in the ribs $B(V,\gamma)$'s to a new operation $\oplus$ in $S(V)$.
    \begin{definition}
        Let $\mathbf{v}=(\gamma,v),\mathbf{v}'=(\gamma',v')$ be two non-trivial elements of the skeleton $S(V)$. We set
        \[\mathbf{v} \oplus \mathbf{v}' := \begin{cases}
        \mathbf{v} &\text{ if } \gamma<\gamma',\\
        \mathbf{v}' & \text{ if } \gamma'<\gamma,\\
        (\gamma,v+v') =(\gamma,v)+_\bullet(\gamma',v') &\text{ if } \gamma=\gamma'.
        
        \end{cases}\]
        We set also $\mathbf{0}\oplus \mathbf{v} = \mathbf{v} \oplus \mathbf{0}= \mathbf{v}$ and $\mathbf{0} \oplus \mathbf{0}=\mathbf{0}$.
        We say that the sum is \textit{not well defined} if $(\gamma=\gamma')$ and $v+v'=0$, and that it is \textit{well defined} otherwise. We denote by $\WD(x,y)$ the binary predicate for well defined sum $x\oplus y$. 
    \end{definition}
    In the next proof, we will use the following:
    \begin{observation}
        If $a,b\in V\setminus \{0\}$, then $\WD(\pi(a),\pi(b))$ holds if and only if $\pi(a+b)=\pi(a) \oplus \pi(b)$.
    \end{observation}

    Notice again that this new operation $\oplus$ is not an addition (for instance, it has absorbing element). Also notice that $\oplus$ and $\WD$ are definable in the language $\mathrm{L}_{S}$ and $\mathrm{L}_{S,K}$.

    \begin{remark}
        \begin{enumerate}
             \item       \[\{(S(V),+_\bullet,\mathbf{0}, (\scal_{\bullet,k})_{k\in K}), (\Gamma,<), \val_S, \pi\} \]
        
        and 
        
        \[\{S(V),\oplus,\mathbf{0},(\scal_{\bullet,k})_{k\in K} \}\]
         are bi-interpretable.  
 \item
        \[\{(S(V),+_\bullet,\mathbf{0}), (K,+,\cdot,0,1), (\Gamma,<) ,\scal_{\bullet}, \val_S, \pi\} \]
        
        and 
        
        \[\{(S(V),\oplus,\mathbf{0}), (K,+,\cdot,0,1) ,\scal_{\bullet},  \pi\}\]
         are bi-interpretable.
         \end{enumerate}
    \end{remark}
    \begin{proof}
        This is the analogous of \cite[Proposition 2.8]{Fle11} for valued vector spaces.
        It remains to show that $+,\val_S, \Gamma$ and $<$ can be interpreted with $\oplus$. We have indeed for $a,b$ and $c$ in $S(V)$: 
        
    \begin{itemize}
        \item $\val_S(a)=\val_S(b)$ if and only if $a \oplus b \neq a$ and $a \oplus b \neq b$.
        \item $\val_S(a)<\val_S(b)$ if and only if $a \oplus b = a$.
        \item $a+_\bullet b=c$ if and only if $\val_S(a)=\val_S(b)$ and $a \oplus b = c$ or $c=\mathbf{0}$.
    \end{itemize}
    
    \end{proof}
    
    Now, one can directly transpose the proof of \cite[Theorem  3.2]{Tou20a} in the setting of valued vector spaces. In many aspects, the proof become simpler and hopefully, easier to follow. Even if the main ideas remain, we rewrite the proof for this purpose of clarification.
    
    We will need the notion of pseudo-convergent sequences for valued groups:
    \begin{definition}
        Let $I$ be a sequence without end point. We say that a sequence $(a_i)_{i\in I}$ in $V$ is \textit{eventually pseudo-Cauchy} or simply \textit{pseudo-Cauchy} if for there is $i_0 \in I$ such that for all $i>j>k>i_0$, we have 
        \[\val(a_i-a_j)>\val(a_j-a_k).\]
        An element $a$ such that for all $i>j>i_0$ 
        \[\val(a-a_j)=\val(a_i-a_j)\]
        is called pseudo-limit of $(a_i)_{i\in I}$. We write \[(a_i)_{i\in I} \Rightarrow a.\]
    \end{definition}

    \begin{proof}
    We denote by $\bar{\mathbb{Z}}$ the set of natural numbers with extremal points $ \mathbb{Z}\cup \lbrace \pm \infty\rbrace$. Let $M$ be a cardinal (possibly infinite) bigger than 1. Let 
    \[\lbrace \tilde{\phi}_i(x,y_i),(c_{i,j})_{j \in \bar{\mathbb{Z}}},k_i \rbrace_{i < M}\] be an inp-pattern in $V$ of depth $M \geq 2$ with $\vert x \vert = 1$, where $c_{i,j}=a_{i,j}\mathbf{b}_{i,j} \in V^{l_1} \times S(V)^{l_2}$ (notice that the set of indices is $\bar{\mathbb{Z}}$: we will make use of one of the extreme elements $\lbrace a_{i,-\infty},a_{i,+\infty}\rbrace$ later). We have to find an inp-pattern of depth $M$ in $S(V)$. Without loss of generality, we take $(c_{i,j})_{i,j}$ mutually indiscernible. By Theorem \ref{TheoremKuhlmann} and mutual indiscernibility, we can assume the formulas $\tilde{\phi}_i$, $i<M$, are of the form 
\[\tilde{\phi}_i(x,c_{i,j}) = \phi_i(\pi(x-a_{i,j;1}), \ldots , \pi(x-a_{i,j;l_i});\mathbf{b}_{i,j}),\]
for some integer $l_i$ and where $\phi_i$  is an $S(V)$-formulas. 

Let $d \models \lbrace \phi_i(\pi(x-a_{i,0;1}), \ldots , \pi(x-a_{i,0;l_i});\mathbf{b}_{i,0}) \rbrace_{i <M}$ be a solution of the first column. Before we give a general idea of the proof, let us reduce to the case where only one term $\pi(x-a_{i,j})$ occurs in the formula $\tilde{\phi}_i$.

\begin{claim}\label{claim1}
We may assume that for all $i<M$, $\tilde{\phi}_i(x,c_{i,j})$ is of the form $\phi_i(\pi(x-a_{i,j;1});\mathbf{b}_{i,j})$, i.e. $\vert a_{i,j} \vert = l_i=1$.
\end{claim}
\begin{proof}
We will first replace the formula $\tilde{\phi}_0(x,c_{0,j})$ by a new one with an extra parameter. 

One can see that at least one of the following two cases occurs
\begin{enumerate}

    \item $\WD \left(\pi(d-a_{0,0;1}), \pi(a_{0,0;1}-a_{0,0;2}) \right)$ or
    \item $\WD \left(\pi(d-a_{0,0;2}),\pi(a_{0,0;2}-a_{0,0;1}) \right).$
\end{enumerate}
 
According to the case, we respectively define a new formula $\psi_0(x,c_{0,j}\upwedge \pi(a_{0,j;2}-a_{0,j;1}) )$ as follow:
\begin{enumerate}
    \item \begin{multline*}\phi_0(\pi(x-a_{0,j;1}),\pi(x-a_{0,j;1})\oplus \pi(a_{0,j;1}-a_{0,j;2}),\pi(x-a_{0,j;3}), \ldots , \\ \pi(x-a_{0,j;l_0});\mathbf{b}_{0,j}) \wedge \WD(\pi(x-a_{0,j;1}), \pi(a_{0,j;1}-a_{0,j;2})),
\end{multline*}

  \item \begin{multline*}\phi_0(\pi(x-a_{0,j;2})\oplus\pi(a_{0,j;2}-a_{0,j;1}),\pi(x-a_{0,j;2}), \ldots , \pi(x-a_{0,j;l_0});\mathbf{b}_{0,j}) \\ \wedge \WD(\pi(x-a_{0,j;2}), \pi(a_{0,j;2}-a_{0,j;1})).
\end{multline*}
\end{enumerate}
We will prove that the pattern where $\tilde{\phi}_0$ is replaced by $\psi_0$:
$$\left\lbrace \psi_0(x,y_0 \upwedge z),(c_{0,j}\upwedge \pi(a_{0,j;2}-a_{0,j;1}))_{j \in \bar{\mathbb{Z}}}, k_0 \right \rbrace \cup \left \lbrace \tilde{\phi}_i(x,y_i),(c_{i,j})_{j \in \bar{\mathbb{Z}}},k_i
\right\rbrace_{1 \leq i <M}$$
is also an inp-pattern. First note that we have added $\pi(a_{0,j;2}-a_{0,j;1})$ to the parameters $\mathbf{b}_{0,j}$, and it still forms a mutually indiscernible array. Clearly, $d$ is still a realisation of the first column: 
$$d\models \left\lbrace \psi_0(x,c_{0,0}\upwedge \pi(a_{0,0;2}-a_{0,0;1}) ) \right\rbrace \cup \left\lbrace \tilde{\phi}_i(x,c_{i,0})\ \vert \ 1 \leq i <M  \right\rbrace.$$ By mutual indiscernibility of the parameters, every path is consistent. Since $\psi_0(x,c_{0,j}\upwedge \pi(a_{0,0;2}-a_{0,0;1})) \vdash \subseteq \tilde{\phi}_0(x,c_{0,j})$, inconsistency of the first row is also clear. By induction on the number of terms $l_0$, we may assume that $\tilde{\phi}_0(x,c_{0,j})$ has only one term of the form $\pi(x-a_{0,j})$. We can do the same for all formulas $\tilde{\phi}_i$, $0 < i < M$.
\end{proof}

If the array $(a_{i,j})_{i<M, j<\omega}$ is constant equal to some $a\in V$, then the following pattern:
 \[\lbrace \phi_i(\mathbf{x},z_i),(\mathbf{b}_{i,j})_{j \in \bar{\mathbb{Z}}},k_i \rbrace_{i <M},\]
 where $\mathbf{x}$ is a variable in $S$ is an inp-pattern of depth $M$ in $S(V)$. Indeed, consistency of the path is clear. If a row is satisfied by some $\mathbf{d}\in S(V)$  , any $d\in V$ such that $\pi(d-a)=\mathbf{d}$ will satisfy the corresponding row of the initial inp-pattern, which is absurd. Hence, the rows are inconsistent.  

The idea of the proof is to reduce the general case (where the $a_{i,j}$'s are distinct) to this trivial case by the same method as above: removing the parameters $a_{i,j} \in V$ and adding new parameters from $S(V)$ to the tuple $\mathbf{b}_{i,j}$ and by replacing the formula $\tilde{\phi}_i$ a  formula of the form 
\[\psi_i(x,b_{i,j}) \vee \WD(\pi(x-a),\pi(a-a_{i,j})).\]

The main challenge is to find a suitable $a\in V$ as a `center' of our pattern. Recall that $d \models \lbrace \phi_i(\pi(x-a_{i,0});b_{i,0}) \rbrace_{i <M}$ is a solution to the first column.

\begin{claim}\label{claim2}
For all $j < \omega $, and $i,k < M$ with $k\neq i$, we have ${\val(d-a_{i,j}) \leq \val(d-a_{k,0}) }. $
\end{claim}
 \begin{proof}
 	Assume not: for some $j < \omega$, and $i,k < M$ with $k\neq i$:
 	$$\val(d-a_{i,j}) > \val(d-a_{k,0})$$
Then, $\pi(a_{i,j}-a_{k,0}) = \pi(d-a_{k,0})$. By mutual indiscernibility, we have
$$a_{i,j} \models \lbrace \phi_k(\pi(x-a_{k,l});\mathbf{b}_{k,l}) \rbrace_{l< \omega}.$$
This contradicts inconsistency of the row $k$.
\end{proof}

In particular, for all $i,k < M$, we have $  \val(d-a_{k,0}) = \val(d-a_{i,0})$. Let us denote  by $ {\gamma:=\val(d-a_{0,0}) }$ this value. By definition, we have the following for all $i,k<M$:
\begin{equation}
 \ \val(a_{i,0}-a_{k,0}) \geq \min \lbrace \val(d-a_{i,0}), \val(d-a_{k,0}) \rbrace \geq \gamma.\tag{$\star$}
\end{equation}

The following lemma give us a correct choice for a `center' $a$.

\begin{lemma}\label{mainprop}
We may assume that there is $i <M $ such that for all $k < M$, the following holds:
\[\WD(\pi(d-a_{i,\infty}),\pi(a_{i,\infty}-a_{k,0})).\]
\end{lemma}

We may conclude from this lemma.
Indeed, assume $i=0$ satisfies the conclusion above. For every $k< M$, we have  
\[ \WD(\pi(d-a_{0,\infty}),\pi(a_{0,\infty}-a_{k,0})).\]
Set $\tilde{b}_{i,j}:=b_{i,j} \upwedge \pi(a_{0,\infty}-a_{i,j}) $ for $i < M, j<\omega$ and 
$$\psi_i(\tilde{x},\tilde{b}_{i,j}):=\phi_i\left(
\tilde{x}\oplus\pi(a_{0,\infty}-a_{i,j}));b_{i,j}\right) \wedge \WD(\tilde{x},\pi(a_{0,\infty}-a_{i,j}))$$
where $\tilde{x}$ is a variable in $S(V)$.

These formulas form an inp-pattern. Indeed, clearly, $\pi(d-a_{0,\infty}) \models \lbrace \psi_i(\tilde{x},\tilde{b}_{i,0})\rbrace_{i < M}$. By mutual indiscernibility of $(\tilde{b}_{i,j})_{i<M,j < \omega}$, every path is consistent. It remains to show that, for every $i < M$, $\lbrace \psi_i(\tilde{x},\tilde{b}_{i,j}) \rbrace_{j<\omega}$ is inconsistent.
Assume there is $\alpha^\star  \models \lbrace \psi_i(\tilde{x},\tilde{b}_{i,j}) \rbrace_{j<\omega}$ for some $i< M$, and let $d^\star$ be such that $\pi(d^\star - a_{0,\infty})= \alpha^\star$. Then, since $\WD(\alpha^\star,\pi(a_{0,\infty}-a_{i,j}))$ holds for every $j<\omega$, $d^\star$ satisfies ${\lbrace \phi_i(\pi(x-a_{i,j}), b_{i,j})\rbrace_{j<\omega}}$, which is a contradiction. All rows are inconsistent. This will concludes our proof. It thus remains to prove Lemma \ref{mainprop}. 
    \end{proof}

\begin{proof}[proof of Lemma \ref{mainprop}]
It is enough to find $i<M$ such that the following holds:
\begin{itemize}
    \item $\gamma  \leq \val(d-a_{i,\infty})$, or
    \item for $\vert M \vert$-many $k \in M$, $ \gamma \leq   \val(a_{i,\infty}-a_{k,0})$. 
\end{itemize}
Indeed in both case, we have for $\vert M \vert$-many $k$ that \[ \gamma = val(d-a_{k,0}) =\min \lbrace \val(d-a_{i,\infty}), \val(a_{i,\infty}-a_{k,0})\rbrace \]
which in particular gives 
\[\WD(\pi(d-a_{i,\infty}),\pi(a_{i,\infty}-a_{k,0})).\]

If $M$ is infinite , we might have to remove the lines $k$ where the above sum is not well defined, which we can do without lost. The first case will correspond to Case A, the second to Case B. \\

	\textbf{Case A }: There are $0 \leq i, k < M$ with $i \neq k$ such that $\val(a_{i,j}-a_{k,l})$ is constant for all $j,l \in \mathbb{N}\cup\{+\infty\}$, equal to some $\epsilon$.
	Note that $(\star)$ gives $\epsilon \geq \gamma$.
	\\
	
\begin{center}
\begin{minipage}{0.70\linewidth}

	\begin{tikzpicture}[line cap=round,line join=round,>=triangle 45,x=1.0cm,y=1.0cm,scale = 0.5]
\clip(-2.226855653004253,-4.195898595071177) rectangle (17.415560940241033,5.607988346530954);
\draw [line width=1.2pt] (2.,5.)-- (6.,-3.);
\draw [line width=1.2pt] (6.,-3.)-- (10.,5.);
\draw [line width=1.2pt] (6.,-3.)-- (7.,-5.);
\draw [line width=1.2pt] (2.2477554844708085,4.504489031058383)-- (2.5,5.);
\draw [line width=1.2pt] (2.5,4.)-- (3.,5.);
\draw [line width=1.2pt] (2.741734673539459,3.5165306529210825)-- (3.5,5.);
\draw [line width=1.2pt] (3.498797557358775,2.0024048852824508)-- (5.,5.);
\draw [line width=1.2pt] (3.748010012984067,1.5039799740318665)-- (5.5,5.);
\draw [line width=1.2pt] (6.,5.)-- (4.,1.);
\draw [line width=1.2pt] (6.5,5.)-- (4.244320737961419,0.509786121281313);
\draw [line width=1.2pt] (9.,3.)-- (8.,5.);
\draw [line width=1.2pt] (8.244978810717003,4.510042378565994)-- (8.5,5.);
\draw [line width=1.2pt] (8.5,4.)-- (9.,5.);
\draw [line width=1.2pt] (8.748951172574003,3.502097654851993)-- (9.5,5.);
\draw [line width=1.2pt] (3.,3.)-- (4.,5.);
\draw [line width=1.2pt] (7.,5.)-- (4.5,0.);
\draw [line width=0.8pt,dotted] (3.,-3.)-- (10.,-3.);
\draw (7.9061370339873625,-2.966082388615786) node[anchor=north west] {$\val(a_{i,j}-a_{k,l})= \epsilon \geq \gamma $};
\begin{scriptsize}
\draw [fill=black] (2.,5.) circle (1.5pt);
\draw [fill=black] (6.,-3.) circle (1.0pt);
\draw [fill=black] (10.,5.) circle (1.5pt);
\draw[color=black] (10.521661642082634,5.174954471018492) node {$a_{k,\infty}$};
\draw [fill=black] (7.,-5.) circle (1.0pt);
\draw [fill=black] (2.2477554844708085,4.504489031058383) circle (1.0pt);
\draw [fill=black] (2.5,5.) circle (1.5pt);
\draw[color=black] (2.579820365184078,5.30959718105949) node {$a_{i,0}$};
\draw [fill=black] (2.5,4.) circle (1.0pt);
\draw [fill=black] (3.,5.) circle (1.5pt);
\draw [fill=black] (2.741734673539459,3.5165306529210825) circle (1.0pt);
\draw [fill=black] (3.5,5.) circle (1.5pt);
\draw [fill=black] (3.498797557358775,2.0024048852824508) circle (1.0pt);
\draw [fill=black] (5.,5.) circle (1.5pt);
\draw [fill=black] (3.748010012984067,1.5039799740318665) circle (1.0pt);
\draw [fill=black] (5.5,5.) circle (1.5pt);
\draw [fill=black] (6.,5.) circle (1.5pt);
\draw [fill=black] (4.,1.) circle (1.0pt);
\draw [fill=black] (6.5,5.) circle (1.5pt);
\draw [fill=black] (4.244320737961419,0.509786121281313) circle (1.0pt);
\draw [fill=black] (9.,3.) circle (1.0pt);
\draw [fill=black] (8.,5.) circle (1.5pt);
\draw [fill=black] (8.244978810717003,4.510042378565994) circle (1.0pt);
\draw [fill=black] (8.5,5.) circle (1.5pt);
\draw[color=black] (8.746222752481538,5.20959718105949) node {$a_{k,0}$};
\draw [fill=black] (8.5,4.) circle (1.0pt);
\draw [fill=black] (9.,5.) circle (1.5pt);
\draw [fill=black] (8.748951172574003,3.502097654851993) circle (1.0pt);
\draw [fill=black] (9.5,5.) circle (1.5pt);
\draw [fill=black] (3.,3.) circle (1.0pt);
\draw [fill=black] (4.,5.) circle (1.5pt);
\draw[color=black] (4.259187384908164,5.30959718105949) node {$a_{i,\infty}$};
\draw [fill=black] (7.,5.) circle (1.5pt);
\draw [fill=black] (4.5,0.) circle (1.0pt);
\end{scriptsize}
\end{tikzpicture}

\end{minipage}
\end{center}
	
	Then, we have:
	$$\val(d-a_{i,\infty})\geq \min\lbrace \val(d-a_{i,0}),\val(a_{i,0}-a_{i,\infty})\rbrace \geq \gamma$$
	Indeed, $\val(a_{i,0}-a_{i,\infty}) \geq \min\lbrace \val(a_{i,0}-a_{k,0}),\val(a_{k,0}-a_{i,\infty})\rbrace = \epsilon \geq \gamma$. 
	Hence, we have for every $0 \leq k < M$:
	$$\val(d-a_{i,\infty}) \geq \val(d-a_{k,0}) = \gamma . $$
	
	\textbf{Case B}: For all $0\leq i,  k < M$ with $i \neq k$, $\left(\val(a_{i,j}-a_{k,l})\right)_{j,l}$ is not constant. \\
	We prove the following claim:
	\begin{claim}
	The indiscernible sequences are totally linearly order by the relation $(a_{k,l})_l \leq (a_{i,l})_l$ defined as follow:
	\[(a_{k,l})_l \preceq (a_{i,l})_l \quad \Leftrightarrow \quad (a_{k,l})_{l<\omega} {\Rightarrow} a_{i,0} \text{  or }(a_{k,-l})_{l<\omega} {\Rightarrow} a_{i,0}.\] 
	\end{claim}

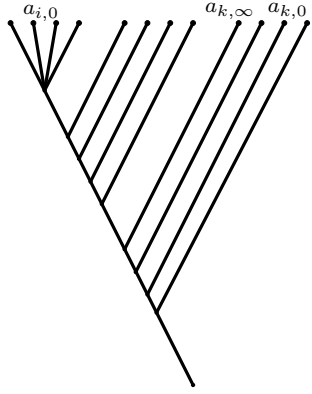
\begin{center}
\begin{minipage}{0.70\linewidth}

	\begin{tikzpicture}[line cap=round,line join=round,>=triangle 45,x=1.0cm,y=1.0cm,scale=0.6]
\clip(-2.6813490653752179,-3.244744554454015) rectangle (16.301104990983102,6.8306249304654555);
\draw [line width=1.2pt] (6.,-3.)-- (2.,5.);
\draw [line width=1.2pt] (2.745113232048317,3.509773535903366)-- (3.5,5.);
\draw [line width=1.2pt] (3.2572757384796804,2.485448523040639)-- (4.5,5.);
\draw [line width=1.2pt] (3.5,2.)-- (5.,5.);
\draw [line width=1.2pt] (3.7524098570521747,1.4951802858956507)-- (5.5,5.);
\draw [line width=1.2pt] (4.,1.)-- (6.,5.);
\draw [line width=1.2pt] (7.,5.)-- (4.5,0.);
\draw [line width=1.2pt] (7.5,5.)-- (4.74424472152211,-0.5033601209440769);
\draw [line width=1.2pt] (8.,5.)-- (5.,-1.);
\draw [line width=1.2pt] (8.5,5.)-- (5.2,-1.4);
\draw [line width=1.2pt] (2.5,5.)-- (2.745113232048317,3.509773535903366);
\draw [line width=1.2pt] (3.,5.)-- (2.745113232048317,3.509773535903366);
\begin{scriptsize}
\draw [fill=black] (6.,-3.) circle (1.0pt);
\draw [fill=black] (2.,5.) circle (1.5pt);
\draw [fill=black] (2.5,5.) circle (1.5pt);
\draw[color=black] (2.671546559378185,5.2231835461590554) node {$a_{i,0}$};
\draw [fill=black] (3.,5.) circle (1.5pt);
\draw [fill=black] (2.745113232048317,3.509773535903366) circle (1.0pt);
\draw [fill=black] (3.5,5.) circle (1.5pt);
\draw [fill=black] (3.2572757384796804,2.485448523040639) circle (1.0pt);
\draw [fill=black] (4.5,5.) circle (1.5pt);
\draw [fill=black] (3.5,2.) circle (1.0pt);
\draw [fill=black] (5.,5.) circle (1.5pt);
\draw [fill=black] (3.7524098570521747,1.4951802858956507) circle (1.0pt);
\draw [fill=black] (5.5,5.) circle (1.5pt);
\draw [fill=black] (4.,1.) circle (1.0pt);
\draw [fill=black] (6.,5.) circle (1.5pt);
\draw [fill=black] (7.,5.) circle (1.5pt);
\draw[color=black] (6.81473270324649,5.27561091638033) node {$a_{k,\infty}$};
\draw [fill=black] (4.5,0.) circle (1.0pt);
\draw [fill=black] (7.5,5.) circle (1.5pt);
\draw [fill=black] (4.74424472152211,-0.5033601209440769) circle (1.0pt);
\draw [fill=black] (8.,5.) circle (1.5pt);
\draw[color=black] (8.07896881089334,5.27561091638033) node {$a_{k,0}$};
\draw [fill=black] (5.,-1.) circle (1.0pt);
\draw [fill=black] (8.5,5.) circle (1.5pt);
\draw [fill=black] (5.2,-1.4) circle (1.0pt);
\end{scriptsize}
\end{tikzpicture}

\end{minipage}
\end{center}

\begin{proof}[Proof of the claim]
   We show first that 
   if two mutually indiscernible sequences $(a_i)_{i\in \bar{\mathbb{Z}}}$ and $(b_i)_{i\in \bar{\mathbb{Z}}}$ such that $(\val(a_i-b_j))_{i,j}$ is not constant, then one of the following occurs:
   \begin{itemize}
       \item $(a_i)_{i<\omega} \Rightarrow b_0$,
       \item $(a_{-i})_{i<\omega} \Rightarrow b_0$,
       \item $(b_i)_{i<\omega} \Rightarrow a_0$,
       \item $(b_{-i})_{i<\omega} \Rightarrow a_0$.
   \end{itemize}
   Assume not of the above cases hold, then $(\val(a_i-b_0))_{i\in \bar{\mathbb{Z}}}$ and  $(\val(b_i-a_0))_{i\in \bar{\mathbb{Z}}}$ are constant. It follows by mutual indiscernibility that $(\val(a_i-b_j))_{i,j}$, contradiction.
   We can conclude by showing the following: given two mutually indiscernible sequences $(a_i)_{i\in \bar{\mathbb{Z}}}$, $(b_i)_{i\in \bar{\mathbb{Z}}}$  and $c_0$, 
   
   \[ \text{If }(a_i)_i \Rightarrow b_0 \text{ and }(b_i)_i \Rightarrow c_0 \text{ then } (a_i)_i \Rightarrow c_0.\]
   By mutual indiscernability, for all $i\in \bar{\mathbb{Z}}$,
   $\val(a_i-b_0)=\val(a_i-a_{i+1}) = \val(a_i-b_j)$ for all $j\in \bar{\mathbb{Z}}$.    Since for all $j\in \bar{\mathbb{Z}}$ \[\val(a_i-c_0) \geq  \min\{\val(a_i-b_{j}), \val(b_{j}-c_0)\}, \]
   and $\val(b_{j}-c_0)$ is increasing with $j$, we must have $\val(a_i-c_0) = \val(a_i-b_0)$. This shows $(a_i)_i \Rightarrow c_0$ and conclude our proof.
\end{proof}

If needed, one can flip the indices and assume that for all $k\neq i$,  either $(a_{k,l})_{l < \omega} {\Rightarrow} a_{i,0}$ or $(a_{i,l})_{l < \omega} {\Rightarrow} a_{k,0}$

We can distinguish two cases:
 \begin{enumerate}
     \item[(B.1)] there is $i \in M$ such that for $\vert M \vert -1 $ many $k\in M\setminus \{i\}$, $(a_{k,l})_l \Rightarrow a_{i,0}$,
     \item[(B.2)] there is $i \in M$ such that for $\vert M \vert $ many $k\in M\setminus \{i\}$, $(a_{i,l})_l \Rightarrow a_{k,0}$. 
 \end{enumerate}
 Note that in case (B.1) only $(a_{i,j})_j$ could be a fan. Case (B.2) is possible only if $M$ is an infinite cardinal.

	Assume that $(B.1)$ holds. For all $k$ such that $(a_{k,l})_l {\Rightarrow} a_{i,0}$, we have \[\val(a_{k,0}-a_{i,\infty})=\val(a_{k,0}-a_{i,0}) \geq \gamma,\]
	since $(a_{k,l}) {\Rightarrow} a_{i,\infty}$ as well.
It remains to prove the inequality for $k=i$. Take  $j\in  M \setminus \{i\}$ such that $(a_{j,l}) {\Rightarrow} a_{i,0}$. We have:
$$\val(a_{i,0}-a_{i,\infty}) \geq \min\lbrace \val(a_{i,0}-a_{j,0}),\val(a_{j,0}-a_{i,\infty}) \rbrace  \geq \gamma.$$
	Assume that $(B.2)$ holds.  For all $k\in M $ such that $(a_{i,l}) {\Rightarrow} a_{k,0}$ we have:
	\[\val(a_{k,0}-a_{i,\infty})\geq \val(a_{k,0}-a_{i,0}) \geq \gamma.\]
It remains to prove the inequality for $k=i$. Take $j \neq i$, such that $(a_{i,l}) \Rightarrow a_{j,l}$. We have:
$$\val(a_{i,0}-a_{i,\infty}) \geq \min\lbrace \val(a_{i,0}-a_{j,0}),\val(a_{j,0}-a_{i,\infty}) \rbrace  \geq \gamma.$$
\end{proof}

\subsection{Quantifier elimination relative to the base field $K$ and the value set $\Gamma$ (variable base field)}

In this subsection, we work with valued vector spaces with variable base field. Similarly to the case of (non-valued) vector spaces, we define $\mathrm{L}_{W,\lambda}$ to be the language of valued vector spaces $\mathrm{L}$ enriched with symbols for $n$-ary predicates $W_n$ and functions $\lambda_{n,i}(x_v,\bar{x_V}): V^{n+1} \rightarrow V$, for $1\leq i\leq n \in \mathbb{N}$:
      \[\mathrm{L}:= \{(V,+,0),(K,+,\cdot, 0, 1),(\Gamma,<,\infty), \scal : K \times V\rightarrow V, \val:V \rightarrow \Gamma\},\]
    \[\mathrm{L}_{W,\lambda}:= \mathrm{L} \cup \lbrace W_n \subset V^n\}_{n\in \mathbb{N}}\cup \{\lambda_{n,i}:V^{n+1} \rightarrow K\rbrace_{1\leq i \leq n \in \mathbb{N}}.\]
    
    We interpret $W_n$ as the set of \textit{$K$-linearly value-independent} $n$-tuples of vectors $(v_1,\dots,v_n)$:
    \[W_n(\bar{v})\equiv \left(\val(v_1)=\cdots=\val(v_{n}) \right) \wedge \left( \forall k_1,\dots,k_{n} \in K \ \val(\sum_{i} k_iv_i)>\val(v_1)  \rightarrow k_1=\cdots=k_{n}=0 \right),\]
    and $\lambda_{n,i}: V^{n+1}\rightarrow K$ as the coordinate projections ( that we also called \textit{lambda-functions}): \\
    for $v,v_1,\cdots,v_n \in V$, if $W_n(\bar{v}) \wedge \neg W_{n+1}(v,\bar{v})$, then $(\lambda_{n,i}(v,\bar{v}))_{1\leq i \leq n}$ is the unique tuples of elements in $K$ such that: 
    \[\val(\sum_{i} \lambda_{n,i}(v,\bar{v}) v_i - v) > \val(v).\]
    
    Notice that by definition, the lambda-functions $\lambda_{n,i}$'s and the predicates $W_n$'s factor to $S(V)$: for $0<i\leq n \in \mathbb{N}$, we may define functions $\lambda_{\bullet,n,i}: S(V)^{n+1}\rightarrow S(V)$ and predicates $W_{\bullet,n} \subset S(V)^n$ such that
    \[\forall v,v_1,\cdots,v_n \in V \ \lambda_{n,i}(v,v_1,\dots,v_n) = \lambda_{\bullet,n,i}(\pi(v),\pi(v_1),\dots,\pi(v_n)),\]
         and
    \[\forall v_1,\cdots,v_n \in V \ W_n(v_1,\dots,v_n) \Leftrightarrow W_{\bullet,n}\left(\pi(v_1),\dots,\pi(v_n)\right).\]
    
    Of course, the valuation $\val:V \rightarrow \Gamma$ also factorises into $\val_S \circ \pi$. Finally, remark that neither the scalar multiplication nor the addition factorise to $S(V)$.

    \begin{remark}\label{RemarkSkeletonEliminationSymbolScal}
        The skeleton equipped with its induced structure (see Corollary \ref{CorollaryInducedStructureSkeleton})
        \[\{(S(V),+_\bullet,\mathbf{0}), (K,+,\cdot, 0,1), (\Gamma,<), \val_S, \scal_{\bullet}\}\]
        is an enriched lexicographic product in the sense of Theorem \ref{EquationEliminationEQuantifierEnrichedLexPro}, that is: it is bi-interpretable with the structure
        \[\{S(V),(\Gamma,<), (K,+,\cdot, 0,1), \val_S:S(V) \rightarrow \Gamma, \lambda_{\bullet,n,i}:S(V)^{n+1}  \rightarrow K, W_{\bullet,n})\}.\]
    \end{remark}     
    \begin{proof}
         The main point of this remark is to show that the function symbols $\scal_{\bullet}$ and $+_\bullet$ can be expressed using the lambda-maps $\lambda_{\bullet,n,i}$ and the predicates $W_{\bullet,n}$.
         Consider an atomic formula $\phi(\bar{\mathbf{x}})$ in the first structure. Let us write 
         \[\bigvee_{\Delta\text{ complete $s$-diagram of }\bar{\mathbf{x}}}\phi \wedge \Delta(\bar{\mathbf{x}}).\]
         As in Remark \ref{RemarkBiinterpretabilityLWLambda}, one sees that each disjunct $\phi \wedge \Delta_s(\mathbf{x})$ is equivalent to a formula $\psi_\Delta$ using $\lambda_{\bullet,n,i}$ and $W_{\bullet,n}$ and with no occurrence of terms of the form $\scal(k,\mathbf{x}_i)$ or $\mathbf{x}_i+_\bullet\mathbf{x}_j$, where $k\in K$, $\mathbf{x}_i,\mathbf{x}_j \in \bar{\mathbf{x}}$. 
    \end{proof}

Since each $K$-vector spaces $B_{\gamma}$, $\gamma \in \Gamma$, eliminates quantifiers relative to $K$, we can show that the skeleton eliminates quantifiers relative to $K$ and $\Gamma$: 

\begin{thmx}\label{TheoremQuantifierEliminationRelativeVariableBaseField}
    As in Example \ref{ExampleEnrichedLexicographicProductVectorSpaces}, consider the theory $T_0$ of lexicographic sums $\mathcal{S}= \{\Gamma[B_\gamma]_\gamma, +_\bullet, (W_{\bullet,n})_n\}$ of $K$-vector spaces $(B_\gamma,0,+,\{W_n\}_{n \in \mathbb{N}})$ in the one sorted language $\mathrm{L}_{W}$  with respect to a linear order $(\Gamma,<)$. We enrich the theory $T_0$ in a theory $T$ as follows:
    \begin{itemize}
        \item with predicates in the sort $(\Gamma,<)$ for the definable sets
    $D_n:=\{\gamma \in \Gamma \ \vert \ \dim_K(B_{\gamma})=n\}$,
        \item with functions $\lambda_{\bullet,i,n}: \Gamma[B_\gamma]_\gamma \times (\Gamma[B_\gamma]_\gamma)^n \rightarrow K$ which to $(n+1)$-uples of the form\\ $(\gamma,\sum_{i}k_iv_i),(\gamma,v_1),\dots, (\gamma,v_{n})$ associates the scalar $k_i$, and $0$ for every other $(n+1)$-uples.
    \end{itemize}
    Then the theory $T$ eliminates quantifiers relative to $K$ and $\Gamma$.

    \begin{center}
                \begin{tikzpicture}
                    \node {$(\Gamma[B_\gamma]_{\gamma}, W_{\bullet,n})$ }
                        child { 
                        node { $(\Gamma,D_n,<)$ } 
                        }
                        child [missing]
                        child { node { $(K,+,\cdot,0,1) $ }
                        };
                \end{tikzpicture}
            \end{center}

\end{thmx}
This is to say that every formula $\psi(\mathbf{x}_S,x_\Gamma,x_K)$ with variables in $\mathbf{x}_S,x_\Gamma,x_K \in S^{\vert \mathbf{x}_S\vert}\times \Gamma^{\vert x_\Gamma\vert} \times K^{\vert x_K\vert}$ is equivalent to a disjunction of formulas of the form:
\[\psi_{\Gamma}(\val_S(\mathbf{x}_S), x_\Gamma) \wedge \psi_K(\lambda(\mathbf{x}_S),x_K) \wedge W(\mathbf{x}_S) ,\]
where $W(\mathbf{x}_S)$ is a conjunction of $W_{\bullet,n}(\mathbf{x}_S')$ for subtuples $\mathbf{x}_S' \subseteq \mathbf{x}_S$ and their negation;  $\psi_\Gamma(y_\Gamma,x_\Gamma)$ is a formula in the language of $(\Gamma,D_n,<)$; $\psi_K(y_K,x_K)$ is a formula in the language of $(K,+,\cdot,0,1)$, and $\lambda$ is a $\vert y_K \vert$-tuple of lambda functions. 

\begin{proof}
    Since the theory of a two-sorted vector spaces is determined by its dimension and the theory of the base field, we get that for all formulas $\phi$ in the language $\mathrm{L}$ of two-sorted vector spaces, the set $\{\gamma \in \Gamma \ \vert \ B_\gamma \models \phi \}$ is definable in $(\Gamma,D_n,<)$.
    Then, the theorem follows  from Theorems  \ref{TheoremRelativeQE} and  \ref{TheoremRelativeQuantifierEliminationEnrichedLexProd}. To the reader who want to avoid the formalism of Section 2,  we give a sketch of a direct proof independent of Theorem \ref{TheoremRelativeQuantifierEliminationEnrichedLexProd}. We proceed by back and forth. Consider two models $\mathcal{M}$ and $\mathcal{N}$ with is partial isomorphism $f=(f_S,f_\Gamma,f_K): A \rightarrow B$ between finitely generated substructure $A=(A_S,A_\Gamma,A_K)$ and $B=(B_S,B_\Gamma,B_K)$. We assume that $f_\Gamma$ and $f_K$ are elementary. 
    Since $A$ is a substructure, $\lambda_{\bullet,i,n}(A_S^{n+1}) \subset A_K$ and $\val(A_S) \subset A_\Gamma$.
    We may extends $f_\Gamma$ and $f_K$ to isomorphims $f_\Gamma': \Gamma^\mathcal{M} \rightarrow f_\Gamma'(\Gamma^\mathcal{M})$ and  $f_K': K^\mathcal{M} \rightarrow f_K'(K^\mathcal{M})$. 
    Then $f\cup f_K' \cup f_\Gamma' $ is a partial isomorphism and we may reset the notation and assume that $A_K=K^\mathcal{M}$ and that $A_\Gamma = \Gamma^\mathcal{M}$.
    Consider now $x \in S_M \setminus A$. \\
    Recall that if $x_1,\cdots,x_{n}$ are elements of the skeleton, $W_{\bullet,n}(x_1,\cdots,x_{n})$ means that they have the same valuation and the vector space components are independent.\\
    Case $1$: there are $x_1,\dots,x_{n} \in A_S \cap \val_S^{-1}(\val_S(x))$ satisfying  $W_{\bullet,n}(x_1,\dots,x_{n})$ and $\neg W_{\bullet,n+1}(x,x_1,\cdots,x_{n})$. Then let $y\in \mathcal{N}$ such that $\val_S(y)=\val_S(f(x_1))$ and  $\lambda_{\bullet,i,n}(y,f_K(x_1),\cdots,f_K(x_{n}))= f_S(\lambda_{\bullet,i,n}(x,x_1,\cdots,x_{n}))$. We can see that $f\cup (x,y)$ is a partial isomorphism.
    \\
    Case $2$: there are no such $x_0,\dots,x_{n-1}$ in  $A_S \cap \val_S^{-1}(\val_S(x))$. Take a set $\{x_0,\dots,x_{n-1}\} \subset A_S\cap \val_S^{-1}(\val_S(x))$ maximal such that $W_{\bullet,n}(x_0,\cdots,x_{n-1})$. It is indeed finite: $A$ is finitely generated and so $A_S$ must be finite. Notice that it can be empty. Then $\Gamma^\mathcal{M} \models D_{n+1}(\val_S(x))$, and so $\Gamma^\mathcal{M} \models D_{n+1}(f_\Gamma^\mathcal{M}(\val_S(x)))$. Take any $y\in N_S$ such that $W_{\bullet,n+1}(y,f_S(x_1),\dots,f_S(x_{n}))$. Then, $f\cup (x,y)$ is a partial isomorphism.
\end{proof}

 \begin{thmx}[Relative quantifier elimination]\label{TheoremQuantifierEliminationRelativeVariableBaseField2}
        The theory of valued vector spaces eliminates quantifiers relative to the value set $(\Gamma,D_n,<)$ , where $D_n:=\{\gamma \in \Gamma \ \vert \ \dim_K(B_\gamma)=n\}$, and to the base field $(K,+,\cdot,0,1)$ in the language $\mathrm{L}_{W,\lambda} \cup \{D_n \subset \Gamma\}_{n\in \mathbb{N}}$.
    \end{thmx}
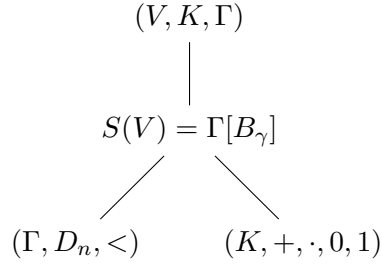
        \begin{center}

                    \begin{tikzpicture}
                    \node {$(V,K, \Gamma)$ }
                        child{ 
                        node {$S(V)=\Gamma[B_\gamma]$}
                            child { 
                            node { $(\Gamma,D_n,<)$ } 
                            }
                            child[missing]
                            child { node { $(K,+,\cdot,0,1) $ }
                            }
                        };
                \end{tikzpicture}
    \end{center}

    The proof only consists of gathering the previous results (Theorems \ref{TheoremQuantifierEliminationRelativetoSkeleton} and \ref{TheoremQuantifierEliminationRelativeVariableBaseField}), and show that since at each knot of the `reduction diagram', we have quantifier elimination relative to children, we can deduce that the structure at the root eliminates quantifiers relative to the structures at the leaves (in an appropriate language). We try to clarify in the following proof. 
    
    \begin{proof}
        Consider a formula $\phi(x_V,x_\Gamma,x_K)$ in $\mathrm{L}_{W,\lambda}$ with variables $x_V,x_\Gamma,x_K \in V^{\vert x_V\vert}\times \Gamma^{\vert x_\Gamma\vert} \times K^{\vert x_K\vert}$.
        Recall that $\mathrm{L}_S$ denotes the language for the skeleton. 
 \[\mathrm{L}_S:=\{(S(V),+_\bullet,\mathbf{0}), \val_S,\scal_{\bullet}\},\]
        and  $\pi: V \rightarrow S(V)$ denotes the projection map. 
               We add to the language $\mathrm{L}_{W,\lambda}$ a sort for the skeleton: 
        \[\mathrm{L}_{S,W,\lambda} : = \{S(V),+_\bullet,\mathbf{0}\}\cup \{K,+,\cdot,0,1)\} \cup \{\Gamma,D_n,<\} \cup \{\val_S,\lambda_{\bullet,n,i},W_{\bullet,n}\}.\]
        Notice that we also removed $\lambda_{n,i}$ and $W_n$.
        By replacing each occurrence of terms of the form $\lambda_{n,i}(t,t_1,\dots, t_n)$ (where the $t_i$'s are vector space terms) by $\lambda_{\bullet,n,i}(\pi(t),\pi(t_1),\dots, \pi(t_n))$, and each occurrence of atomic formula of the form $W_{n}(t_1,\dots, t_n)$ by $W_{\bullet,n}(\pi(t_1),\dots, \pi(t_n))$, we see that our formula $\phi(x_V,x_\Gamma,x_K)$ is equivalent to a formula \[\psi(\pi(t(x_V,x_K)),x_\Gamma,x_K)\]
        where $t(x_V,x_K)$ is a tuples of vector space terms and 
        $\psi(\mathbf{x}_S,x_V,x_K)$ is an $\mathrm{L}_{W,\lambda}\cup \mathrm{L}_{S,W,\lambda}$-formula. 
        
        Now, notice that $\mathrm{L}_{S,W,\lambda}$ is formally an enrichment of $\mathrm{L}_{S}$. By Theorem \ref{TheoremQuantifierEliminationRelativetoSkeleton} and resplendence,  $\phi(x_V,x_\Gamma,x_K)$ is equivalent to a formula of the form
        \[\psi_S(\pi(t(x_V,x_K)),x_\Gamma,x_K)\]
        where $t(x_V,x_K)$ is a tuple of vector space term and  $\psi_S(\mathbf{x}_S,x_\Gamma,x_K)$ is a formula in the language          
        \[\left\{(S(V),+_\bullet,\mathbf{0}), (\Gamma, <), (K,+,\cdot, 0,1), \val_S, \scal_{\bullet}, W_{\bullet,n}, \lambda_{\bullet,n,i}, \val_S \right\}.\]
        
        By Remark \ref{RemarkSkeletonEliminationSymbolScal} $\psi_S(\mathbf{x}_S,x_\Gamma,x_K)$ is equivalent to a formula without any occurrence of $\scal_{\bullet}$ and $+_\bullet$. Now we can use Theorem \ref{TheoremQuantifierEliminationRelativeVariableBaseField}: $\psi_S(\mathbf{x}_S,x_\Gamma,x_K)$ is equivalent to a disjunction of formulas of the form
        \[\psi_{\Gamma}(\val_S(\mathbf{x}_S), x_\Gamma) \wedge \psi_K(\lambda(\mathbf{x}_S),x_K) \wedge W_S(\mathbf{x}_S) ,\]
where $W_S(\mathbf{x}_S)$ is a conjunction of $W_{\bullet,n}(\mathbf{x}_S')$ for subtuples $\mathbf{x}_S' \subseteq \mathbf{x}_S$ and their negation;  $\psi_\Gamma(y_\Gamma,x_\Gamma)$ is a formula in the language of $(\Gamma,D_n,<)$; $\psi_K(y_K,x_K)$ is a formula in the language of $(K,+,\cdot,0,1)$, and $\lambda$ is a $\vert y_K \vert$-tuple of lambda functions $\lambda_{\bullet,n,i}$.
    
    Putting all these results together and using that $\val=\val_S\circ \pi$, $W_n=W_{n,s}\circ \pi$ and $\lambda_{n,i}=\lambda_{\bullet,n,i}\circ \pi$ , we get that our original formula $\phi(x_V,x_\Gamma,x_K)$ is equivalent to a disjunction of formulas of the form
    \[\phi_{\Gamma}(\val(t(x_V,x_K)), x_\Gamma) \wedge \phi_K(\lambda(t(x_V,x_K)),x_K) \wedge W(t(x_V,x_K))\]
    
where $W(y_V)$ is a conjunction of $W_{n}(y_V')$ for subtuples $y_V' \subseteq y_V$ and their negation;  $\phi_\Gamma(y_\Gamma,x_\Gamma)$ is a formula in the language of $(\Gamma,D_n,<)$; $\phi_K(y_K,x_K)$ is a formula in the language of $(K,+,\cdot,0,1)$ with $\lambda$ a $\vert y_K \vert$-tuple of lambda functions $\lambda_{n,i}$; and $t(x_V,x_K)$ is a $\vert y_V \vert$ tuple of vector space terms.  This concludes our proof.
    \end{proof}
    
From Theorem \ref{TheoremQuantifierEliminationRelativeVariableBaseField2}, we deduce the following:
\begin{corollary}
    Let $V$ be a valued vector space with variable base field $K$ and value set $(\Gamma,<)$. 
    \begin{itemize}
        \item The induced structure on the value set $\Gamma$ is given by  $(\Gamma,D_n,<)$ where $D_n:=\{\gamma \in \Gamma \ \vert \ \dim_K(B_\gamma)=n\}.$
        \item The induced structure on $K$ is the pure structure of field  $(K,+,\cdot, 0,1)$.
        \item The induced structure on ${\gamma}\times B_\gamma$ is the pure structure of (two-sorted)  $K$-vector space.
    \end{itemize}
    
\end{corollary}

    \subsection{Transfer principles relative to the base field $K$ and the value set $\Gamma$ (variable base field)}

    As usual, let us state first relative model completeness:
    
    \begin{proposition}
        (Relative model completness) Assume $(V,\Gamma,K)$ is a valued vector space with variable base field $K$ and value set $(\Gamma,<)$ in the language $\mathrm{L}_{W,\lambda}$. Consider an extension $(W,\Delta,L) \supseteq (V,\Gamma,K)$. Then the following are equivalent:
            \begin{itemize}
                \item $(W,\Delta,L) \succeq (V,\Gamma,K)$,
                \item $L\succeq K$ as fields and $(\Delta,(E_n)_{n\in \mathbb{N}},<) \succeq (\Gamma,(D_n)_{n\in \mathbb{N}},<)$ as colored linear order where $E_n:=\{\gamma \in \Delta \ \vert \ \dim_L(B_\gamma(W))=n\}$ and $ \ D_n:=\{\gamma \in \Gamma \ \vert \ \dim_K(B_\gamma(V))=n\}$.
            \end{itemize}
    \end{proposition}
    
   The proof is deduced from Theorem \ref{TheoremQuantifierEliminationRelativeVariableBaseField2} and Remark \ref{RemarkModelCompletnessRQECLoseSort}.

    \begin{proposition}
        ($\NIP$ transfer) A valued vector space $\mathcal{V}$ with variable base field $K$ and value set $(\Gamma,<)$ is $\NIP$ if and only if $K$ is NIP (and $\Gamma$ is $\NIP$).
    \end{proposition}
    
    As we notice earlier, equipped with its full induced structure, $(\Gamma,<)$ is a linear order with (countably many) colors $D_n$. By \cite[Proposition A.2]{Sim15}, such a structure is dp-minimal or in other words, NIP and of burden $1$. We keep mentioning the value set in the propositions above and below, as these results hold for any enrichment of $K$ and $\Gamma$.

    \begin{proof}
    Assume $K$ (and $\Gamma$) to be NIP, we must show that $\mathcal{V}$ is NIP. By Proposition \ref{PropositionNIPTransferSkeleton}, it is enough to show that $S(V)$ is NIP. We see the Skeleton in the following language:
            \[\{ S,  (K,+,\cdot,0,1),(\Gamma, D_n, <), \lambda_{\bullet,n,i}, W_{\bullet,n}, \val_S\},\]
    Then it follows immediately from Corollary \ref{CorollaryNIPTransferLexicoProduct} (enriched version) and Corollary \ref{CorollaryNIPtransferVectorSpace}. For a proof independant of Section 2, we can argue as follow: we have to show that any atomic formula $\phi(x;y_S,y_K,y_\Gamma)$ with variables $y_K\in K, \ y_\Gamma \in \Gamma , \ x,y_S\in S(V)$ is $\NIP$.
    These are of the form:
        \begin{itemize}
            \item $W_{\bullet,n}(x,y_S)$ , 
            \item $\phi_\Gamma(\val_S(x),y_\Gamma)$ where $\phi_\Gamma(x_\Gamma,y_\Gamma)$ is a $\Gamma$-sorted formula, 
            \item $\phi_K(\lambda(x,y_\Gamma),y_K)$, where $\phi_K(x_K,y_K)$ is a $K$-sorted formula and $\lambda(x,y_\Gamma)$ is a tuple of $\lambda_{\bullet,n,i}$-terms.
        \end{itemize}
        
        As for vector spaces, the first formula is $\NIP$ of VC-dimension at most $n$.  The second and the third are $\NIP$ since $\phi_\Gamma(x_\Gamma,y_\Gamma)$ and $\phi_K(x_K,y_K)$ are $\NIP$.
    \end{proof}

    \begin{proposition}\label{PropositionValuedVectorSpaceBurdenTranfer}
        Let $(V, K, \Gamma, \val)$  be a valued vector space with base field $K$ and value set $\Gamma$. We also equipped $\Gamma$ with its induced structure, i.e. with predicates for $D_n:=\{ \gamma \in \Gamma \ \vert \ \dim_K(B_{\gamma})=n\}$ for $n\in \mathbb{N}$.
        
        Then  $\bdn(V)= \sup_n(\kappa^{n}_{inp}(K),\bdn(\Gamma))$ if the $K$-dimension of the $B_\gamma$'s are not bounded by any integer, that is to say if $\Gamma \neq \bigcup_{n\leq d} D_n.$ for all $d$.

        Otherwise, $\bdn(\mathcal{V}) = \max \{\kappa^{d}_{inp}(K), \bdn(\Gamma)\}$ where  $d$ is the smallest integer such that $\Gamma= \bigcup_{n\leq d} D_n $ and  
    \end{proposition}
    We deduce it from Theorem \ref{TheoremBurdenVectorSpace}, Proposition \ref{PropositionBdnValuedVectorSpaceRelativeSkeleton} and (the enriched version of) Theorem \ref{ThmBdnLexicoProd}.
    \begin{proof} 
     By Proposition \ref{PropositionBdnValuedVectorSpaceRelativeSkeleton}, we only have to compute the burden of the skeleton. 
     
    By definition, the ribs $\{ (B_\gamma,+,0), (K,+,\cdot,0,1), \scal \}$ equipped with their structure of two-sorted vector spaces interpret  uniformly the base field. Recall that since the theory of a two-sorted vector space is determined by its dimension and the theory of the base field, we get that for all formula $\phi$ in the language $\mathrm{L}$ of two-sorted vector spaces, $\{\gamma \in \Gamma \ \vert \ B_\gamma \models \phi \}$ is definable in $(\Gamma,D_n,<)$.
     
    By Corollary \ref{CorollaryInducedStructureSkeleton}, the skeleton equipped with its induced structure can be seen as the enriched lexicographic sum $\Gamma[B_\gamma]_\gamma$ of the ribs $\{ (B_\gamma,+,0), (K,+,\cdot,0,1), \scal \}$ with respect to the spine $(\Gamma,D_n,<)$.  A rib in an elementary extension is equivalent to $\{(B_\rho,+,0),(K,+,\cdot,0,1), \scal\}$ where $B_\rho=\prod_\rho B_\gamma$ and $\rho$ is a one type in $(\Gamma,D_n,<)$. By Theorem \ref{TheoremBurdenVectorSpace}, the burden of $B_\rho$ depend of its $K$-dimension, which depend only on $\rho$. We get:
        \[ \rho(x) \vdash x\in D_d \wedge x\notin D_{d+1} \quad \Leftrightarrow \quad  \bdn(B_\rho) = \kappa^{d}_{inp},\]
     and
     \[ \rho(x) \vdash \{x\notin D_{n}\}_{n\in \mathbb{N}} \quad \Leftrightarrow \quad  \bdn(B_\rho) = sup_n\{\kappa^{n}_{inp}\}.\]
     
    By the enriched version of Theorem \ref{ThmBdnLexicoProd}, 
    we get the result. Notice once again that by \cite[Proposition A.2]{Sim15}, $\bdn(\Gamma,D_n,<)$ is of burden $1$.
    
    \end{proof}
    
    With Adler's convention, we see that we have 
    \[ \bdn^\star(\mathcal{V})\leq  \max({\kappa^{d}_{inp}}^\star(K),\bdn^\star(\Gamma)),\]
    if ribs have finite K-dimension at most $d$ and
    \[ \bdn^\star(\mathcal{V})=  {\sup_n}^\star({\kappa^{n}_{inp}}^\star(K),\bdn^\star(\Gamma)),\]
    if the K-dimension of the ribs is not bounded by a finite number.

    \section*{Acknowledgements}
    I would like to thank the logic team of the university of Campania Luigi Vanvitelli, and in particular Lorna Gregory for all the useful discussions and Paola D'Aquino for her support and guidance.
    I am also very gratefull to the Oberwolfach institute for providing me a research environment and financial support when I started this project and to the university of Campania Luigi Vanvitelli for giving me the opportunity to continue this work.
    
    \bibliographystyle{plain}
\bibliography{bibtex}

\end{document}